\documentclass[12pt]{article}
\usepackage{amsmath}
\usepackage{amssymb}
\usepackage{amscd}
\usepackage{epsfig}
\usepackage{amsthm}

\def\t{\tilde}

\newcommand{\length}{\operatorname{\ul{length}_{\Delta}}}
\newcommand{\Length}{\operatorname{length_{\Delta}}}

\def\ora{\overrightarrow}

\newtheorem{dfn}{Definition}[section]
\newtheorem{rem}[dfn]{Remark}
\newtheorem{thm}[dfn]{Theorem}
\newtheorem{defn}[dfn]{Definition}
\newtheorem{lem}[dfn]{Lemma}
\newtheorem{lemma}[dfn]{Lemma}

\newtheorem{prop}[dfn]{Proposition}

\newtheorem{cor}[dfn]{Corollary}

\newtheorem{conjecture}[dfn]{Conjecture}

\newtheorem{conv}[dfn]{Convention}
\newtheorem{ex}[dfn]{Example}

\newtheorem{question}[dfn]{Question}

\def\proof{\par\medskip\noindent{\it Proof: }}

\def\lra{\longrightarrow}

\def\ra{\rightarrow}

\def\CR{\curvearrowright}
\def\acts{\CR}

\def\embed{\hookrightarrow}

\def\P{{\mathbb P}}
\def\C{{\mathbb C}}
\def\R{{\mathbb R}}

\def\Z{{\mathbb Z}}
\def\K{{\mathbb K}}

\def\O{{\mathcal O}}

\def\Q{{\mathbb Q}}

\def\N{{\mathbb N}}

\def\eps{\epsilon}
\def\al{\alpha}
\def\be{\beta}
\def\ga{\gamma}

\def\Ga{\Gamma}

\def\De{\Delta}

\def\Del{\Delta}
\def\Si{\Sigma}
\def\si{\sigma}

\def\la{\lambda}
\def\La{\Lambda}

\def\>{\rangle}
\def\<{\langle}
\def\del{\delta}

\def\3{\ss}
\def\8{\infty}

\def\ol{\overline}
\def\ul{\underline}
\def\ov{\overrightarrow}
\newcommand{\restr}{\mbox{\Large \(|\)\normalsize}}

\begin{document}

\title{A path model for geodesics in Euclidean buildings and its
applications to  representation theory}
\author{Misha Kapovich and John J. Millson}
\date{April 12, 2008}
\maketitle

\begin{abstract}
In this paper we give a combinatorial characterization of projections of geodesics in Euclidean buildings to
Weyl chambers. We apply these results to the representation theory of
complex reductive Lie groups and to spherical Hecke rings associated with split nonarchimedean reductive
Lie groups. Our main application is a generalization of the
saturation theorem of Knutson and Tao for $SL_n$ to other complex semisimple
Lie groups.
\end{abstract}

\section{Introduction}

Let $\ul{G}$  be a $\Q$-split reductive algebraic group defined over
$\Z$ and $\ul{G}^\vee$ be its Langlands' dual. In this paper we
continue our study (which we began in \cite{KLM3}) of the interaction
between the representation theory of the group
$G^\vee:=\ul{G}^\vee(\C)$ and geometry of the Bruhat--Tits building
associated with the nonarchimedean group
$G=\ul{G}(\K)$, where $\K$ is a complete field with discrete
valuation. We restrict ourselves to the case when $\K$ is a local field, in which case, algebraically speaking,
we will be studying the relation between
the representation ring of the group
$G^\vee$ and the {\em spherical Hecke algebra}
$\mathcal{H}_G$ associated with  $G$.

In his papers \cite{Littelmann1}, \cite{Littelmann2},
P.~Littelmann introduced a {\em path model} for the
representations  of complex reductive Lie groups $G^\vee$. The
Littelmann path model gives a method to compute the structure
constants of the representation ring of $G^\vee$ by counting
certain piecewise linear paths, called LS paths.

In this paper we define a class of piecewise-linear paths in $\De$, a Weyl chamber of the Weyl group $W$ of $G^\vee$.
These paths will be called {\em Hecke paths}  (see Definition \ref{Hecke}) because of their
connection with $\mathcal{H}_G$. We will prove that a path $p$ in $\De$ is a Hecke path if and only if $p$
is the projection into $\De$ of a geodesic segment in the
Euclidean (Bruhat-Tits) building $X$ associated with $G$.  Thus, unlike LS paths which had to be {\em invented}, the
Hecke paths appear very naturally as projections of geodesic
segments. Hecke paths are defined by eliminating one of the axioms for
LS paths, therefore {\em each LS path for $G^\vee$ is a Hecke path
for $G$}.

The converse relation is more subtle and is discussed later in the
introduction. To state our main results we need a definition of
$k_R$, the {\em saturation factor} of the root system $R$ of the
group $\ul{G}$. Let $\alpha_1,...,\alpha_l\in R$ be the simple
roots (corresponding to $\De$). Let $\theta$ be the highest root
and define positive integers $m_1,...,m_l$ by
$$
\theta = \sum_{i=1}^l m_i \alpha_i.$$
Then $k_R$ is the least common multiple of the numbers $m_i$, $i=1,...,l$.
We refer to section \ref{kfactors} for the computation of $k_R$.

\medskip
Below, $L$ is the character lattice of a maximal torus in $G^\vee$ (so that $\De\subset L\otimes \R$),
$Q(R^\vee)$ is the root lattice of $R^\vee$. Our main result is the following
theorem (see section \ref{saturationsection}), which, in a weaker
form, has been conjectured by S.~Kumar:

\begin{thm}
[Saturation theorem]\label{main}  Let $G^{\vee}$ and $L$ be as above. Suppose that
$\al,\be,\ga\in L$ are dominant characters such that $\al+\be+\ga\in
Q(R^\vee)$ and that there exists $N\in \N$ so that
$$
(V_{N\al}\otimes V_{N\be} \otimes V_{N\ga})^{G^\vee}\ne 0.
$$
Then for $k=k_R^2$ we have:
$$
(V_{k\al}\otimes V_{k\be} \otimes V_{k\ga})^{G^\vee}\ne 0.
$$
\end{thm}

Here and in what follows $V_\la$ is the irreducible representation
of $G^\vee$ associated with the dominant weight $\la$ of $G^\vee$.
Also it will be convenient to introduce integers $n_{\al,\be}(\ga^*)$,
the structure constants of the representation ring of the group $G^\vee$. (Here and in what follows, $\ga^*$
is the dominant weight contragredient to $\ga$.) Hence
$n_{\al,\be}(\ga^*)$ are defined by the equation
$$
V_{\alpha} \otimes V_{\be} = \bigoplus_{\ga} n_{\al,\be}(\ga^*) V_{\ga^*}.$$
Here the right-hand side is the decomposition of the tensor product
of the irreducible representations $V_{\alpha}$ and $V_{\be}$
into a direct sum of irreducible representations.
We can then formulate the above theorem as
$$n_{N\alpha, N\beta}(N \gamma^*)\ne 0 \Rightarrow
n_{k\alpha, k\beta}(k \gamma^*)\ne 0.
%\ \text{for all dominant weights} \ \alpha,\beta,\gamma.
$$

We will freely move back and forth between the symmetric formulation of
saturation in Theorem \ref{main} and the asymmetric formulation immediately
above.

As an immediate corollary of Theorem \ref{main} we obtain a new
proof of the saturation theorem of A.~Knutson and T.~Tao
\cite{KnutsonTao}:

\begin{cor}
Suppose that $R=A_l$, i.e. the semisimple part of $G^\vee$ is
locally isomorphic to $SL_{l+1}$. Suppose that $\al,\be,\ga$ are
dominant characters such that $\al+\be+\ga\in Q(R^\vee)$ and that
there exists $N\in \N$ so that
$$
(V_{N\al}\otimes V_{N\be} \otimes V_{N\ga})^{G^\vee}\ne 0.
$$
Then
$$
(V_{\al}\otimes V_{\be} \otimes V_{\ga})^{G^\vee}\ne 0.
$$
\end{cor}

Another proof of this theorem was given by H.\ Derksen and J.\ Weyman in \cite{DW}. However
the proofs of \cite{KnutsonTao} and \cite{DW} do not work for root systems different from $A_l$.

\begin{question}
\label{simplyconjecture}
Is it true that  if $G$ is a simple simply-laced group, then in Theorem \ref{main}
one can always take $k=1$ and in the case of non-simply laced groups the
smallest $k$ which suffices is $k=2$?
\end{question}

The affirmative answer to this question is supported by the odd orthogonal groups, symplectic groups and $G_2$,
when one gets the saturation constant $k=2$ rather than $2^2$ and $6^2$ \cite{KLM3, BK},
the group $Spin(8)$, when the saturation constant equals $1$
\cite{KKM}, as well as by a number of computer
experiments with the exceptional root systems and the root systems $D_l$.

Using the results of \cite{KLM2, KLM3} one can reformulate  Theorem
\ref{main} as  follows:

\begin{thm}
\label{reform} There exists a convex homogeneous cone $D_3\subset \De^3$,
defined by the {\em generalized triangle inequalities},
which depends only on the Weyl group $W$, so that the following hold:

1. If a triple $(\al,\be,\ga)\in (\De\cap L)^3$ satisfies
$$
(V_{\al}\otimes V_{\be} \otimes V_{\ga})^{G^\vee}\ne 0,
$$
then $(\al,\be,\ga)\in D_3$ and
$$
\al+\be+\ga\in Q(R^\vee).
$$

2. ``Conversely'', if $(\al,\be,\ga)\in k^2_R \cdot L^3 \cap D_3$ and
$\al+\be+\ga\in k^2_R \cdot Q(R^\vee)$, then
$$
(V_{\al}\otimes V_{\be} \otimes V_{\ga})^{G^\vee}\ne 0.
$$
\end{thm}

We now outline the steps required to prove the Theorem \ref{main}. We will need several facts about Hecke rings
${\mathcal H}$. Let $\mathcal{O}$ denote the valuation ring of $\K$, then
$K:=\ul{G}(\mathcal{O})$ is a maximal compact subgroup in $G$.
The lattice $L$ defined above is the cocharacter lattice of a maximal torus $T\subset G$.
The Hecke ring ${\mathcal H}$, as a $\Z$-module, is freely generated by the characteristic functions
$\{c_\la: \la \in L \cap \Delta\}$. The
multiplication on ${\mathcal H}$ is defined via the convolution
product $\star$. Then the  structure constants $m_{\al,\be}(\ga)$ of ${\mathcal H}$
are defined by
$$
c_\al \star c_\be= \sum_\ga m_{\al,\be}(\ga) c_\ga.
$$
We refer the reader to  \cite{Gross},
\cite{KLM3} and to section \ref{heck} of this paper for more details.

Let $o\in X$ be the special vertex of the building $X$ which is
fixed by $K$. In section \ref{distances} we define the notion
$d_\De(x,y)$ of the $\De$-valued distance between points $x, y\in
X$. Given a piecewise-geodesic path $p$ in $X$ we define its
$\De$-length as the sum of $\De$-distances between the consecutive
vertices.

The structure constants for ${\mathcal H}$ are related to the geometry of $X$ via the following

\begin{thm}
\cite[Theorem 8.12]{KLM3}.
\label{Hecke=triangles}
The number $m_{\al,\be}(\ga)$ is equal to the product
of a certain positive constant  %(depending only on $G$)
by the  number of geodesic triangles $T\subset X$ whose vertices are {\em special vertices}
of $X$ with the first vertex equal to $o$ and whose $\De$-side
lengths are $\al, \be, \ga^*$.
\end{thm}

Theorem \ref{main} is essentially Statement 3, which follows from
Statements 1 and 2, of the following

\begin{thm}
\label{3->4}
Set $\ell:=k_R$.

1. Suppose that
$(\al,\be,\ga)\in D_3\cap L^3$ and $\al+\be+\ga\in Q(R^\vee)$.
Then the structure constants $m_{\cdot,\cdot}(\cdot)$ of the Hecke ring of the group $G$ satisfy
$$
m_{\ell\al,\ell\be}(\ell\ga)\ne 0.$$

2. Suppose that $\al, \be, \ga$ are dominant coweights of $\ul{G}$, such that
$m_{\al,\be}(\ga)\ne 0$. Then

$$ n_{\ell\al,\ell\be}(\ell\ga)\ne 0.$$

3. As a consequence of 1 and 2 we have: Suppose that
$(\al,\be,\ga)\in D_3\cap L^3$ and $\al+\be+\ga\in Q(R^\vee)$.
Then
$$
n_{\ell^2\al,\ell^2\be}(\ell^2\ga)\ne 0.$$
\end{thm}

\begin{rem}
a. Part 1 of the above theorem was proven in \cite{KLM3}. Thus the point
of this paper is to prove Part 2 of the above theorem.

b. Examples in \cite{KLM3} show that both implications
$$
(\al,\be,\ga)\in D_3\cap L^3, \al+\be+\ga\in Q(R^\vee) \Rightarrow
m_{\al,\be}(\ga)\ne 0
$$
and
$$
m_{\al,\be}(\ga)\ne 0 \Rightarrow n_{\al,\be}(\ga)\ne 0
$$
are {\bf false} for the groups $G_2$ and $SO(5)$. Therefore the dilation by $k_R$ in both cases is necessary at least
for these groups.
\end{rem}

More generally, we prove (Theorem \ref{mini3->4})

\begin{thm}
\label{minus}
Suppose that $\al, \be, \ga\in L$ are dominant weights so that one of them
is the sum of minuscule weights. Then
$$
m_{\alpha, \beta}(\gamma)\ne 0 \Rightarrow
n_{\alpha, \beta}(\gamma)\ne 0.
$$
\end{thm}

The proof of Part 2 of Theorem \ref{3->4} proceeds as follows.
In section \ref{characterizationsection} we prove a characterization
theorem for folded triangles which implies:

\begin{thm}
\label{char}
 There exists a geodesic triangle $T\subset X$ whose
vertices are {\em special vertices} of $X$ and whose $\De$-side
lengths are $\al, \be, \ga^*$ if and only if there exists a Hecke
path $p: [0,1]\to \De$ of $\De$-length $\be$ so that
$$
p(0)=\al, p(1)=\ga.
$$
\end{thm}

\begin{figure}[tbh]
\centerline{\epsfxsize=4.5in \epsfbox{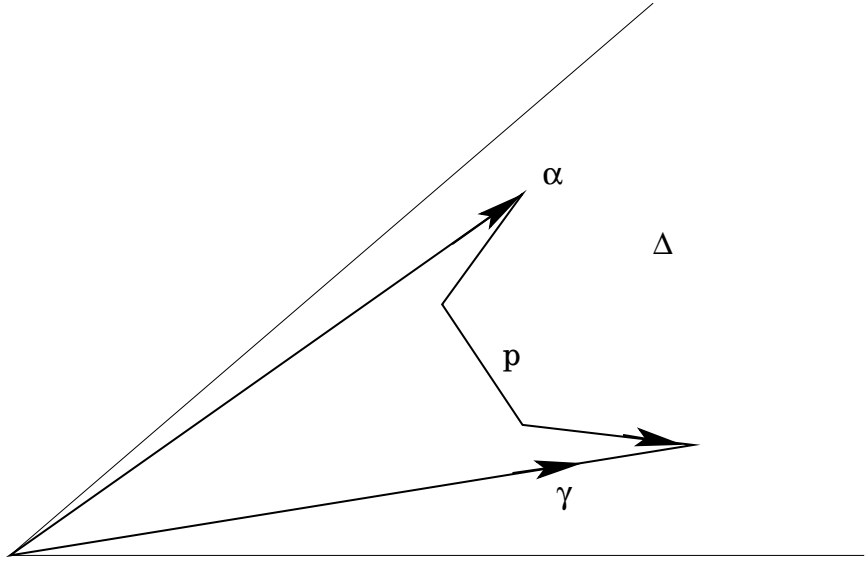}}
\caption{\sl A broken triangle.}
\label{triangle.fig}
\end{figure}

The directed segments $\pi_\alpha$, the Hecke path $p$ and the (reversed) directed
segment $\pi_\gamma$ fit together to form a ``broken triangle'',
see Figure \ref{triangle.fig}. Here and in what follows $\pi_\la$ is the geodesic path
parameterizing the directed segment $\ov{ox}=\la$.

Then, by combining theorems \ref{Hecke=triangles} and \ref{char}, we
obtain

\begin{thm}
\label{comp}
$m_{\al,\be}(\ga)\ne 0$ if and only if there exists a Hecke path
$p: [0,1]\to \De$ of $\De$-length $\be$  so that
$$
p(0)=\al, p(1)=\ga.
$$
\end{thm}

This statement is an analogue of Littelmann's theorem which relates structure constants $n_{\al,\be}(\ga)$
of the representation ring with LS paths. The problem however is that not every Hecke path is an LS
path (even for the group $SL(3)$).

In order to prove Theorem \ref{main} our ``path model'' for the nonvanishing of
the Hecke structure constants must be generalized to a model
where Hecke paths are replaced by {\em generalized Hecke
paths} (we also have to replace the LS-paths by {\em generalized LS paths}). More precisely,
begin with a geodesic triangle $[x,y,z]\subset X$ whose vertices are special vertices of $X$ and whose $\De$-side lengths
are $\be, \ga^*, \al$. Now replace the geodesic segment $\ol{xy}$
with a certain piecewise-geodesic path $\t{p}\subset X$
connecting $x$ and $y$ in $X$, and
which is contained in the 1-skeleton of a single apartment in $X$. We then show that the projection $p$ of $\t{p}$
to $\De$ is a {\em generalized Hecke path}. The $\De$-length of the path $p$ still equals $d_\De(x,y)=\be$.

The generalized Hecke paths have the advantage over Hecke
paths that  {\em their break points occur only at vertices of the building}.
By using this observation we obtain

\begin{thm}\label{dilate}
If $p$ is a generalized Hecke path with $p(0)=0$, then $k_R \cdot p$, the image
of $p$ under dilation by $k_R$, is a generalized LS path of $\De$-length equal to $k_R \Length(p)$.
\end{thm}

Lastly, for {\em generalized LS paths} we prove, by modifying slightly Littelmann's arguments, the following:

\begin{thm}
\label{Li}
$n_{\al,\be}(\ga)\ne 0$ if and only if there exists a generalized LS path $q$ with $\De$-length $\be$ so that
$$
\al+q(1)=\ga
$$
and the concatenation $\pi_\al*q$ is contained in $\De$.
\end{thm}

Theorem \ref{main} is now obtained in 2 steps, the first of which is contained in \cite{KLM3} and the second
is at the heart of the present paper:

Step 1. It was shown in \cite[Theorem 9.17]{KLM3} that
$$
n_{N\al,N\be}(N\ga)\ne 0 \Rightarrow m_{N\al,N\be}(N\ga)\ne 0.
$$
Therefore, the vector $(N\al, N\be, N\ga^*)$ belongs to the homogeneous cone $D_3$.
Hence, by Part 1 of Theorem \ref{3->4}
we conclude that $m_{\ell\al,\ell\be}(\ell\ga)\ne 0$.

Step 2. $m_{\ell\al,\ell\be}(\ell\ga)\ne 0$ implies existence of a geodesic triangle
$[z, x, y]\subset X$ with the special vertices and $\De$-side lengths $\ell\al, \ell\be, \ell\ga^*$.
Then by projecting the corresponding path $\t{p}$ to $\De$ and dilating it by $k_R$ we obtain
a path $k\al+q(t)$ in $\De$ connecting $k\al$ to $k\ga^*$, where $q(t)$ is a generalized LS path of
$\De$-length $k\be$. (Recall that $k=\ell^2$.) Therefore, by appealing to Theorem \ref{Li}, we see that
$$
n_{k\al,k\be}(k\ga)\ne 0
$$
which concludes the proof.

\medskip
This paper is organized as follows. Preliminary material is
discussed in section \ref{prelim}, where we review the concepts of
Coxeter complexes, buildings and piecewise-linear paths in buildings, as well as
the {\em generalized distances} in buildings. In section
\ref{sectionchains} we define the notion of {\em chains}, which is
essentially due to Littelmann. This concept allows one to define
both LS paths and Hecke paths, as well as to relate Hecke paths and
foldings of geodesic paths, which is discussed in the next section.

The main technical tool of this paper is the concept of {\em
folding} of geodesics in a building into an apartment, which is done
via retraction of the entire building to an apartment or chamber.
\footnote{A similar idea is used by S. Gaussent and P. Littelmann in
\cite{GL}, where they fold galleries in a building to galleries in
an apartment.} Properties of foldings are discussed in section
\ref{sectionfolding}. In section \ref{euclideanfolding} we prove
that the image of each geodesic segment in $X$ under {\em folding}
$f: X\to \De$ of $X$ to a Weyl chamber, is a Hecke path, Theorem \ref{chaincondition}.
We then prove a partial converse to this result, i.e. that each
Hecke path which satisfies the {\em simple chain condition} can be
unfolded in $X$. We also find a necessary and sufficient condition
for  unfolding of a path $p$ which is {\em local}, i.e. it depends
only on germs of the path $p$ at its break-points.

In section \ref{LSpolygons} we review Littelmann's
path model for the representation theory of complex semisimple Lie
groups, in particular we discuss LS paths and generalized LS paths
as well as raising and lowering operators.

We use {\em approximation of LS paths by paths satisfying the simple
chain condition} to unfold LS paths in $X$, Theorem \ref{LSfolded}
in section \ref{unfoldingLS}. Although there are Hecke paths which
are not LS paths, since the unfolding condition is local, by restricting
the root system we reduce the general unfolding problem to the case
of the LS paths. We thus establish that a path in $\De$ is
unfoldable in $X$ it and only if it is a Hecke path, this is done in
section \ref{characterizationsection}, Theorem
\ref{chain-unfolding}. The reader interested only in the proof of
the saturation theorem can omit this section.

We prove the saturation theorem in section \ref{saturationsection},
Corollary \ref{corQ3->Q4}. As we stated above, the idea of the proof is to replace Hecke
paths with piecewise-linear paths contained in the 1-skeleton of the Euclidean
Coxeter complex. We show that dilation by $k_R$ of such a path
results in a {\em generalized LS path} which in turn suffices for
finding nonzero invariant vectors in triple tensor products.

\paragraph{Hecke paths and folded galleries.} It is interesting to ask what the relationship is between our paper
using Hecke paths and the results of S. Gaussent and P.
Littelmann in \cite{GL} and others using positively folded galleries:

(a) {\em  Hecke paths} correspond to the {\em positively folded galleries} defined in \cite{GL}.
However this correspondence can be established only {\bf after} the folding--unfolding results of the present
paper are proven.
Therefore it appears that
one cannot prove our results characterizing Hecke paths as projections of geodesic segments
in building using the results of \cite{GL} and vice versa.

(b) In \cite{Schwer}, C. Schwer has used the gallery approach to
compute the Hecke structure constants $m_{\alpha,\beta}(\gamma)$.
In an earlier version of this paper we applied our theory to
compute some Hecke structure coefficients. Thus  both Hecke path
and gallery models can be used  to compute the structure constants
of the spherical Hecke algebra.

(c) One can prove that
$$n_{\al,\be}(\ga)\ne 0\Rightarrow m_{\al,\be}(\ga)\ne 0$$
\cite[Theorem 9.17]{KLM3} using both Hecke paths (as it is done in section \ref{characterizationsection})
and positively folded galleries.

(d) There does not  seem to be a direct way to carry over our
proof of the main Theorem \ref{main} or part 2 of Theorem
\ref{3->4} to a proof using positively folded galleries. It is
clear from the above, that one of the critical steps in our
argument is provided by Theorem \ref{dilate} stating that dilation
by $k_R$ converts generalized Hecke paths to generalized LS paths.
But there does not appear to be a way to ``stretch'' a folded
gallery by the factor $k_R$ --- or, more generally, to produce an
LS gallery from a positively folded gallery. This is an advantage
of the Hecke paths, which are {\em geometric} objects, over
positively folded galleries, which are {\em combinatorial}
objects. Even when $k_R=1$, there is no obvious (at least to us)
reason why existence of a positively folded gallery would imply
existence of an LS gallery. We note that if $k_R \neq 1$ then the
existence of a positively folded gallery does not imply the
existence of an LS gallery, there are counterexamples for $SO(5)$
and $G_2$.

\medskip
{\bf Acknowledgments.} This paper grew out of our joint work
\cite{KLM1, KLM2, KLM3} with Bernhard Leeb. We are grateful to him
for all his ideas that he contributed to those papers, especially
the idea of folding that came out of his contribution to
\cite{KLM3}. We are also grateful to Tom Haines for valuable discussions. We thank the referee for his/her remarks.
Part of the work on this paper was done when the first author was
visiting the Max Plank Institute (Bonn), during this work he was also supported by
the NSF grant DMS-02-03045 and DMS-04-05180. The second author was
supported by the NSF grant DMS-01-04006 and DMS-04-05606. Together, the authors were suported by the NSF grant
DMS-05-54349. The authors gratefully acknowledge support of these institutions.

\tableofcontents

\section{Definition and notation}
\label{prelim}

\subsection{Root systems and Coxeter complexes}

A (discrete, nonnegatively curved) {\em Coxeter complex} is a pair
$(A,W)$, where $A$ is either a Euclidean space or the unit sphere
and $W$ is a discrete reflection group acting on $A$. The {\em rank} of
the Coxeter complex is the dimension of $A$. The group $W$ is called {\em Weyl group} of the Coxeter complex.
It is called an {\em affine Weyl group} if $A$ is a Euclidean space.

%\begin{rem}
%There are also Euclidean Coxeter complexes which are nondiscrete (i.e. the
%group $W$ is nondiscrete, but has discrete linear part), as well as hyperbolic Coxeter complexes,
%where one uses reflection groups of isometries of the hyperbolic
%space. However we do not need them in the present paper.
%\end{rem}

An isomorphism of Coxeter complexes $(A,W), (A',W')$ is an isometry $\iota: A\to A'$ so that
$$
\iota W \iota^{-1}=W'.
$$

{\em Walls} in
$(A,W)$ are fixed point sets of reflections $\tau\in W$. A point
$x\in A$ is called {\em regular} if it does not belong to any wall
and {\em singular} otherwise. The closure of each connected
component of the set of regular points is called an {\em alcove} in
the Euclidean case and a {\em chamber} in the spherical case.

In the case when $W$ acts cocompactly on $A$, each alcove (resp. chamber) is a product (resp. join) of simplices
and $(A,W)$ determines structure of a {\em polysimplicial complex} on $A$. In general
there exists a totally-geodesic subspace
$A'\subset A$ which is $W$--invariant and such that $A'/W$ is compact.
Therefore each alcove in $A$ is a product
of simplices and a Euclidean subspace in $A$. Thus much of the discussion
 of Euclidean Coxeter complexes can be reduced to the case when $A/W$ is compact.

\begin{rem}
\label{tri}
A  triangulation of a fundamental alcove (chamber) in $A$ determines a
$W$--invariant simplicial complex on $A$. Thus
we can always think of $A$ as a simplicial complex.
\end{rem}

A {\em half-apartment} in $A$ is the closure of a
connected component of $A\setminus H$, where $H$ is a wall in $A$.

``Most'' Coxeter complexes are associated with root systems as we
describe below. Suppose that $R$ is a root system on a vector space
$V$ (i.e. each element $\al\in R$ is a linear functional on $V$).
The {\em rank} of $R$ is the number of simple roots in $R$, i.e. is the rank of the free abelian
subgroup in $V^*$ generated by $R$.
Let $A$ denote the Euclidean affine space corresponding to $V$. This
data defines a finite Coxeter group $W_{sph}$, which is a reflection
group generated by reflections in the hyperplanes $H_{\al}=\{ x:
\al(x)=0\}$. {\em Weyl chambers} of $W_{sph}$ are closures of the
connected components of the complement to
$$
\cup_{\al\in R} H_\al.
$$
In what follows we fix a {\em positive Weyl chamber} $\De$, it
determines the subset of positive roots  $R^+\subset R$ and of {\em
simple roots} $\Phi\subset R^+$. We also have the group of coweights
$P(R^\vee)$ associated with $R$:
$$
\la\in P(R^\vee) \iff \forall \al\in R, \quad \al(\la)\in \Z.
$$
Let $W_{aff}$ denote the {\em affine Coxeter group} determined by
the above data, this group is generated by reflections in the
hyperplanes ({\em affine walls})
$$
H_{\al,t}=\{\al(x)=  t\}, t\in \Z.$$

Given a vector $\ga\in \De$, we define the {\em contragredient vector} $\ga^*$
as $w_0(-\ga)$, where $w_0$ is the longest element of $W$.
In other words, $\ga^*$ is the intersection of the $W$-orbit
$W\cdot (-\ga)$ with the chamber $\De$.

The translation subgroup in $W_{aff}$ is the {\em coroot lattice}
$Q(R^\vee)$,  it is generated by the coroots $\al^\vee, \al\in R$.
The group  $P(R^\vee)$ is the normalizer of $W_{aff}$ in the group
of translations $V$. Note that $V/P(R^\vee)$ is compact. The
dimension of this quotient is the same as the dimension of $V$
provided that $rank(R)$ equals the dimension of $V$.

 The {\em special vertices} of a Euclidean Coxeter complex are the
points whose stabilizer in $W_{aff}$ is isomorphic to $W_{sph}$.
Equivalently, they are the points in the $P(R^\vee)$-orbit of the
origin.

\begin{rem}
If $A/W_{aff}$ is compact, the  vertices of $(A,W_{aff})$ are the
vertices of the polysimplicial complex determined by $W_{aff}$.
\end{rem}

Let $A^{(0)}$ denote the {\em vertex
set} of $(A, W_{aff})$, which consists of points of maximal
intersection of walls in $A$. If $R$ spans $V^*$, the set $A^{(0)}$
equals the vertex set of the polysimplicial complex in $A$ defined
by tessellation of $A$ via alcoves of $W_{aff}$.

Given a Coxeter complex $(A,W)$ and a point $x\in A$ we define a new
Coxeter complex $(S_x,W_x)$ where $S_x$ is the unit tangent sphere
at $x$ and $W_x$ is the stabilizer of $x$ in $W$.

For a nonzero vector $\nu\in V$ we let $\bar{\nu}:= \nu/|\nu|$
denote the {\em normalization} of $\nu$. We define {\em rational}
elements of the unit sphere $S$ to be the unit vectors of the form
$$
\eta=\bar\nu, \quad \nu\in P(R^\vee).
$$
The next lemma follows immediately from compactness of
$V/P(R^\vee)$:

\begin{lem}
\label{density} Rational points are dense in $S$.
\end{lem}

Suppose that $(A, W)$ is a Euclidean Coxeter complex. A {\em
dilation} of $(A, W)$ is a dilation $h$ (i.e. a composition of
translation and similarity $v\mapsto \la v, \la>0$) in the affine
space $A$ so that
$$
h W h^{-1} \subset W.
$$
We let $Dil(A, W)$ denote the semigroup of dilations of the complex
$(A, W)$. We will refer to the number $\la$ as the {\em conformal
factor} of the dilation $h$.

Given a point $x\in A$ and a dilation $h\in Dil(A,W)$, we can define
a new spherical Coxeter complex $(S_x, W'_x)$ on the unit tangent
sphere $S_x$ at $x$ via pull-back
$$
W'_x:= h^*(W_{t(x)}),
$$
where $W_{h(x)}$ is the stabilizer of $h(x)$ in $W$.

\begin{defn}
\label{relation} Suppose that $W$ is a finite Coxeter group acting
on a vector space $V$. Define a (nontransitive) relation $\sim_W$ on
$V\setminus \{0\}$ by
$$
\mu \sim_W \nu \iff
$$
\centerline{$\mu, \nu$ belong to the same Weyl chamber of $W$.}

\noindent We will frequently omit the subscript $W$ in this
notation.
\end{defn}

\begin{defn}
\cite[page 514]{Littelmann2}
\label{porder} We say that
nonzero vectors $\nu, \mu\in V$ satisfy $\nu \rhd_W \mu$ (for short,
$\nu \rhd \mu$) if for each positive root $\al$,
$$
\al(\nu)< 0 \Rightarrow \al(\mu)\leq 0.
$$
\end{defn}

\begin{lem}
\label{difference}
Suppose that $\nu, \mu\in P(R^\vee)$, $w\in W=W_{aff}$ is such that $w(\nu)=\mu$. Then
$$
\mu-\nu\in Q(R^\vee).
$$
\end{lem}
\proof The mapping $w$ is a composition of reflections $\tau_i\in W$.
Therefore it suffices to prove the assertion in case when $w$
is a reflection $\tau$. This reflection is a composition of
a translation $t$  and a reflection $\si\in W_o$.
The translation $t$ belongs to the translation subgroup $Q(R^\vee)$ of $W$,
therefore it suffices to consider the case when $\tau=\si\in W_o$.
Then $\tau=\tau_\be$, where $\be$ is a root and we have
$$
\mu-\nu= -\be(\nu)\be^\vee.
$$
Since $\be(\nu)\in \Z$ and $\be^\vee\in Q(R^\vee)$, the assertion follows. \qed

\subsection{Paths}
\label{paths}

Suppose that $A, V, W_{aff}$, etc., are as in the previous section.

Let $\tilde{\mathcal P}$ denote the set of all piecewise-linear paths $p: [a,b]\to
V$.  We will be identifying paths that differ by
orientation-preserving re-parameterizations $[a,b]\to [a',b']$.
Accordingly, we will always (re)parameterize a piecewise-linear path with
constant speed. We let $p'_-(t), p'_+(t)$ denote the derivatives of the function $p$
from the left and from the right. The space $\t{\mathcal P}$ will be
given the topology of uniform convergence.

If $p, q: [0,1]\to A$ are piecewise-linear paths in a simplicial complex such that
$p(1)=q(0)$, we define their {\em composition} $r=p\cup q$ by
$$
r(t)=\left\{
\begin{array}{c}
p(t), ~~t\in [0,1],\\ q(t-1), ~~t\in [1, 2].
 \end{array}\right.
$$
Let ${\mathcal P}\subset \tilde{\mathcal P}$ denote the set of paths
$p: [0,1]\to V$ such that $p(0)=0$. Given a path $p\in {\mathcal P}$
we let $p^*\in {\mathcal P}$ denote the {\em reverse} path
$$
p^*(t)= p(1-t)-p(1).
$$
For a vector $\la\in V$ define a
geodesic path $\pi_\la\in {\mathcal P}$ by
$$
\pi_\la(t)=t\la, \quad t\in [0,1].
$$
Given two paths $p_1, p_2\in {\mathcal P}$ define their {\em
concatenation} $p=p_1* p_2$ by
$$
p(t)= \left\{
\begin{array}{c}
p_1(2t), \quad t\in [0, 1/2],\\ p_1(1)+ p_2(2t-1), \quad t\in [1/2,
1].
 \end{array}\right.
$$

\medskip
Suppose that $p\in {\mathcal P}$ and $J=[a,b]$ is nondegenerate
subinterval in $I=[0,1]$. We will use the notation  $p|J\in
\tilde{\mathcal P}$ to denote the function-theoretic restriction of
$p$ to $[a,b]$. We will use the notation $p|_J$ to denote the path
in ${\mathcal P}$ obtained from $p|J$ by pre-composing $p|J$ with an
increasing linear bijection $\ell: I\to J$ and post-composing it
with the translation by the vector $-p(a)$.

\medskip
Fix a positive Weyl chamber $\De\subset V$; this determines the set of positive roots $R^+\subset R$,
the set of simple roots $\Phi\subset R^+$. We define the
subset ${\mathcal P}^+\subset {\mathcal P}$ consisting of the paths
whose image is contained in $\De$.

For a path $p\in {\mathcal P}$ and a positive root $\al\in R^+$ define the
{\em height} function
$$
h_\al(t)= \al(p(t))
$$
on $[0,1]$. Let $m_\al=m_\al(p)\in \R$ denote the minimum of
$h_\al$. Clearly $m _\al(p)\le 0$ for all $p\in {\mathcal P}$.
 We define the set of ``integral paths''
$$
{\mathcal P}_\Z:= \{ p\in {\mathcal P}: \forall \al\in \Phi,
m_\al(p)\in \Z\}.
$$
More restrictively, we define the set ${\mathcal P}_{\Z,loc}$ of paths $p\in {\mathcal
P}$  which satisfy the following {\em local integrality condition}:

For each simple root $\al\in \Phi$ the function $h_\al$  takes
integer values at the points of local minima.

\subsection{The saturation factors associated to a root system}
\label{kfactors}

In this section we define and compute {\em saturation factors}
associated with root systems.  Let $o\in A$ be a
special vertex, which we will identify with
$0\in V$.

\begin{dfn}
We define the {\em saturation factor} $k_R$ for the root system $R$
to be the least natural number $k$ such that $k\cdot A^{(0)}\subset
P(R^\vee)\cdot o$. The numbers $k_R$ for the irreducible root
systems are listed in the table (\ref{ta}).
\end{dfn}

Note that the condition that $k\cdot A^{(0)}\subset P(R^\vee)\cdot
o$ is equivalent to that each point of $k\cdot A^{(0)}$ is a special
vertex.

Below we explain how to compute the saturation factors $k_R$
following \cite{KLM3}. First of all, it is clear that if the root
system $R$ is reducible and $R_1,...,R_s$ are its irreducible
components, then $k_R= LCM(k_{R_1},...,k_{R_s})$, where  $LCM$
stands for the {\em
  least common multiple}. Henceforth we can assume that the system
$R$ is reduced, irreducible and has rank $n=\dim(V)$. Then the
affine Coxeter group $W_{aff}$ acts cocompactly on $A$ and its
fundamental domain (a {\em Weyl alcove}) is a simplex.

Let $\{\al_1,...,\al_{\ell}\}$ be the collection of simple roots in $R$
(corresponding to the positive Weyl chamber $\De$) and $\theta$ be
the highest root. Then

\begin{equation}
\label{highestroot} \theta= \sum_{i=1}^{\ell} m_i \al_i.
\end{equation}

We have

\begin{lem}
\cite[Section 2]{KLM3}
\label{kexists} $k_R=
LCM(m_1,...,m_n)$.
\end{lem}

Below is the list of saturation factors:

\begin{equation}
\label{ta}
\begin{array}{|c|c|c|}
\hline ~ & ~ & ~  \\ \hbox{Root system} & \theta &  k_R \\ \hline
A_\ell & \al_1+...+\al_\ell & 1\\ \hline B_\ell & \al_1+
2\al_2+...+2\al_\ell & 2\\ \hline C_\ell & 2\al_1+
2\al_2+...+2\al_{\ell-1}+\al_\ell & 2\\ \hline D_\ell & \al_1+
\al_2+\al_3+ 2\al_4+...+2\al_\ell & 2 \\ \hline G_2 & 3\al_1+2\al_2&
6 \\ \hline F_4 & 2\al_1+ 3\al_2+4\al_3 +2\al_4 &  12\\ \hline E_6 &
\al_1+\al_2+ 2\al_3+2\al_4+2\al_5+3\al_6 &  6 \\ \hline E_7 & \al_1+
2\al_2+2\al_3+2\al_4+3\al_5+&  12 \\ ~ & +3\al_6+ 4\al_7 & ~ \\
\hline E_8 & 2\al_1+2\al_2+3\al_3+3\al_4+ 4\al_5+ &  60\\ ~&
+4\al_6+ 5\al_7 +6\al_8 & ~  \\ \hline
\end{array}
\end{equation}

\subsection{Buildings}

Our discussion of buildings follows \cite{KleinerLeeb}. We refer the reader
to \cite{Brown}, \cite{Ronan}, \cite{Rousseau} for the more combinatorial discussion.

Fix a spherical or Euclidean (discrete) Coxeter complex $(A,W)$, where $A$ is a
Euclidean space $E$ or a unit sphere $S$ and $W=W_{aff}$ or $W=W_{sph}$
is a discrete Euclidean or a spherical Coxeter group acting on $A$.

A metric space $Z$ is called {\em geodesic} if any pair of points $x, y$ in $Z$ can be connected by a
geodesic segment $\ol{xy}$.

Let $Z$ be a metric space. A {\em geometric structure} on $Z$ {\em
modeled on} $(A,W)$ consists of an atlas of isometric embeddings
$\varphi:A\embed Z$ satisfying the following compatibility
condition: For any two charts $\varphi_1$ and $\varphi_2$, the
transition map $\varphi_2^{-1}\circ\varphi_1$ is the restriction of
an isometry in $W$. The charts and their images,
$\varphi(A)=a\subset Z$, are called {\em apartments}. We will
sometimes refer to $A$ as the {\em model apartment}. We will require
that there are {\em plenty of apartments} in the sense that any two
points in $Z$ lie in a common apartment. All $W$-invariant notions
introduced for the Coxeter complex $(A, W)$, such as rank, walls, singular
subspaces, chambers etc., carry over to geometries modeled on $(A,
W)$. If $a, a'\subset X$ are alcoves (in the Euclidean case) or chambers
(in the spherical case) then there exists an apartment $A'\subset X$ containing $a\cup a'$:
Just take regular points $x\in a, x'\in a'$ and an apartment $A'$ passing through $x$ and $x'$.

A geodesic metric space $Z$ is said to be a $CAT(\kappa)$-space if
{\em geodesic triangles in $Z$ are ``thinner'' than geodesic
triangles in a simply-connected complete surface of the constant
curvature $\kappa$}. We refer the reader to \cite{Ballmann}  for the
precise definition. Suppose that $Z$ is a (non-geodesic) metric space with the discrete metric:
$$
d(x,y)=\pi \iff x\ne y.
$$
We will regard such a space as a $CAT(1)$ space as well.

\begin{dfn}
A {\em spherical building} is a $CAT(1)$-space modeled on a spherical
Coxeter complex.
\end{dfn}

Spherical buildings have a natural structure of polysimplicial
piecewise spherical complexes. We prefer the geometric to the
combinatorial view point because it appears to be more appropriate
in the context of this paper.

\begin{defn}
A Euclidean building is a $CAT(0)$-space modeled on a (discrete)
Euclidean Coxeter complex.
\end{defn}

A building is called {\em thick} if every wall is an intersection of
apartments. A non-thick building can always be equipped with a
natural structure of a thick building by reducing the Coxeter group  \cite{KleinerLeeb}.

\medskip
Let $\K$ be a local field with a (discrete) valuation $\nu$ and valuation ring $\O$.
Throughout this paper, we will be working with $\Q$-split reductive algebraic groups (i.e. Chevalley groups) $\ul{G}$ over $\Z$.
We refer the reader to \cite{Borel, Demazure} for a detailed discussion of such groups. A reader
can think of $\K=\Q_p$ (the $p$-adic numbers) and a {\em classical} Chevalley group $\ul{G}$,
e.g. $GL(n)$ or $Sp(n)$, in which case  Bruhat-Tits buildings
below can be described rather explicitely, see e.g. \cite{Garrett}.

Given a group $\ul{G}$ as above, we get a {\em nonarchimedean Lie
group} $G=\ul{G}(\K)$, to which we can associate a Euclidean building (a
Bruhat-Tits building) $X=X_G$. We refer the reader to \cite{BT},
\cite{KLM3} and \cite{Rousseau} for more detailed discussion of the
properties of $X$. Here we only recall that:

1. $X$ is thick and locally compact.

2. $X$ is modeled on a Euclidean Coxeter complex $(A,W_{aff})$ whose
dimension equals the rank of $\ul{G}$, and the root system is
isomorphic to the root system of $\ul{G}$.

3. $X$ contains a special vertex $o$ whose stabilizer in $G$ is $\ul{G}(\O)$.

\begin{ex}
Let $X$ be a (discrete) Euclidean building, consider the  {\em
spaces of directions} $\Si_x X$. We will think of this space as the
space of germs of non-constant geodesic segments $\ol{xy}\subset X$.
As a polysimplicial complex $\Si_xX$ is just the link of the point
$x\in X$. The space of directions has the structure of a spherical
building modeled on $(S,W_{sph})$, which is thick if and only if $x$
is a special vertex of $X$  \cite{KleinerLeeb}. The same applies
in the case when $X$ is a spherical building.
\end{ex}

If $X$ is a Euclidean building modeled on $(A,W)$, for each point
$x\in X$ the space of directions $\Si_x(X)$ has {\em two structures}
of a spherical building:

1. The {\em restricted building structure} which is modeled on the
Coxeter complex $(S, W_{x})$, where $S=S_x(A)$ is the unit tangent
sphere at $x$ and $W_{x}$ is the stabilizer of $x$ in the Coxeter
group $W$. This building structure is thick.

2. The {\em unrestricted building structure} which is modeled on the
Coxeter complex $(S, W_{sph})$, where $S=S_x(A)$ is the unit tangent
sphere at $x$ and $W_{sph}$ is the linear part of the affine Coxeter
group $W_{aff}$. This building structure is not thick, unless $x$ is
a special vertex.

Let $B$ be a spherical building modeled on a spherical Coxeter
complex $(S, W_{sph})$. We say that two points $x, y\in B$ are
antipodal, if $d(x,y)=\pi$; equivalently, they are antipodal points
in an apartment $S'\subset B$ containing both $x$ and $y$. The
quotient map $S\ra S/W_{sph}\cong\De_{sph}$ induces a canonical
projection $\theta:B\ra \De_{sph}$ folding the building onto its
model Weyl chamber. The $\theta$-image of a point in $B$ is called
its {\em type}.

\begin{rem}
To define $\theta(x)$ pick an apartment $S'$ containing $x$ and a
chart $\phi: S\to S'$. Then $\theta(x)$ is the projection of
$\phi^{-1}(x)$ to  $S/W_{sph}\cong\De_{sph}$. We note that this is
clearly independent of $S'$ and $\phi$.
\end{rem}

\begin{lem}
\label{thetaproperties}
 1. If $h: A\to A'$ is an isomorphism of apartments in $B$ (i.e.
${\phi'}^{-1}\circ h\circ \phi\in W$) then $\theta\circ h=\theta$.

2. If $x, x'\in B$ which belong to apartments $A, A'$ respectively
and $-x\in A, -x'\in A'$ are antipodal to $x, x'$, then
$\theta(x)=\theta(x')$ implies $\theta(-x)=\theta(-x')$.
\end{lem}
\proof (1) is obvious, so we prove (2). Pick an isomorphism $h: A\to
A'$. Then (since $\theta(x)=\theta(x')$) there exists $w\in W\acts
A'$ such that $w(h(x))=x'$. Hence $w\circ h(-x)=-w\circ(x)=-x'$. The
claim now follows from 1. \qed

\medskip
We will regard {\em $n$-gons} $P$ in a building $X$ as maps $\nu:
\{1,...,n\}\to X$, $\nu(i)=x_i$, where $x_i$ will be the vertices of
$P$.  If $rank(X)\ge 1$ we can connect the consecutive vertices
$x_i, x_{i+1}$ by shortest geodesic segments $\ol{x_i
x_{i+1}}\subset X$ thus creating a {\em geodesic polygon} $[x_1,
x_2,\ldots, x_n]$ in $X$ with the edges $\ol{x_i x_{i+1}}$. Observe
that in case $x_i, x_{i+1}$ are antipodal, the edge $\ol{x_i
x_{i+1}}$ is not unique.

We say that two subsets $F, F'$ in a building $X$ are {\em
congruent} if there exist apartments $A, A'$ in $X$ containing $F,
F'$ resp., and an isomorphism $A\to A'$ of Coxeter complexes which
carries $F$ to $F'$.

\begin{conv}
Suppose that $X$ is a spherical building. We will be considering
only those geodesic triangles $T$ in $X$ for which the length of
each geodesic side of $T$ is $\le \pi$.
\end{conv}

\medskip
Let $X, Y$ be buildings and $f: X\to Y$ a continuous map satisfying
the following: For each alcove (in Euclidean case) or spherical
chamber (in the spherical case) $a\subset X$, the image $f(a)$ is
contained in an apartment of $Y$ and the restriction $f|a$ is
an isometry or similarity. Then we call $f$ {\em differentiable} and
define the {\em derivative} $df$ of $f$ as follows. Given a point
$x\in X$ and $y=f(x)$, the derivative $df_x$ is a map $\Si_x(X)\to
\Si_y(Y)$. For each $\xi\in \Si_x(X)$ let $\ol{xz}\subset X$ be a
geodesic segment whose interior consists of regular points only and
so that $\xi$ is the unit tangent vector to  $\ol{xz}$. Then $f$
sends $\ol{xz}$ to a nondegenerate geodesic segment $\ol{y f(z)}$
contained in an apartment $A\subset Y$. Then we let $df_x(\xi)\in
\Si_y(Y)$ be the unit tangent vector to $\ol{y f(z)}$.

We will be also using the above definition in the setting when a building $Y$ is a Euclidean
Coxeter complex $(A,W)$ and $h\in Dil(A,W)$. Then after letting $X:= h^*(A,W)$ we get an isometry $h: X\to Y$.

\begin{conv}
Throughout the paper we will be mostly using roman letters $x, y,
z$, etc., to denote points in Euclidean buildings and Greek
letters $\xi, \eta, \zeta$, etc., to denote points in spherical
buildings. Sometimes however (e.g., in section \ref{converting})
we will be working simultaneously with a spherical building $X$
and its links $\Si_x(X)$, which are also spherical buildings. In
this case we will use roman letters for points in $X$ and Greek
letters for points  in $\Si_x(X)$.
\end{conv}

\subsection{Generalized distances and lengths in buildings}
\label{distances}

Let $(A,W)$ be a spherical or Euclidean Coxeter complex. The
complete invariant of a pair of points $(x,y)\in A^2$ with respect
to the action $W\acts A$, is its image $d_{ref}(x,y)$ under the
canonical projection to $A\times A/W$. Following \cite{KLM2} we call $d_{ref}(x,y)$ the
{\em refined distance from $x$ to $y$}. This notion carries over to
buildings modeled on the Coxeter complex $(A,W)$: For a pair of
points $(x,y)$ pick an apartment $A'$ containing $x,y$ and, after
identifying $A'$ with the model apartment $A$, let $d_{ref}(x,y)$ be
the projection of this pair to $A\times A/W$.

If points $\xi, \eta$ in a spherical building are antipodal we will
use $\pi$ for the refined distance $d_{ref}(\xi,\eta)$: This does
not create much ambiguity since given apartment contains unique
point antipodal to $\xi$.

In the case of Euclidean Coxeter complexes there is an extra
structure associated with the concept of refined length. Given a
Euclidean Coxeter complex $(A, W_{aff})$, pick a special vertex
$o\in A$. Then we can regard $A$ as a vector space $V$, with the
origin $0=o$. Let $\De\subset A$ denote a Weyl chamber
 of $W_{sph}$, the tip of $\De$ is at $o$.

Then following \cite{KLM2}, we define the {\em $\De$-distance} between points of $(A,
W_{aff})$ by composing $d_{ref}$ with the natural forgetful map
\[ A\times A/W_{aff}\ra A/W_{sph}\cong\De .\]
To compute the $\De$-distance $d_\De (x,y)$ we regard the oriented
geodesic segment $\ol{xy}$ as a vector in $V$ and project it to
$\De$. Again, the concept of $\De$-distance, carries over to the
buildings modeled on $(A,W_{aff})$.

\begin{defn}
\label{maindef} Let $X$ be a thick Euclidean building. Define the
set $D_n(X)\subset \De^n$ of $\De$-side lengths which occur for
geodesic $n$-gons in $X$.
\end{defn}

It is one of the results of \cite{KLM2} that $D_n:=D_n(X)$ is a convex
homogeneous polyhedral cone in $\De^n$, which depends only on $(A,W_{sph})$.
The polyhedron $D_3$ in Theorem \ref{reform} is the polyhedron $D_3(X)$.
The set of {\em stability} inequalities defining $D_n$
is determined in \cite{KLM1} and
\cite{BS}.

\begin{thm}
\label{klm2}\cite[Corollary 8.4]{KLM3}.
 Let $X$ be a thick
Euclidean building modeled on $(A,W_{aff})$. Suppose that $\al, \be,
\ga\in P(R^\vee)$, $\al+\be+\ga\in Q(R^\vee)$ and $(\al,\be,\ga)\in
D_3(X)$. Then there exists a geodesic triangle $T\subset X$ whose
vertices are vertices of $X$ and the $\De$-side lengths are $\al,
\be, \ga$.
\end{thm}

\medskip
Suppose that $p$ is a piecewise-linear path in a Euclidean building $X$. We say
that $p$ is a {\em billiard path} if for each $t, s\in [0,1]$ the
tangent vectors $p'(t), p'(s)$ have the same projection to the
chamber $\De$ in the model apartment. If $p$ is a path which is the
composition
$$
\ol{x_0 x_1}\cup ...\cup \ol{x_{m-1} x_{m}}
$$
of geodesic paths, then the $\De$-length of $p$ is defined as
$$
\Length(p):=\sum_{i=1}^m d_\De (x_{i-1}, x_i)
$$
where $d_\De (x,y)$ is the $\De$-distance from $x$ to $y$

Each piecewise-linear path $p$ admits a unique representation
$$
p=p_1\cup...\cup p_n
$$
as a composition of maximal billiard subpaths so that
$$
\la_i= \Length (p_i).
$$
We define
$$
\length(p):= \ul{\la}= (\la_1,...,\la_n).
$$
Clearly, $\Length(p)$ is the sum of the vector components of
$\ul{\la}$.

\subsection{The Hecke ring}
\label{heck}

In this section we review briefly the definition of spherical
Hecke rings and their relation to the geometry of Euclidean buildings; see
\cite{Gross, KLM3} for more details.

Let $\K$ denote a locally compact field with discrete valuation
$v$, and (necessarily) finite residue field of the order $q$. Let  $\mathcal{O}$
be the subring of $\K$ consisting of elements with nonnegative
valuation. Choose a uniformizer $\pi\in \mathcal{O}$.

Consider a connected reductive algebraic group $\ul{G}$ over $\K$.
We fix a maximal split torus $\ul{T}\subset \ul{G}$ defined over
$\mathcal{O}$. We put $G := \ul{G}(\K), K:= \ul{G}(\mathcal{O})$ and
$T :=\ul{T}(\K)$. We let $\ul{B} \subset \underline{G}$ be a Borel subgroup
normalized by $\underline{T}$ and set $B:= \underline{B}(\K)$.

Let $X$ denote the Bruhat-Tits building
associated with the group $\ul{G}$; $o\in X$ is a distinguished special vertex
stabilized by the compact subgroup $K$.

We also have  free abelian cocharacter group
$X_*(\underline{T})$  of $\ul{T}$ whose rank equals $\dim(\underline{T})$.
This group contains the set of coroots $R^{\vee}$ of the group $\ul{G}$.
The roots are the characters of $\underline{T}$
that occur in the adjoint representation on the Lie algebra of $\underline{G}$.
The subset $R^+$ of the roots that occur in representation on the Lie algebra
of $\underline{B}$ forms a positive system and the indecomposable elements
of that positive system form a system of simple roots $\Phi$. We let $W$ denote
the
corresponding (finite) Weyl group.

The set of positive roots $\Phi$ determines a positive Weyl chamber $\Pi^+$ in
$X_*(\underline{T})$, by
$$
\Pi^+ = \{\lambda \in X_*(\underline{T}): \ \<\lambda,\alpha\> \geq 0, \alpha
\in \Phi\}.
$$
This chamber is a fundamental domain for the action of $W$ on
$X_*(\underline{T})$.

We define a partial ordering on $\Pi^+$ by $\la \succ \mu$ iff the difference
$\la-\mu$
is a sum of positive coroots.

\begin{dfn}
The (spherical) Hecke ring $\mathcal{H} = \mathcal{H}_{G}$ is
the ring of all locally constant, compactly supported functions $f: G\lra \Z$
which are $K$--biinvariant. The multiplication in $\mathcal{H}$ is by the
convolution
$$
f \star g (z) = \int_G f(x) \cdot g(x^{-1}z)dx$$
where $dx$ is the Haar measure on $G$ giving $K$ volume $1$.
\end{dfn}

The ring $\mathcal{H}$ is commutative and associative (see e.g. \cite{Gross}).
For $\lambda \in X_*(\ul{T})$ let
$c_{\lambda}$ be the characteristic function of the corresponding
$K$-double coset $\la(\pi)\in K\backslash G/K$.
Then the functions $c_\la, \la\in \De$ freely generate ${\mathcal H}$ as an
(additive) abelian group. Deep result of Satake \cite{Satake} relates the Hecke ring of $G$ and
the representation ring of $G^\vee$.

The structure constants  for ${\mathcal H}$ are defined by the formula

\begin{equation}
\label{heckesum}
c_{\lambda} \star c_{\mu} = \sum_\nu
m_{\lambda,\mu}(\nu) c_{\nu} =
c_{\lambda + \mu} + \sum
m_{\lambda,\mu}(\nu) c_{\nu},
\end{equation}
where the last sum is taken over all $\nu\in \Pi^+$ such that $\lambda + \mu
\succ \nu$ and therefore is
finite.

It turns out that the structure constants $m_{\lambda,\mu}(\nu)$
are nonnegative integers which are polynomials in $q$ with integer
coefficients. These constants are determined by the geometry of
the building as follows. Given $\al, \be, \ga$ let
$\mathcal{T}=\mathcal{T}_{\al,\be}(\ga)$ denote the (finite) set
of geodesic triangles $[o, x, y]$ in the building $X$ which have
the $\De$-side lengths $\al, \be, \ga^*$, so that $y$ is the
projection of the point
$$
\gamma(\pi)\in T\subset G
$$
into $X$ under the map $g\mapsto g\cdot o$.
Recall that $\ga$ as a cocharacter and therefore it defines a homomorphism
${\mathbb K}^*\to T$.

\begin{thm}
\cite[Theorem 9.11]{KLM3}.
\label{trcount}
$m_{\alpha,\beta}(\ga)$ equals the cardinality of $\mathcal{T}$.
\end{thm}

\begin{rem}
Instead of relating geometry of locally compact buildings to representation theory of $G^\vee$
as it is done in this paper, one can use the non-locally compact building
associated with the group $\ul{G}(\C((t)))$,
as it is done, for instance, in \cite{GL}. It was shown in \cite{KLM2} that the choice of
a field with discrete valuation (or, even, the entire Euclidean building)
is irrelevant as far as the existence of triangles with the
given $\De$-side lengths is concerned: What is important is the affine Weyl group.
For our purposes, moreover, it is more convenient to work with local fields and locally-compact buildings.
In particular, it allows us to compare structure constants of Hecke and representation rings.
\end{rem}

\section{Chains}
\label{sectionchains}

\subsection{Absolute chains}
\label{absolute}

Let $R$ be a root system on a Euclidean vector space $V$,
$W=W_{sph}$ be the finite Coxeter group associated with $R$, let
$W_{aff}$ denote the affine Coxeter group associated to $R$. Our
root system $R$ is actually the {\em coroot system} for the one
considered by Littelmann in \cite{Littelmann2}. Accordingly, we will
switch weights to coweights, etc. We pick a Weyl chamber $\Del$ for
$W$, this determines the positive roots and the simple roots in $R$.
Let $-\Del$ denote the negative chamber.

We get the Euclidean Coxeter complex $(A,W_{aff})$, where $A=V$ and
the spherical Coxeter complex $(S,W)$ where $S$ is the
unit sphere in $V$. By abusing notation we will also use the
notation $\De, -\De$ for the positive and negative chambers in
$(S,W)$. We will use the notation $(A,W,-\De)$ for a
Euclidean/spherical Coxeter complex with chosen negative chamber.
More generally, we will use the notation $(A,W_{aff},a)$ for a Euclidean
Coxeter complex with chosen alcove $a$.

\begin{figure}[tbh]
\centerline{\epsfxsize=4.5in \epsfbox{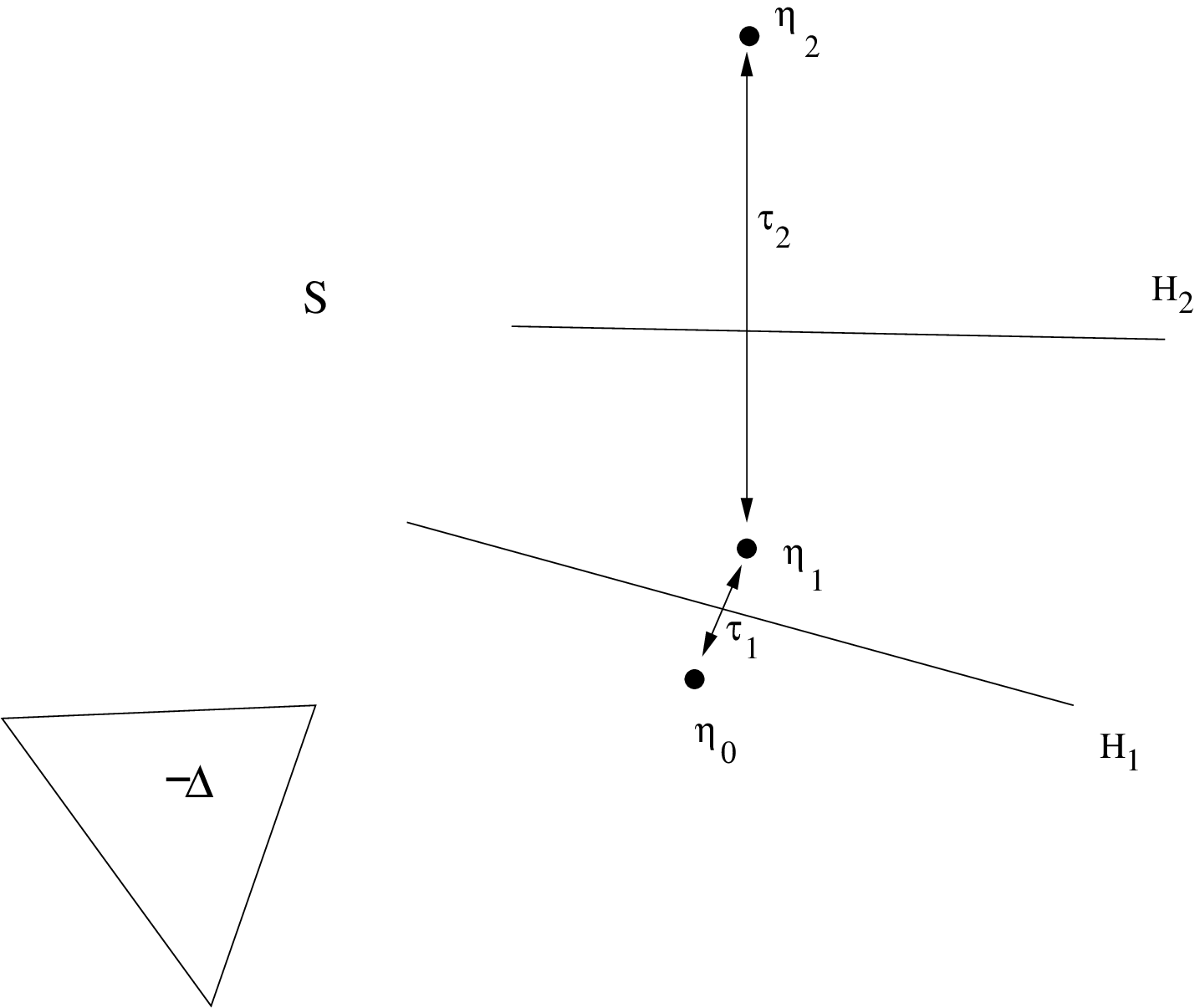}} \caption{\sl A
chain.} \label{Fc}
\end{figure}

\begin{defn}\label{D0}
 A $W$--{\em chain} in
$(A,W,-\De)$ is a finite sequence $(\eta_0,...,\eta_m)$ of nonzero
vectors in $V$ so that for each $i=1,...,m$ there exists a positive
root $\be_i\in R^+$ so that the corresponding reflection
$\tau_i:=\tau_{\be_i}\in W$ satisfies

1. $\tau_i(\eta_{i-1})= \eta_{i}$.

2. $\be_i(\eta_{i-1}) <0$.

\noindent Sometimes we will refer to a chain as a $(V,W,-\De)$-chain to
emphasize the choice of $V$, $W$ and $\De$. When the choice of $W$
is clear we will frequently refer to $W$-chains as {\em chains}.
\end{defn}

\begin{rem}
One could call these chains $(V,W,\De)$-chains instead,
but this definition would not generalize to affine chains.
\end{rem}

Recall that $\tau_i$ is the reflection in the wall
$H_i=\{\be_i=0\}$. More geometrically one can interpret the
condition (2) by saying that the wall $H_i$ separates the negative
chamber $-\De$ from the vector $\eta_i$. In other words, the
reflection $\tau_i$ moves the vector $\eta_{i-1}$ ``closer'' to the
positive chamber.

\bigskip
The concept of a chain generalizes naturally to Euclidean Coxeter complexes $(A,W_{aff})$. Pick an
alcove $a$ in $A$.

\begin{defn}
\label{DA}
An affine {\em chain} in $(A,W_{aff},a)$ is a finite sequence $(\eta_0,...,\eta_m)$ of
elements in $A$ so that for each $i=1,...,m$ there exists a reflection $\tau_i\in W_{aff}$ such that

1. $\tau_i(\eta_{i-1})= \eta_{i}$.

2. The hyperplane $H_i\subset A$ fixed by $\tau_i$ separates $a$ from $\eta_i$.
\end{defn}

We now return to the chains as in definition \ref{D0}.
By restricting vectors $\eta_i$ to have unit length we define chains
in the spherical Coxeter complex, see Figure \ref{Fc}.

\begin{defn}
The points $\eta_i$ as in Definition \ref{D0} will be called {\em vertices} of the
chain. We will say that the chain begins at $\eta_0$ and {\em ends}
at $\eta_m$, or that this chain is {\em from} $\eta_0$ {\em to}
$\eta_m$. We refer to a subsequence $(\eta_i,
\eta_{i+1},...,\eta_m)$ as a {\em tail} of the chain.
\end{defn}

We will refer to the number $m$ as the {\em length} of the chain. A
chain $(\eta_i)$ is called {\em simple} if it has length 1. Set
$$
dist_W(\nu,\mu)=dist(\nu, \mu)$$ to be the maximal length $m$ of a
$W$-chain which begins at $\nu$ and ends at $\mu$.

Given a chain $(\eta_0,...,\eta_m)$ we define a {\em subdivision} of
this chain to be a new chain in $(A,W)$ which is still a chain
from $\eta_0$ to $\eta_m$ and which contains all the vertices of
the original chain.

The concept of chain determines a partial order on the $W$-orbits in
$V$:

\begin{defn}
\cite[page 509]{Littelmann2}
\label{order}
 For a pair of nonzero vectors $\nu, \mu\in V$ which
belong to the same $W$-orbit, write $\nu\ge_{W} \mu$ (or simply
$\nu\ge \mu$) if there exists a $W$-chain from $\nu$ to $\mu$.
Accordingly, $\nu >\la$ if $\nu\ge \la$ and $\nu \ne \la$.
\end{defn}

\begin{lem}
\label{easy} Suppose that $\nu \ge \mu$ and $\al$ is a positive root
such that $\al(\mu)\le 0$. Then
$$
\nu\ge \tau_\al(\mu).
$$
\end{lem}
\proof If $\al(\mu)=0$ then $\tau_\al(\mu)=\mu$ and there is nothing
to prove. Thus we assume that $\al(\mu)<0$. Consider a chain
$(\nu=\nu_0,...,\nu_s=\mu)$. Then, since $\al(\mu)<0$, we also have
$\al(\tau_\al(\mu))>0$ and thus we get a longer chain
$$
(\nu=\nu_0,...,\nu_s=\mu, \tau(\mu)). \qed
$$

The word metric $d_W$ on the finite Coxeter group $W$ defines the
{\em length function}
$$
\ell: W\cdot \la \to \N
$$
by
$$
\ell(\mu):= \min\{ d_W(w,1):  w\in W, w^{-1}(\mu)\in \De\}.
$$

\begin{prop}
\label{wordlength} If $\nu > \mu$ then $\ell(\nu)> \ell(\mu)$.
\end{prop}
\proof It suffices to prove the assertion in the case when
$$
\mu=\tau(\nu), \tau=\tau_\be,
$$
where $\be$ is a positive root, $\be(\nu)<0, \be(\mu)>0$. If $W\cong
\Z/2$ the assertion is clear, so we suppose that it is not the case.
Then we can embed the Cayley graph $\Ga$ of $W$ as a dual graph to
the tessellation of $V$ by the Weyl chambers of $W$. Suppose that
$\nu \in w^{-1}(\De)$, then the wall $H_{\be}=\{\be=0\}$ separates
$w^{-1}(\De)$ from $\De$, where $w\in W$ is the shortest element
 such that $w(\nu)\in \De$. Let $p: [0,1]\to \Ga$
denote the shortest geodesic from $1$ to $w$ in $\Ga$. The path $p$
crosses the wall $H$ at a point $x=p(T)$. We construct a new path
$q$  by
$$
q|[0,T]=p|[0,T], \quad q|[T, 1]= w\circ p|[T,1].
$$
The path $q$ connects $1\in \De$ to the Weyl chamber $\tau w(\De)$
containing $\mu$. This path has a break-point at $x$, which is not a
vertex of the Cayley graph. Therefore, by eliminating the
backtracking of $q$ at $x$, we obtain a new path which connects $1$
to $\tau w(\De)$ and whose length is one less than the length of
$p$. \qed

\begin{cor}
\label{maxi} The length of a chain in $V$ does not exceed the
diameter of the Cayley graph of $W$.
\end{cor}

\begin{cor}
\label{end} Suppose that $\nu\in \De$ and $\nu\ge \mu$. Then
$\mu=\nu$.
\end{cor}
\proof Since $\nu\in \De$, $\ell(\nu)=0$. Hence by Proposition \ref{wordlength},
$\ell(\mu)= 0$, which implies  that $\mu=\nu$. \qed

\begin{lem}
\label{2orders}
Suppose that $\nu=w(\mu)\ne 0$ for some $w\in W$ and $\mu\rhd \nu$. Then
$\nu \ge \mu$.
\end{lem}
\proof Let $\De$ be the positive Weyl chamber.
Suppose that $\De_0$ is a chamber containing $\nu$ and $\De'$
is a chamber containing $\mu$. These chambers are non-unique, but
we can choose them in such a way that for all $\xi\in \De_0, \xi'\in \De'$
$$
\xi'\rhd \xi.
$$
In other words, if a wall separates $\De'$ from $\De$ then it also separates $\De_0$ from $\De$.
Let
$$
(\De_0, \De_1,...,\De_m=\De')
$$
be a gallery of Weyl chambers, i.e. for each $i$, $\De_i\cap
\De_{i+1}$ is a codimension 1 face $F_i$ of $\De_i, \De_{i+1}$. We choose this gallery to have
the shortest length, i.e. so that the number $m$ is minimal.
Let $H_i$ be the wall containing $F_i$ and $\tau_i$ be the reflection in $H_i$.
We claim that the sequence
$$
\nu=\eta_0, \eta_1=\tau_1(\eta_0),..., \mu=\eta_m=
\tau_m(\eta_{m-1}),
$$
after deletion of equal members, is a chain.

Our proof is by induction on $m$. If $m=0$ and $\nu=\mu$, there is
nothing to prove. Suppose the assertion holds for $m-1$, let us
prove it for $m$. We claim that for all points $\xi_1\in \De_1,
\xi_0\in \De_0$, $\xi_1 \rhd \xi_0$. Indeed, otherwise the wall
$H_1$ does not separate $\De_0$ from $\De$, but separates $\De_1$
from $\De$. Then $H_1$ does not separate $\De_m$ from $\De$
either. Thus, as in the proof of Proposition \ref{wordlength}, we can
replace the gallery $(\De_0, \De_1,...,\De_m)$ with a shorter
gallery connecting $\De_0$ to $\De'$, contradicting minimality of
$m$. Now, clearly,
$$
\eta_0\ge \eta_1, \eta_0\rhd \eta_1.
$$
Therefore, by the induction
$$
\eta_0\ge \eta_1 \ge \eta_m \Rightarrow \nu=\eta_0\ge \eta_m=\mu.  \qed
$$

\begin{figure}[tbh]
\centerline{\epsfxsize=4.5in \epsfbox{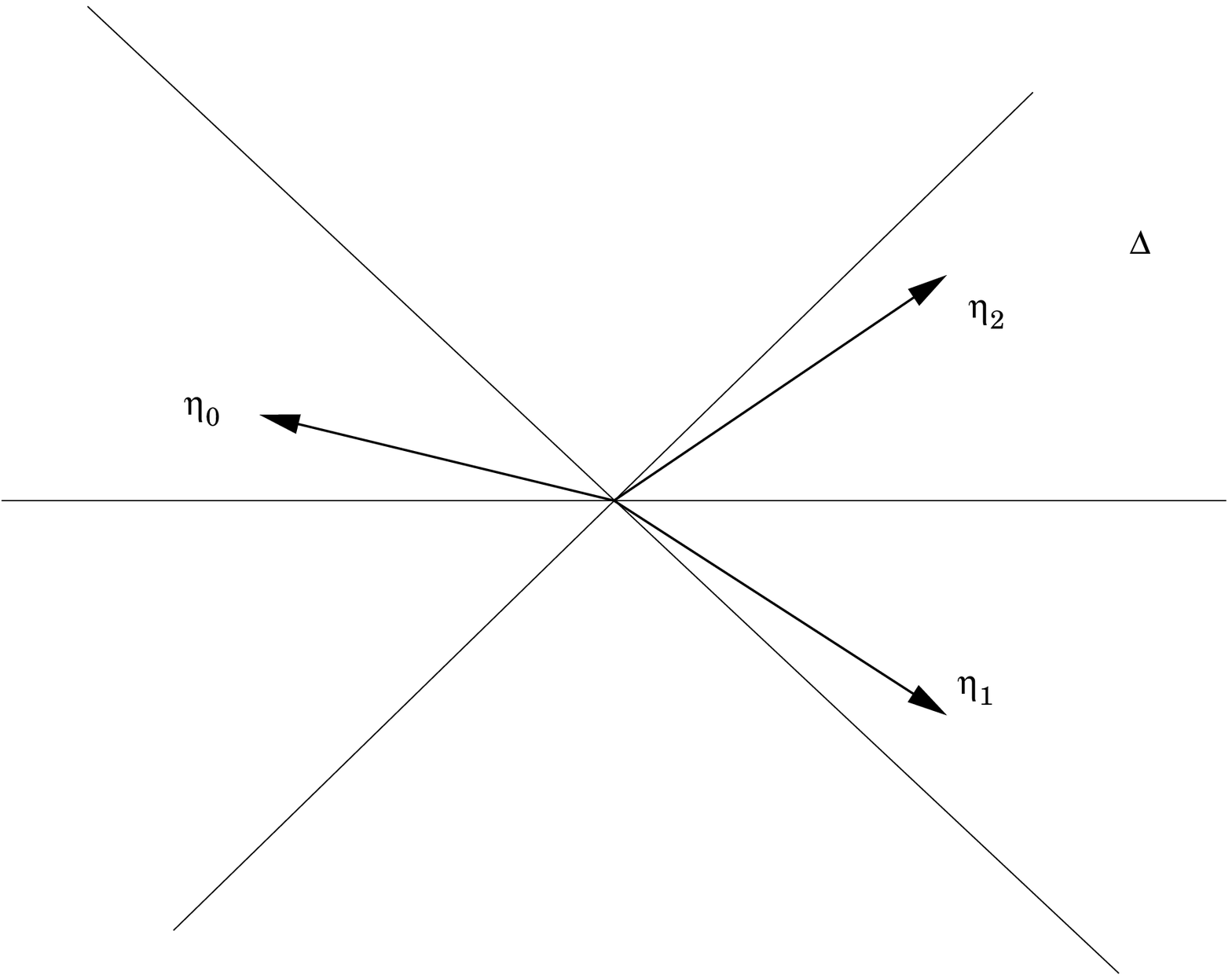}} \caption{\sl A
chain.} \label{Fchain}
\end{figure}

\begin{rem}
The converse to the above lemma is false for instance for
the root system $A_2$. See Figure \ref{Fchain}, where $\eta_0\ge \eta_1$ but  $\eta_0\ntriangleright \eta_1$.
\end{rem}

As a corollary of Lemma \ref{2orders} we obtain:

\begin{lem}
\label{constructingachain}
Let $\nu\in V\setminus \{0\}$ and let $\mu$ be the unique vector in
$W\cdot \nu$ which belongs to $\De$. Then $\nu\ge \mu$.
\end{lem}
\proof Clearly, $\mu \rhd \nu$. Then the assertion follows from
Lemma \ref{2orders}.\qed

\begin{defn}
[Maximality condition] \label{maximality}
  We say that a chain
$(\eta_0, \eta_1,..., \eta_m)$ is  {\em maximal} if it cannot be
{\em subdivided} into a longer $W$-chain. Equivalently, \newline
$dist(\eta_i,\eta_{i+1})=1$ for each $i=0,...,m-1$.
\end{defn}

\begin{lem}
\label{simplemax}
Suppose that $\nu\ge \mu$, and there exists a simple root $\be$ such that
$\tau_\be(\nu)=\mu$ and  $\be(\nu)<0$. Then $dist(\nu,\mu)=1$.
\end{lem}
\proof Consider a chain from $\nu$ to $\mu$, i.e. a sequence of vectors $\nu=\nu_0, \nu_1,...,
\nu_s=\mu$ and positive roots $\be_1,...,\be_s$ so that
$$
\nu_i= \tau_{\be_i}(\nu_{i-1}) \quad \hbox{and}\quad
\be_i(\nu_{i-1}) <0, i=1,...,s.
$$
Then
$$
\mu-\nu= -\be(\nu) \be^\vee,
$$
and
$$
\mu-\nu= \sum_{i=1}^s -\be_{i-1} (\nu_i)\be^\vee_i.
$$
Thus
$$
\be= \sum_{i=1}^s \frac{ \< \nu_i, \be_i^\vee\>}{\<\nu, \be^\vee\>}
\be_{i-1}
$$
i.e. the simple root $\be$ is a positive linear combination of
positive roots. It follows that $s=1$ and $\be_1=\be$. \qed

\subsection{Relative chains}

{\bf 1. Chains relative to a root subsystem.}

\medskip
Let $R$ be a root system on $V$ with the set of simple roots $\Phi$,
$W$ be the corresponding Weyl group. Let $\Phi'\subset \Phi$ be a
subset, $W'\subset W$ the corresponding reflection subgroup and
$\De'$ the positive chamber for $W'$, defined by the property that
all simple roots $\al\in \Phi'$ are nonnegative on $\De'$. Thus we
will be (frequently) considering $(V,W',-\De')$-chains rather than $(V,W,-\De)$-chains.
In this paper we will be using subgroups $W'$ which are stabilizers
of points $x$ in $(A,W_{aff})$ (where $W=W_{sph}$). Then we will
think of a relative $(V, W')$-chain as a chain in the tangent space
$T_x A$ of the point $x$.

\begin{lem}
\label{simplesubsystem} Given a nonzero vector $\eta\in V$ there
exists a $W'$-chain \newline $(\eta=\eta_0,...,\eta_m)$ with
$\eta_m\in \De'$, so that the chain $(\eta_i)$ is maximal with
respect to the original root system $R$.
\end{lem}
\proof The proof is by induction on the number $r=r(\eta)$ of simple
roots in $\Phi'$ which are negative on the vector $\eta$. If
$r(\eta)=0$ then $\eta\in \De'$ and there is nothing to prove.
Suppose the assertion holds for all $\eta$ with $r(\eta)\le k$. Let
$\eta$ be such that $r(\eta)=k+1$. Pick a root $\be\in \Phi'$ such
that $\be(\eta)<0$. Then for the vector
$$
\zeta:= \tau_{\be}(\eta)
$$
we have:
$$
\zeta= \eta- \be(\eta)\be^\vee.
$$
Clearly, $\beta(\zeta) >0$. Thus the pair $(\eta,\zeta)$ is a
$W'$-chain; this chain is maximal as a $W$-chain by Lemma
\ref{simplemax}, since it is defined using a simple reflection in
$W$.

For each simple root $\al\in \Phi'\setminus \{\be\}$ which is
nonnegative on $\eta$
$$
\al(\zeta) = \al(\eta)- \be(\eta) \al(\be^\vee)\ge \al(\eta) \ge 0.
$$
Therefore $r(\zeta)< r(\eta)$ and we are done by the induction. \qed

\begin{lem}
\label{ifff} Suppose that $W'\subset W$ is a reflection subgroup as
above. For any two vectors $\al, \del$ the following are equivalent:

1. There exists $\be, \ga$ so that $\al \ge_{W'} \be \sim_{W} \ga
\ge_{W'}  \del$.

2. There exists $\be$ so that $\al \ge_{W'}  \be \sim_{W} \del$.

3. There exists $\ga$ so that $\al \sim_{W} \ga \ge_{W'}  \del$.

\noindent Here $\sim$ is the relation from Definition \ref{relation}.
\end{lem}
\proof It is clear that 2 $\Rightarrow$ 1 and 3 $\Rightarrow$ 1. We
will prove that 1 $\Rightarrow$ 2, since the remaining implication
is similar. We have chains
$$
(\al=\eta_0,...,\eta_m=\be), \quad \eta_i=\tau_i(\eta_{i-1}), \ i=1,...,m,
$$
$$
(\ga=\eta_0',...,\eta_s'=\del), \quad \eta_i'=\tau_i'(\eta_{i-1}'), \ i=1,...,s.
$$
Then we can extend the first chain to
$$
(\al=\eta_0,...,\eta_m=\be, \tau_1'(\be)=\eta_{m+1},...,
\tau_s'(\eta_{m+s-1})=\eta_{m+s}=:\eps).
$$
After discarding equal members of this sequence we obtain a chain
from $\al$ to $\eps$. Since $\be, \ga$ belong to the same chamber,
the vectors $\eps$ and $\del=\tau_s'\circ ... \circ \tau_1'(\ga)$
also belong to the same chamber. Therefore we obtain
$$
\al\ge \eps \sim \del. \qed
$$

\begin{defn}
We will write $\al\gtrsim_{W'}  \del$ if one of the equivalent
conditions in the above lemma holds. We will frequently omit the
subscript $W'$ in this notation when the choice of the subgroup $W'$
is clear or irrelevant.
\end{defn}

The reason for using the relation $\sim$ with respect to $W$ rather
than its subgroup $W'$ (as may seem more natural) is that we will be
taking limits under which the subgroup $W'$ is increasing (but the
limit is still contained in $W$). Such limits clearly preserve the
relation $\gtrsim_{W'}$, but not the relation defined using
$\sim_{W'}$.

The next corollary immediately follows from Lemma \ref{ifff}:

\begin{cor}
\label{iff}
 $\al \gtrsim \del \iff -\del \gtrsim -\al$.
\end{cor}

\begin{lem}
\label{hard} Suppose that $\nu \gtrsim \mu$ and $\al\in \Phi$ is
such that $\al(\nu)<0$, $\al(\mu)\ge 0$. Then
$$
\tau_\al(\nu)\gtrsim \mu.
$$
\end{lem}
\proof Let $\la$ be such that
$$
\nu \ge \la \sim \mu.
$$
Then $\al(\la)\ge 0$ and it follows that $\la\ne \nu$, i.e. $\nu >
\la$. By applying \cite[Lemma 4.3]{Littelmann2} we get
$$
\tau_\al(\nu) \ge \la \sim \mu. \qed
$$

\medskip
{\bf 2. Chains relative to positive real numbers ($a$-chains in the
sense of Littelmann).}

Let $a$ be a positive real number
and let $\nu, \mu\in V$ be nonzero vectors in the same $W$-orbit.

\begin{defn}
[P.~Littelmann, \cite{Littelmann2}] \label{achain}  {\em An
$a$-chain} for $(\nu,\mu)$ is a chain $(\la_0,\ldots,\la_s)$ which
starts at $\nu$, ends at $\mu$ and satisfies

(i) For each $i>0$ we have
$$
t_i:=\be_i(a\la_{i-1})\in \Z,
$$
where $\la_i=\tau_{\be_i}(\la_{i-1})$ as in the Definition
\ref{D0}.

(ii) For each $i$, $dist(\la_{i-1}, \la_i)=1$.
\end{defn}

\begin{rem}
Our root system $R$ is the coroot system for the one considered by
Littelmann.
\end{rem}

Our goal is to give this definition a somewhat more geometric
interpretation. In particular, we will see that the concept of an
$a$-chain is a special case of the concept of a chain relative to a
root subsystem.

The root system $R$ defines an affine Coxeter complex $(A, W_{aff})$
on $A$. Let $x\in P(R^\vee)$ be a special vertex; set $x_i:=
x+a\la_i$, $\tau_i:=\tau_{\be_i}$, $i=0,...,s$. Thus $t_i=
\be_i(x_{i-1})$. Note that $t_i\in \Z$ iff $\be_i(x_{i-1})\in \Z$.

\begin{prop}
$\be_i( x_0)\in \Z$ for each $i=1, 2,...,s$.
\end{prop}
\proof It suffices to consider the case $x=0$. We have:
$$
x_{i-1}= \tau_{i-1}(x_{i-2})= x_{i-2} - \be_{i-1}( x_{i-2})
\be_{i-1}^\vee= x_{i-2} - t_{i-1} \be_{i-1}^\vee= \ldots = x_0-
\sum_{j=1}^{i-1}  t_j \be_{j}^\vee.
$$
Hence
$$
t_i= \be_i( x_{i-1})= \be_i( x_0) -  \sum_{j=1}^{i-1}  t_j
\be_i(\be_{j}^\vee).
$$
Since $t_i\in \Z$ and $\be_i(\be_{j}^\vee)\in \Z$ for all $j$, it
follows that $\be_i( x_0)\in \Z$. \qed

\medskip
We define the integers $k_i:= \be_i(x_0)$ and the affine walls
$$
H_i:=H_{\be_{i},k_{i}}= \{ v\in V: \be_{i}(v)=k_i\}.
$$
The reflection $\si_i$ in  the wall $H_i$ belongs to the group
$W_{aff}$, its linear part is the reflection $\tau_i\in W_{sph}$,
$i=1,...,s$.

The argument in the above proof can be easily reversed
and hence we get

\begin{cor}
\label{C1} The integrality condition (i) is equivalent to the
assumption that the point $x_0$ lies on the intersection of walls
$H_i$ of the Euclidean Coxeter complex $(A,W_{aff})$, where each
$H_i$ is parallel to the reflection hyperplane of $\tau_i$.
Equivalently, the $W$-chain $(\la_0,...,\la_s)$ is actually a
$W'$-chain, where $W'=W_{x_0}$ is the stabilizer of $x_0$ in
$W_{aff}$.
\end{cor}

Therefore, identify the vectors $\la_i$ with vectors in the tangent
space $V':= T_{x_0}(A)$, let $\De'\subset V'$ denote the Weyl
chamber of $W'$ which contains the (parallel transport of the)
positive chamber $\De$. We obtain

\begin{prop}
\label{eq}  Littelmann's definition of an $a$-chain is equivalent to
the conjunction of

1. $(\la_0,...,\la_s)$ is a chain in $(V',W',-\De')$.

2. This chain is maximal as a $W$-chain.
\end{prop}

Thus the choice of the real number $a$ amounts to choosing a
 Coxeter subcomplex $(V', W')$ in $(V,W_{sph})$. The reader will also note
the discrepancy between (1) and (2): The chain condition refers to
the restricted Coxeter complex $(V', W')$, while the maximality
condition refers to the unrestricted one, $(V,W_{sph})$. This is the
key difference between {\em LS paths} and {\em Hecke paths}.

\begin{rem}
\label{R1} Note that both conditions (1) and (2) are vacuous if
$x_0$ is a special vertex in the Euclidean Coxeter complex
(equivalently, if $a\nu$ is a coweight): If
$$
(\la_0, \la_1,...,\la_s)
$$
is a $W'$-chain, since $W'\cong W_{sph}$,  we can subdivide this
chain to get a longer chain
$$
(\la_0=\la'_0, \la_1',...,\la_{m-1}', \la_m'=\la_s)
$$
between $\la_0$ and $\la_s$ which satisfies the unit distance
condition  $dist(\la_i', \la_{i+1}')=1$ for all $i$.
\end{rem}

\subsection{Hecke paths}

The goal of this section is to introduce a class of piecewise-linear paths which
satisfy a condition similar to Littelmann's definition of LS (and
generalized LS) paths. These paths ({\em Hecke paths}) play the role
in the problem of computing structure constants for spherical Hecke
algebras which is analogous to the role that LS paths play in the
representation theory of complex semisimple Lie groups.

\medskip
Let $(A, W_{aff}, -\De)$ be a Euclidean Coxeter complex
corresponding to a root system $R$, with fixed negative chamber
$-\De$.  Let $p\in \tilde{\mathcal P}$ be a path equal to the
composition
$$
\ol{x_{1} x_{2}}\cup ... \cup \ol{x_{n-1} x_{n}}.
$$
For each vertex $x=x_i, i=2,...,n-1$ we define the unit tangent
vectors $\xi, \eta \in S_x$ to the segments $\ol{x_ix_{i-1}},
\ol{x_ix_{i+1}}$.

\begin{defn}
\label{D3} We say that the path $p$ {\em satisfies the chain
condition} if for each $x=x_i, i=2.,...,n-1$ there exists a unit
vector $\mu$ so that

1.
$$
-\xi \ge_{W_x} \mu, \hbox{in the spherical Coxeter complex~}
(S_x,W_x,\De).
$$
2. $\mu \sim \eta$ in the (unrestricted) spherical Coxeter complex
$(S, W_{sph})$, i.e. for each root $\al\in R$ we have
$$
\al(\mu) \ge 0 \iff \al(\eta) \ge 0.
$$
In other words, for each $t\in [0, 1]$ we have
$$
p'_-(t)\gtrsim_{W_{p(t)}} p'_+(t).
$$
\end{defn}

\noindent Intuitively, at each break-point $p(t)$ the path $p$
``turns towards the positive chamber''. In what follows we will use
the notation ${\mathcal P}_{chain}$ for the set of all paths $p\in
{\mathcal P}$ satisfying the chain condition.

\begin{defn}
\label{Hecke} A path $p\in \t{\mathcal P}$ is called a {\em Hecke
path} if it is a billiard path which satisfies the chain condition, i.e.
for each $t$
$$
p'_-(t)\ge_{W_{p(t)}} p'_+(t).
$$
\end{defn}

Below is example of a class of Hecke paths. Suppose that $p\in
{\mathcal P}$ and for each $t\in [0, 1]$ either $p$ is smooth at $t$
or there exists a reflection $\tau\in W_{p(t)}$ so that the
derivative of $\tau$ equals $\tau_\be$, $\be \in R^+$ and

1. $d\tau_{p(t)} (p'_-(t))=p'_+(t)$.

2. $\be(p'_-(t))<0$, $\be(p'_+(t))>0$.

\begin{figure}[tbh]
\centerline{\epsfxsize=3.5in \epsfbox{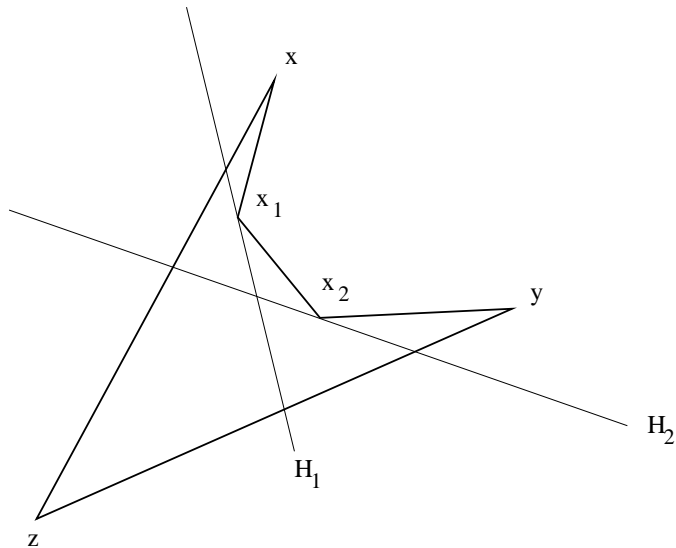}} \caption{\sl A
billiard path satisfying the simple chain condition.} \label{Fa}
\end{figure}

Then $p$ is a Hecke path with the length of each chain in
$(S_{p(t)}, W_{p(t)})$ equal $0$ or $1$. See Figure \ref{Fa}.

\begin{defn}
We say that $p$ satisfies the {\em simple chain} condition if at
each break-point $x=p(t)$ the chain can be chosen to be simple, i.e.
of length $1$.
\end{defn}

In what follows we will also need

\begin{defn}
\label{gH}
Suppose that $p$ is a path satisfying the chain condition. We call $p$ a {\em generalized Hecke path} if
each geodesic segment in $p$ (regarded as a vector)
is an integer multiple of some $w\varpi$ where $w\in W$ and $\varpi=\varpi_i$
is one of the fundamental coweights.
\end{defn}

\subsection{A compactness theorem}

Pick $\eps>0$. We define the subset ${\mathcal P}_{m,\eps}\subset
{\mathcal P}$ consisting of paths $p$ with
$$
\length(p)=\ul{\la}=(\la_1,...,\la_m)$$
 so that for each $i$,
\begin{equation}
\label{bdd}
\eps\le |\la_i|\le \eps^{-1}.
\end{equation}

\begin{thm}
\label{compactness} For each $\eps>0$ the set ${\mathcal
P}_{chain,m,\eps}:= {\mathcal P}_{chain}\cap {\mathcal P}_{m,\eps}$
is compact in ${\mathcal P}_{m,\eps}$.
\end{thm}
\proof Suppose that $p\in {\mathcal P}_{chain}$ is a concatenation
of $m$ billiard paths $p_i$. Then the number of breaks in the broken
geodesic $p_i$ is bounded from above by a constant $c$ equal to the
length of a maximal chain in the Bruhat order of the finite Weyl
group $W_{sph}$, see Corollary \ref{maxi}. This immediately implies
that the subset ${\mathcal P}_{chain,m,\eps}$ is precompact in
${\mathcal P}$. What has to be proven is that this subset is closed.

\bigskip
Let $\la_{ij}\in V$ be vectors of nonzero length so that
$$
p_i= \pi_{\la_{i1}} * \pi_{\la_{i2}}*...*\pi_{\la_{in}}
$$
are Hecke paths, $i=1, 2,...$. We suppose that
$$
\lim_{i\to\infty} \la_{ij}=0, j=2,...,n-1
$$
and
$$
\lim_{i\to\infty} \la_{i1}=\la_{\infty,1}, \lim_{i\to\infty} \la_{in}=\la_{\infty, n}
$$
are nonzero vectors. It is clear that
$$
\lim_i p_i= p_{\infty}:= \pi_{\la_{\infty,1}} * \pi_{\la_{\infty,n}}.
$$

\begin{lem}
\label{1}
Under the above conditions the path $p_{\infty}$ is again a Hecke path.
\end{lem}
\proof Let $x:= \la_{\infty,1}$. We need to check that the unit vectors $\bar\la_{\infty,1}, \bar\la_{\infty,n}$ satisfy
$$
\bar\la_{\infty,1} \ge_{W_x} \bar\la_{\infty,n}.
$$
Here and below, $W_x$ is the stabilizer of $x\in V$ in the Coxeter group $W_{aff}$.

Let $x_{ij}$ denote the break-point of $p_i$ which is the concatenation point between $\la_{i,j}$ and $\la_{i,j+1}$. Then
$$
\bar\la_{i,j} \ge_{W_{x_{ij}}} \bar\la_{i,j+1},
$$
$\lim_i x_{ij}=x$. If $\si_{ij}\in W_{aff}$ is a reflection fixing $x_{ij}$,
then, up to a subsequence,
$$
\si_{\infty,j}= \lim_i \si_{ij}\in W_{aff}
$$
fixes the point $x$. Therefore it follows from the definition of a chain that
$$
\bar\la_{\infty,j} \ge_{W_x} \bar\la_{\infty,j+1}
$$
for each $j$. By putting these inequalities together we obtain
$$
\bar\la_{\infty,1} \ge_{W_x} \bar\la_{\infty,n}. \qed
$$

\medskip
Suppose now that $p_i, q_i$ are Hecke paths,
$$
p_i= \pi_{\la_{i1}} * \pi_{\la_{i2}}*...*\pi_{\la_{in}}, q_i= \pi_{\mu_{i1}} * \pi_{\mu_{i2}}*...*\pi_{\mu_{im}}
$$
and the concatenation $r_i:= p_i*q_i$ satisfies the chain condition. We suppose that $r_{\infty}=p_\infty * q_\infty$
is the limit of the sequence of paths $r_i$ and the paths $p_\infty, q_\infty$ are not constant.

\begin{lem}
\label{2}
Under the above assumptions, the path $r_\infty$ also satisfies the chain condition.
\end{lem}
\proof Lemma \ref{1} implies that the paths $p_\infty, q_\infty$ are Hecke. Therefore it suffices to verify the chain condition
at the concatenation point $x=p_\infty(1)$. Up to passing to a subsequence, we have:
$$
\lim_{i\to\infty} \la_{ij}=\la_{\infty,j}, \lim_{i\to\infty} \mu_{ij}=\mu_{\infty,j},
$$
$$
\lim_{i\to\infty} \bar\la_{ij}=\bar\la_{\infty,j}, \lim_{i\to\infty} \bar\mu_{ij}=\bar\mu_{\infty,j}.
$$

Suppose that
$$
\lim_i \la_{ij}=0, j=k+1,...,n, \quad \lim_i \mu_{ij}=0, j=1,...,l,
$$
and
$$
\la_{\infty,k}\ne 0, \quad \mu_{\infty,l+1}\ne 0.
$$

By Lemma \ref{1},
$$
\bar\la_{\infty,k} \ge_{W_x} \bar\la_{\infty,n}, \quad \bar\mu_{\infty,1} \ge_{W_x} \bar\mu_{\infty,l+1}.
$$
On the other hand, it is clear that
$$
\bar\la_{\infty,n} \gtrsim_{W_x} \bar\mu_{\infty,1}.
$$
Therefore, by Lemma \ref{ifff}, we get
$$
\bar\la_{\infty,k} \gtrsim_{W_x} \bar\mu_{\infty,l+1}. \qed
$$

We now can finish the proof Theorem \ref{compactness}. Suppose that $p_i=p_{i1}*...*p_{im}$ is a sequence of paths
in ${\mathcal P}_{chain,m,\eps}$, where each $p_{ij}$ is a billiard path, and $p_{\infty}$ is the limit of this sequence.
Each sequence of billiard paths $(p_{ij})_{i\in \N}$ converges to a billiard path $p_{\infty,j}$ which is a Hecke path according to Lemma \ref{1}.
Consider now a concatenation point $x$ of the subpaths $p_{\infty}:=p_{\infty,j}, q_{\infty}:= p_{\infty,j+1}$. The paths
$p_\infty, q_\infty$ are non-constant by the inequality (\ref{bdd}). Therefore we
can apply Lemma \ref{2} to conclude that the path $p_\infty$ satisfies chain condition at the point $x$. \qed

\begin{rem}
It is easy to see that the assumption that the length of each billiard subpath $p_{ij}$ is bounded away
from zero, is necessary in Theorem \ref{compactness}.
\end{rem}

\section{Folding}
\label{sectionfolding}

The key tool for proving the main results of this paper is {\em
folding} of polygons in a building $X$ into apartments and Weyl
chambers. The folding construction replaces a geodesic segment $\t{p}$ in $X$
with a piecewise-linear path $p$ in an apartment.
This construction was used in \cite{KLM3} to construct various
counter-examples.  The reader will note that the folding
construction used in the present paper is somewhat different from
the one in \cite{KLM3}.

\subsection{Folding via retraction}

Suppose that $X$ is a Euclidean or spherical building modeled on the
Coxeter complex $(A,W)$, we identify the model apartment $A$ with an
apartment $A\subset X$; let $a\subset A$ be an alcove (or a chamber
in the spherical case). Recall that the {\em retraction, or folding,
to an apartment} $f=Fold_{a,A}: X\to A$ is defined as follows (see
for instance \cite{Rousseau}):

Given a point $x\in X$ choose an apartment $A_x$ containing $x$ and
$a$. Then there exists a (unique) isomorphism $\phi: A_x\to A$
fixing $A\cap A_x$ pointwise and therefore fixing $a$ as well. We
let $f(x):= \phi(x)$. It is easy to see that $f(x)$ does not depend
on the choice of $A_x$. Observe that $f$ is an isometry on each
geodesic $\ol{xy}$, where $y\in a$.

\medskip
The retraction can be generalized as follows.

Suppose that $X$ is a Euclidean building, $a$ is an alcove with a
vertex $v$ (not necessarily special). Let $\De\subset A$ denote a
Weyl chamber with  tip $o$. Choose a dilation $h\in Dil(A,W)$ which
sends $v$ to $o$. Let ${\mathbb P}: A\to \De$ denote the natural
projection which sends points $x\in A$ to $W_{sph}\cdot x \cap
\Del$, where $W_{sph}$ is the stabilizer of $o$ in $W_{aff}$.

We define a folding $g=Fold_{v,h,\De}: X\to \De$ as the composition
$$
\P\circ h\circ Fold_{a,A}.
$$
The mapping $g$ will be called a {\em folding} into a Weyl chamber.
Observe that $g$ (unlike $Fold_{a,A}$) does not depend upon the
choice of the alcove $a$, therefore it will be denoted in what
follows $g=Fold_{v,h,\De}$. Note that
$$
Fold_{v,k\circ h,\De}= k\circ Fold_{v,h,\De}.
$$

In case when $h=Id$ we will abbreviate $Fold_{v,h,\De}$ to
$Fold_{\De}$.

\begin{rem}
The folding maps are Lipschitz and differentiable. The restriction
of the retraction $Fold_{a,A}$ to each chamber (alcove) is a
congruence of two chambers (alcoves).
\end{rem}

%\begin{rem}
%Observe that instead of retracting to an apartment one can retract
%$X$ to the parallel set of a flat in $X$. We will explore this in \cite{HKM}.
%\end{rem}

If $h=Id$ (and thus $v=o$) one can describe $f=Fold_{\De}$ as
follows. Given a point $x\in X$ find an apartment $A_x$ through $v,
x$ and a Weyl chamber $\De_x\subset A_x$ with  tip $v$. Let
$\phi: \De_x\to \De$ be the unique isometry extending to an
isomorphism of Coxeter complexes $A_x\to A$. Then $f(x)=\phi(x)$.

\medskip
Suppose now that $X$ is a Euclidean building, $x\in X$; we give the
link $\Si_x(X)$ structure of an unrestricted spherical building $Y$.
Let $R$ denote the corresponding root system. Let $\del$ be a
chamber in $Y$ and $\xi, \mu\in \del$. Let $f: X\to \Del$ be a
folding of $X$ to a Weyl chamber. Let $x'$, $\xi', \mu'$ denote the
images of $x$, $\xi, \mu$ under $f$ and $df_x$. Then

\begin{lem}
\label{positivity} 1. $d_{ref}(\xi, \mu)= d_{ref}(\xi', \mu')$.

2. For each $\al\in R$,
$$
\al(\mu')\ge 0 \iff \al(\xi')\ge 0.
$$
\end{lem}
\proof The restriction $df|_\del$ is an isometry which is the
restriction of an isomorphism of spherical apartments. This proves
(1). To prove (2) observe that $df$ sends $\del$ to a spherical Weyl
chamber $\del'$ in $\Si_{x'}X$. \qed

\medskip
Let $f$ be a folding of $X$ into an apartment or a chamber.

\begin{lem}
\label{broken} For each geodesic segment $\ol{xy}\subset X$ its
image $f(\ol{xy})$ is a broken geodesic, i.e. it is a concatenation
of geodesic segments.
\end{lem}
\proof We give a proof in the case of a folding into an apartment
and will leave the other case to the reader. Let $A'\subset X$
denote an apartment containing the geodesic segment $\ol{xy}$. Let
$a_1,...,a_m$ denote the alcoves (or chambers in the spherical case)
in $A'$ covering $\ol{xy}$, set $\ol{x_i x_{i+1}}:= \ol{xy}\cap
a_i$. For each $a_i$ there exists an apartment $A_i$ containing the
alcoves $a$ and $a_i$. The restriction of the retraction $f$ to
$A_i$ is an isometry. It is now clear that the path
$f(\ol{xy})$ is a composition of the geodesic paths $f(\ol{x_i
x_{i+1}})$.
 \qed

We let $x_i'=f(x_i)$ denote the break points of $f(\ol{xy})$. For
each $x_i'$ let $\xi_i', \eta_i'$ denote the unit tangent vectors in
$T_{x_i'} A$ which are tangent to the segments $\ol{x_i' x_{i-1}'},
\ol{x_i' x_{i+1}'}$ respectively.

\begin{lem}
\label{billiard}
 The broken geodesic $f(\ol{xy})$ is a {\em billiard
path}, i.e. for each break point $x_i'$ the vectors  $\xi':=\xi_i',
\eta':=\eta_i'$ satisfy
$$
\exists w\in W_{x_i'} : w(\xi')=-\eta'.
$$
\end{lem}
\proof We again present a proof only in the case of a folding into
an apartment. Let $Y=\Si_{x_i} X$ denote the spherical building
which is the space of directions of $X$  at $x_i$. Let $\De_Y$
denote the Weyl chamber of this building and $\theta: Y\to \De_Y$
the canonical projection. The directions $\xi=\xi_i$ and
$\eta=\eta_i$ of the segments $\ol{x_i x_{i-1}}, \ol{x_i x_{i+1}}$
are antipodal in the building $Y$. Since the folding $f$ is an
isomorphism of the apartments $A_i\to A, A_{i+1}\to A$, and
$$
df(\eta)=\eta', \quad df(\xi)=\xi',
$$
we see that
$$
\theta(\xi)=\theta(\xi'), \quad \theta(\eta)=\theta(\eta').
$$
The assertion now follows from Lemma  \ref{thetaproperties}, part
(2). \qed

\begin{lem}
\label{invarianceoflength} Suppose that $X$ is a Euclidean building,
$f=Fold_{a,A}: X\to A$ and $g=Fold_{z,h,\De}$ are foldings to an
apartment and a chamber respectively. Then for each piecewise-linear path $p$ in
$X$ we have:

1. $\length(f(p))=\length(p)$.

2. $\length(g(p))=k\cdot \length(p)$, where $k>0$ is the conformal
factor of the dilation $h$.
\end{lem}
\proof We will prove the first assertion since the second assertion
is similar. It suffices to give a proof in the case when $p$ is a
billiard path. Then, analogously to the proof of Lemma
\ref{billiard}, there exists a representation of $p$ as a composition
of geodesic subpaths
$$
p=p_1\cup ...\cup p_m
$$
so that the restriction of $f$ to each $p_i$ is a congruence.
Therefore
$$
\Length(p_i)= \Length(f(p_i))
$$
and hence
$$
\Length (p)= \sum_i \Length(p_i)= \Length (f(p)). \qed
$$

\medskip
{\bf Derivative of the retraction.} We assume that $rank(X)\ge 1$.
We identify the model apartment $A$ with an apartment in $X$. Pick
$a\subset A$ which is an alcove (in the Euclidean case) or a chamber
(in the spherical case). Given a point $x'\in X$ choose an apartment
$(A',W')$ through $a$ and $x'$ and let $\phi: A\to A'$ denote the
inverse to the retraction $f=Fold_{a,A}: A'\to A$. Set $x=f(x')$ and
let $W_x'$ denote the stabilizer of $x'$ in $W'$. Then the link
$Y=\Si_{x'}(X)$  has a natural structure of a thick spherical
building modeled on $(S,W_x')$. It is easy to see that $(S, W_x')$
is independent of the choice of $A'$. Observe that if $x'$ is
antipodal to a regular point $y\in a$ then $x'$ is regular itself
and therefore $W_x'=\{1\}$. We next define a chamber $s\subset S$:

Given a regular point $y\in a\setminus \{x'\}$ and a geodesic
segment $\ol{x'y}$, let $\zeta=\zeta(y)$ denote the unit tangent
vector to $\ol{x'y}$ at $x'$. Then the set
$$
\{\zeta(y): y \hbox{~is a regular point in~} a\},
$$
is contained in a unique spherical chamber $s\subset S$. (If $x'$ is
antipodal to some $y\in int(a)$ then $s=S$.)

Set $f':= \phi\circ f= Fold_{a,A'}$.

\begin{lemma}
\label{consistent}
 The derivative $d_{x'}f': Y\to S$ equals $Fold_{s,S}$.
\end{lemma}
\proof Given $\eta\in Y$, find an alcove (or a spherical chamber)
$c$ so that $\eta\in \Si_{x'} c$. Then there exists an apartment
$A_\eta\subset X$ containing both $a$ and $c$. Let $S_\eta$ denote
the unit tangent sphere of $A_\eta$ at $x'$. Then $\eta\in S_\eta$
and $s\subset S_\eta$. Now it is clear from the definition that
$$
df'(\eta)=Fold_{s,S}(\eta)$$
 since both maps send $S_\eta$ to $S$ and fix $s$ pointwise. \qed

\medskip
{\bf Folding of polygons.} Suppose now that $X$ is a building and
$\t{P}=[\t{z}, \t{x}_1, \ldots, \t{x}_n]$ is a geodesic polygon in $X$. Pick an
apartment $A\subset X$ which contains $\ol{\t{z} \t{x}_1}$ and an alcove
$a\subset A$ which contains $\t{z}$. Let $\De\subset A$ denote a Weyl
chamber (in case $X$ is Euclidean) with  tip $o$. Let $f$ be a
folding of $X$ of the form
$$
f=Fold_{a,A}$$
 or
$$
f=Fold_{\t{z},h,\De},$$
 where $h$ is a dilation sending $\t{z}$ to $o$. We will
then apply $f$ to $\t{P}$ to obtain a {\em folded polygon} $P:= f(\t{P})$
in $A$ or $\De$ respectively.

Observe that the restriction of $f$ to the edges $\ol{\t{z} \t{x}_1}$ and
$\ol{\t{x}_n \t{z}}$ of $P$ is an isometry or a similarity. The restriction
of $f$ to the path
$$
\t{p}=\ol{\t{x}_1 \t{x}_2}\cup ...\cup \ol{\t{x}_{n-1}\t{x}_n}
$$
preserves the {\em type} of the unit tangent vectors, cf. Lemma
\ref{invarianceoflength}. We will be using foldings into apartments
and chambers to transform geodesic polygons in $X$ into {\em folded
polygons}.

In the special case when $\t{P}=T$ is a triangle (and thus $n=2$), the
folded triangle $P=f(T)$ has two geodesic sides $\ol{z w_1}:=f(\ol{\t{z} \t{x}_1}),
\ol{x_2 z}:= f(\ol{\t{x}_2 \t{z} })$ and one {\em broken side} $p:=f(\ol{\t{x}_1 \t{x}_2})$, so we
will think of $f(T)$ as a {\em broken triangle}.

The next proposition relates folding into a Weyl chamber with the
concept of folding of polygons used in \cite{KLM3}. Let $P=[o,
x_1, x_2, \ldots , x_n]$ be a polygon in $\De$. Triangulate $P$
from the vertex $o$ into geodesic triangles $T_i=[x_i,
x_{i+1},o]$. Suppose that $\t{P}\subset X$ is a geodesic polygon
$$
\t{P}=[o, \t{x}_1, \t{x}_2, \ldots, \t{x}_n], \t{x}_1=x_1,
$$
triangulated into geodesic triangles $\t{T}_i=[\t{x}_i, \t{x}_{i+1},o]$, where
each $\t{T}_i$ is contained in an apartment $A_i$. Assume that for each
$i$ there exists a congruence
$$
\phi_i: \t{T}_i\to T_i
$$
i.e. an isometry sending $\t{x}_j$ to $x_j$ ($j=i, i+1$) which extends
to an isomorphism of Coxeter complexes $\phi_i: A_i\to A$.

\begin{prop}
\label{uniqueness}
 Under the above assumptions, for each $i$,
$Fold_{\De}|\t{T}_i=\phi_i|\t{T}_i$.
\end{prop}
\proof Let $\De_i\subset A_i$ denote the preimage of $\De$ under
$\phi_i$. Then each $\De_i$ is a Weyl chamber, hence
$\phi_i|\De_i=Fold_{\De}|\De_i$, by the alternative description of
$Fold_{\De}$ given earlier in this section. \qed

\medskip
The following lemma shows that unfolding of polygons is a {\em local
problem}. Suppose that $T=[o, x_1,...,x_n]\subset A$ is a
geodesic polygon so that $x_i\ne o$ for each $i$. For each $i=2,...,
n-1$ we define the unit vectors
$$
\xi_i, \eta_i, \zeta_i\in \Si_{x_i} A
$$
which are tangent to the segments $\ol{x_i x_{i-1}},
\ol{x_{i}x_{i+1}}$, $\ol{x_i o}$. Define thick spherical
buildings $Y_i:= \Si_{x_i}(X)$. By combining the above proposition
with \cite[Condition 7.5]{KLM3} we obtain

\begin{lem}
\label{locality} The polygon $T$ can be unfolded in $X$ to a
geodesic triangle $\t{T}$ whose vertices project to $o, x_1, x_n$ if
and only if for each $i=2,...,n-1$, there exists a triangle
$[\t\xi_i, \t\zeta_i, \t\eta_i]\subset Y_i$ so that
$$
d_{ref}(\t\xi_i, \t\zeta_i) = d_{ref}(\xi_i, \zeta_i),
$$
$$
d_{ref}(\t\eta_i, \t\zeta_i) = d_{ref}(\eta_i, \zeta_i),
$$
$$
d(\t\xi_i, \t\eta_i) = \pi.
$$
\end{lem}

\medskip
We will eventually obtain a characterization of the broken triangles
in $\De$ which are foldings of geodesic triangles in $X$ as {\em
billiard triangles satisfying the chain condition}, see section
\ref{last}. The goal of the next section is to give a {\em
necessary} condition for a broken triangle to be unfolded; we also
give a partial converse to this result.

\subsection{Converting folded triangles in spherical buildings into chains}
\label{converting}

Suppose that $X$ is a (thick) spherical or Euclidean building modeled on $(A,W)$. Consider a
triangle $\tilde{T}=[x, y, z]\subset X$ with
$$
\be:=d_{ref}(x, y), \quad \ga:=d_{ref}(y, z).
$$
 Assume that $A$ is embedded in $X$ so that it contains $x$ and $z$.
Let $a\subset A$ be a spherical chamber or a Euclidean alcove containing $z$.
In the spherical case we regard $a$ as the {\em negative chamber} in $A$,
let $\De$ denote the positive chamber $-a$. We have the retraction $f:=
Fold_{a,A}:X\to A$.

\begin{figure}[tbh]
\centerline{\epsfxsize=5in \epsfbox{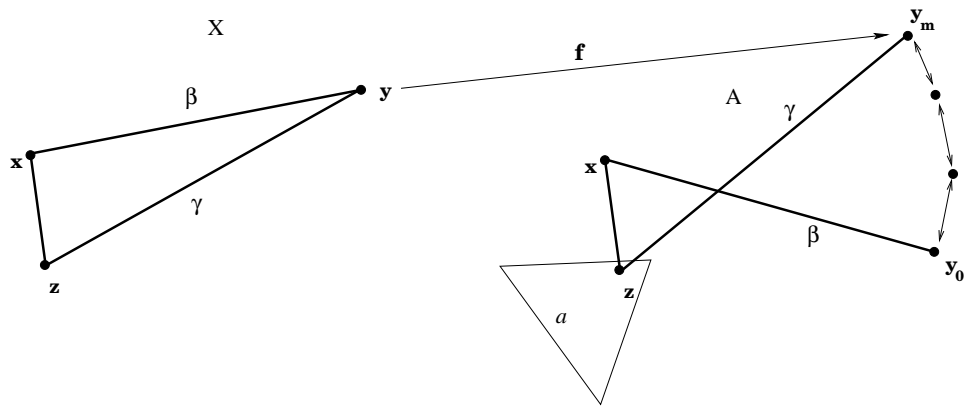}} \caption{\sl
Converting geodesic triangle to a chain.} \label{chain.fig}
\end{figure}

\begin{thm}\label{T2}
There exists a $(A,W,a)$--chain $(y_0,..., y_m)$  such that
$y_m=f(y)$, $d_{ref}(x, y_0)=\be$, $d_{ref}(y_m, z)=\ga$. (In the case when $X$ is a Euclidean building
the above chain is an affine chain.)
See Figure \ref{chain.fig}.
\end{thm}
\proof
We prove the assertion for the spherical buildings as the Euclidean case is completely analogous.
(This is also the only case when this theorem is used in the present paper.)

Our proof is by induction on the rank of the building.
Consider first the case when $rank(X)=0$ (i.e. $A=S^0$ is the
2-point set). If $y$ and $x$ are both distinct from $z$, then
$f(y)\ne z$. This implies that $f(y)=x$ and we take
$$
y_0:= z, y_1:= y, m=1.
$$
In the remaining cases we will use the chain $y_0=f(y)=y_m$.

Suppose now that $rank(X)=r\ge 1$ and the assertion holds for all
(spherical) buildings of rank $r-1$, let's prove it for buildings of
rank $r$.

We let $\tilde{p}: [0, c]\to \ol{xy}$ denote the unit speed parametrization
of $\ol{xy}$ and set $p:= f(\t{p})$. {\em We assume for now that $z\notin \ol{xy}$}.

As in the proof of Lemma \ref{broken}, we ``triangulate'' the geodesic
triangle $\t{T}$ into geodesic triangles $\t{T}_i:= [z, \t{x}_i,
\t{x}_{i+1}]$, where the points $\t{x}_i=\t{p}(t_i)$, $i=1,...,n$, are chosen
so that each triangle $\t{T}_i$ is contained in an apartment
$A_i\subset X$ and the map $f$ restricts to an isometry $f:
\t{T}_i\to f(\t{T}_i)\subset A$. Here $\t{x}_0:=x, \t{x}_{n+1}=y$.
Observe that each side of $\t{T}_i$ has positive length.

Let $\t{S}_i:=\Si_{\t{x}_i}(A_i)$ denote the unit tangent sphere at
$\t{x}_i$. Define $\t{s}_i$ to be the (unique) chamber in $\t{S}_i$
containing all the directions of the geodesic segments from
$\t{x}_i$ to the interior of $a$. This determines the {\em positive
chamber} $\t{\De}_i=-\t{s}_i\subset S_i$.

Set
$$
\t\eta_i:=\t{p}'(t_i), \t\xi_i:= -\eta_i,
$$
and let $\t\zeta_i\in \t{S}_i$ denote the unit tangent vector to
$\ol{\t{x}_i z}$.

Now, applying the retraction $f$ to all this data, we obtain:

1. The folded triangle $T=f(\t{T})$ which has two geodesic sides
$\ol{zx}, \ol{zy_m}$ (where $y_m=f(y)$), and the broken side
represented by the path $p=f(\t{p})$. In particular, $d_{ref}(y_m, z)=d_{ref}(y, z)$ (as required by the theorem).

2. The vertices $x_i=p(t_i)= f(\t{x}_i)$ of the broken geodesic $p$.

3. Unit tangent vectors $\xi_i=df(\t{\xi}_i), \eta_i=df(\t{\eta}_i),
\zeta_i=df(\t{\zeta_i})$ in $\Si_x A$. These vectors are tangent to
the segments $\ol{x_i x_{i-1}}, \ol{x_i x_{i+1}}, \ol{x_i z}$
respectively.

4. The positive chamber $\De_i=df(\t{\De}_i)$ and the negative
chamber $s_i=df(\t{s}_i)$ in the spherical Coxeter complex
$(S_i=\Si_{x_i}(A), W_i=W_{x_i})$. The negative chamber contains the
directions tangent to the geodesic segments from $x_i$ to the
chamber $a\subset A$.

\begin{figure}[tbh]
\centerline{\epsfxsize=5in \epsfbox{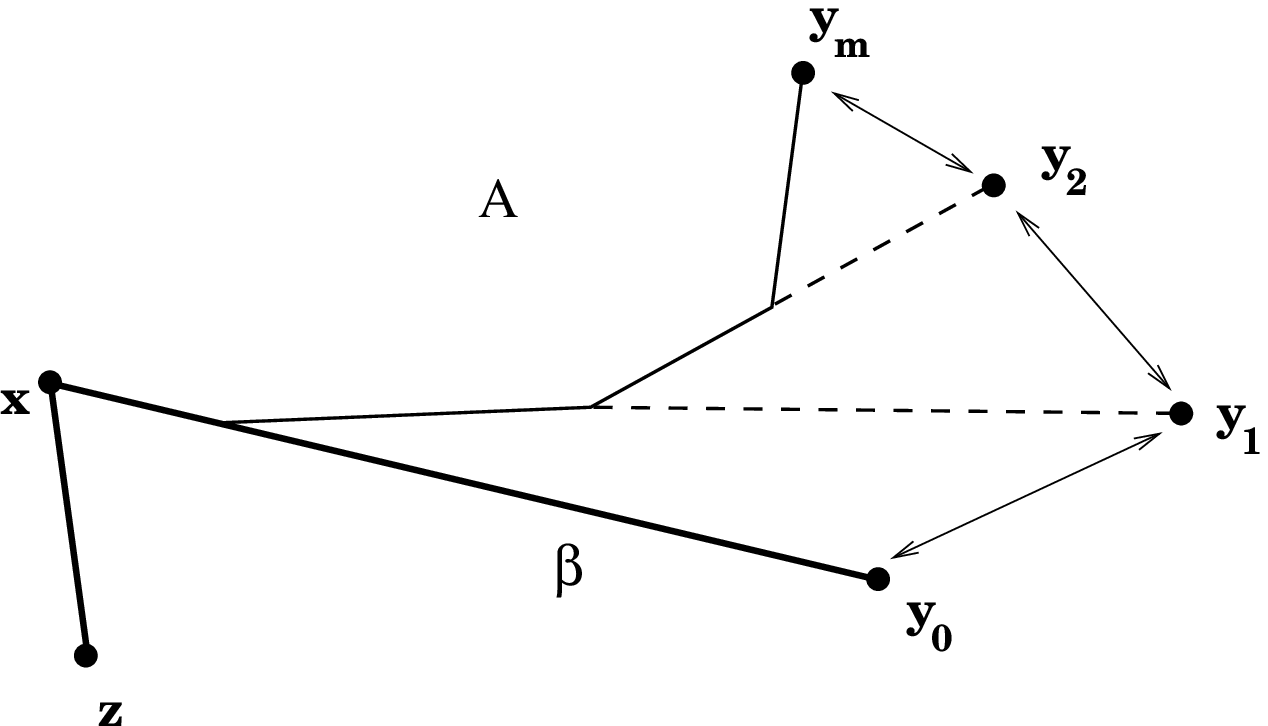}} \caption{\sl
Forming a chain by unbending.} \label{con.fig}
\end{figure}

Our goal is to convert the broken side $p$ of $T$ into a chain in
$A$ by ``unbending'' the broken geodesic $p$ to a geodesic segment
in $A$. See Figure \ref{con.fig}.

\begin{lem}
\label{trivial} The path $p$ satisfies the following:

1. The metric lengths of $p$ and $\tilde{p}=\ol{xy}$ are the same.

2. $p'(0)=\tilde{p}'(0)$.

3. At each break-point $x_i$ there exists an $(S_i, W_{i}, s_i)$-chain from $-\xi_i$ to $\eta_i$.
\end{lem}
\proof The first two assertions are clear from the construction.
Let's prove the last statement. For each $i$ and the point
$\t{v}=\t{x}_i$ we have the spherical building $Y:=\Si_{\t{v}}(X)$
which has rank $r-1$. This building contains the antipodal points
$$
\t\xi_i, \t\eta_i
$$
and the point $\t\zeta_i$. We form the geodesic triangle
$\tau=[\t\xi_i, \t\eta_i, \t\zeta_i]\subset Y$, where we use an
arbitrary shortest geodesic in $Y$ to connect $\xi_i$ to $\eta_i$.
Therefore $\xi_i, \eta_i, \zeta_i$ are vertices of the broken
geodesic triangle $df(\tau)\subset S_i$.

As in Lemma \ref{consistent}, we use the isomorphism $S_i\to
\t{S}_i$ (sending $s_i$ to $\t{s}_i$) to identify these apartments.
Under this identification, $df: Y\to S_i$ is the retraction
$Fold_{s_i, S_i}$ of $Y$ to the apartment $S_i$. Thus $df(\tau)$ is
a folded triangle in $S_i$.

Hence, by the (rank) induction hypothesis, for each $i$ there exists
a chain
$$
(-\xi_i,..., \eta_i)
$$
in the spherical Coxeter complex $(S_i, W_{i}, s_i)$. \qed

\begin{lem}
For each path $p: [0, c]\to A$ satisfying the conclusion of Lemma
\ref{trivial}, there exists a point $y'\in A$ such that
$$
d_{ref}(x, y')= \be %= \Length(p),
$$
and
$$
y' \ge p(c)
$$
in $(A,W,-\De)$.
\end{lem}
\proof We use the second induction, on the number $n$ of vertices in
the broken geodesic $p$. Set $u:= p(c)$.

The metric length of the path $p$ equals the metric length of the
path $\t{p}=\ol{x y}$, the tangent directions of these paths at $x$
are the same. Therefore, if $n=0$ (and hence the path $p$ is
geodesic) there is nothing to prove, one can simply take $y'=u$.

Assume that the assertion holds for all $n\le N-1$, let's prove it
for $N$. We treat the path $p$ as the composition
$$
p|[0, t_N]\cup \ol{x_N u}
$$
Our goal is to replace the geodesic subpath $\ol{x_N u}$ with a
geodesic path $w( \ol{x_N u})$, where $w\in W$ is fixing $x_N$, so
that:

1. $\ol{x_{N-1} x_N}\cup w( \ol{x_N u})$ is a geodesic segment.

2. There exists an $(A, W, -\De)$-chain between $w(u)$ and $u$. Then
we would be done by the induction on $n$. Indeed, the new path
$$
p|[0, t_N]\cup w(\ol{x_N u})
$$
has one less break-point and still satisfies the conclusion of Lemma
\ref{trivial}. Thus, by the induction hypothesis, there exists
$y'\in A$ so that
$$
y' \ge w(u) \ge u \Rightarrow y' \ge u,
$$
$$
d_{ref}(x, y')= d_{ref}(x, y).
$$

{\bf Construction of $w$.} Recall that there exists an $(S_N, W_{N},
-\De_N)$-chain
$$
(-\xi_N=\nu_0,..., \nu_k=\eta_N),
$$
hence we have a sequence of reflections  $r_1,...,r_k\in W_{N}$
(fixing walls $H_i\subset S_N$, $i=1,...,k$) so that:
$$
r_i(\nu_{i-1})=\nu_i, i=1,...,k,
$$
and each wall $H_i$ separates $\nu_{i}$ from the negative chamber
$s_N$. We extend each reflection $r_i$ from $S_N$ to a reflection
$r_i$ in $A$, and each $H_i$ to a wall $H_i$ in $A$.

We therefore define the following points in $A$:
$$
y_k:=u, y_{k-1}:= r_{k}(y_k), y_{k-2}:= r_{k-1}(y_{k-1}),\ldots,
y_{0}:= r_{1}(y_1).
$$
Note that the directions $\nu_i$ are tangent to the segments
$\ol{x_N y_i}$. Thus for each $i$, the wall $H_i$ separates the
point $y_{i}$ from the negative chamber $a\subset A$ and the
sequence
$$
(y_0,...,y_k=u)
$$
forms a chain. We set $w=r_1\circ ... \circ r_k$. The vector $\nu_0$
is antipodal to $-\xi_N$, hence the path
$$
p|[0, T_N] \cup w(\ol{x_N u})
$$
is geodesic at the point $x_N=p(T_N)$. \qed

This concludes the proof of Theorem in the case when $z\notin\ol{x
y}$.

We now consider the special case when the above proof has to be
modified: The triangle $T$ is degenerate, i.e. $z\in \ol{xy}$, but
the alcove $a$ is such that $y\notin A$.  Thus the folding
$Fold_{a,A}$ is not an isometry on $T$. Then the tangent direction
$\t\zeta_i$ is not defined when $x_i=z$. Note that $x_i=z$ then is
the only break-point in the broken side of $T'$.

In this case we replace the vertex $z$ with an arbitrary point $z'$
in the interior of $a$ and repeat the above arguments. The chains
constructed in the process will be independent of the choice of $z'$
and thus, after taking the limit $z'\to z$, we obtain a chain as
required by the assertion of Theorem. \qed

\begin{cor}
Cf. \cite[Theorem 8.2, Part 4]{KLM3}.  Suppose that $X$ is a Euclidean building.
Assume that $\al:=d_\De(z,x), \be:=d_\De(x,y), \ga:=d_\De(y,z)$ are in $P(R^\vee)$
and $x, y, z$ are special vertices of $X$. Then
$$
\al+\be+\ga \in Q(R^\vee).
$$
\end{cor}
\proof
Let $(y_0,...,y_m)$ be an affine chain given by Theorem \ref{T2}.
We regard the point $x$ as the origin $o$ in $A$; thus we will regard $z$, $y_0, y_m$ as vectors in $V$.
Then, according to Lemma \ref{difference},
$$
y_m-y_0\in Q(R^\vee).
$$
Consider the vectors $\be':= y_0-x, \ga':= z-y_m, \al':= x-z$ in $P(R^\vee)$. By the definition of $\De$-length,
$$
\al'\in W_{sph} \al, \quad \be'\in W_{sph} \be, \quad \ga'\in W_{sph} \ga.$$
 Therefore, by applying  Lemma \ref{difference} again we see that
the differences
$$
\al-\al',  \quad \be-\be',  \quad \ga-\ga'
$$
all belong to $Q(R^\vee)$. Since
$$
\al'+\be'+\ga'=y_0-y_m\in Q(R^\vee),
$$
the assertion of lemma follows. \qed

\medskip
The following simple proposition establishes a partial converse to
Theorem \ref{T2}:

\begin{prop}
\label{simplechain} Suppose that $X$ is a thick spherical building and, as before, the point $z$ belongs to a
negative chamber $a=-\De$. Then, for each simple chain $(y_0, y_1)$
such that $d_{ref}(x, y_0)=\pi, d_{ref}(z, y_1)=\be$, there exists a
point $y\in X$ so that
$$
d_{ref}(y, z)= \ga \quad \hbox{and}\quad  d_{ref}(x , y)= \pi.$$
\end{prop}

\begin{rem}
Recall that $d_{ref}(x , y)= \pi$ means that the points $x$ and $y$ are antipodal.
\end{rem}

\proof Let $\tau(y_0)=y_1$ where $\tau$ is a reflection in a wall
$H\subset A$ as in the definition of a chain. Let $A=A^-\cup A^+$ be
the union of half-apartments, where $A^-$ is bounded by $H$ and contains $a$.
By the definition of a chain, $y_1\in A^-, y_0\in A^+$
and hence the antipodal point $x=-y_0$ belongs to $A^+$.

Since $X$ is thick, there exists a half-apartment $B^-\subset X$
which intersects $A$ along $H$. Define the apartment $B:= A^-\cup
B^-$; then there exists  an isomorphism of Coxeter complexes
$$
\phi: A\to B, \phi|A^-=id.
$$
We set $y:= \phi(y_0)$.

Since $\phi$ is an isomorphism of Coxeter complexes which fixes $z$,
it preserves the refined distance to the point $z$ and hence
$$
d_{ref}(y, z)= d_{ref}(y_1, z)=\ga.$$

The union $C:= A^+\cup B^-$ is also an apartment in $X$. Then there
exists an isomorphism  $\psi: B^-\to A^-$ so that
$$
\psi\circ\phi|A^+= \tau|A^+.
$$
The isomorphism $\psi$ extends to an isomorphism $\rho: C\to A$
fixing $A^+$ pointwise and hence fixing the point $x$. Therefore
$$
d_{ref}(x , y)= d_{ref}(x, y_0)= \pi.  \qed
$$

\subsection{Folding polygons in Euclidean buildings}
\label{euclideanfolding}

Our next goal is to show that each folding transforms certain piecewise-linear
paths in Euclidean buildings to paths satisfying the chain condition.

Suppose that $X$ is a Euclidean building with model apartment
$(A,W_{aff})$, $\De\subset A$ is the positive Weyl chamber with the
tip $o$. Consider a piecewise-linear path $\tilde{p}: [0,c]\to \t{A}$, which is
parameterized  with the unit speed, where $\tilde{A}\subset X$ is an
apartment. We assume that for each $t\in [0,c]$
$$
\tilde{p}'_-(t) \sim_{W_{sph}} \tilde{p}'_+(t),
$$
for instance, $\tilde{p}$ could be a geodesic path.

Thus the path $\tilde{p}$ {\em trivially satisfies the chain
condition}. Let $g: X\to \De$ be a folding into $\De$,
$g=Fold_{z,h,\De}$ for a certain $z\in A$ and $h$. Recall that the
folding $g$ is the composition of three maps:
$$
g= \P_{\Del} \circ h \circ f, \quad f=Fold_{a,A},
$$
where $a$ is an alcove in $A$ containing $z$, $h\in Dil(A,W_{aff})$ is a
dilation sending $z$ to the point $o$. Consider
the structure of a Coxeter complex on $A$ given by the pull-back
$$
h^*(A,W_{aff}).
$$
We thus get a new (typically non-thick) building structure for $X$,
the one modeled on $h^*(A,W_{aff})$.

\begin{defn}
We say that a path $\tilde{p}$ is {\em generic} if it is disjoint
from $z$ and from the codimension 2 skeleton of $X$ and the
break-points of $\tilde{p}$ are disjoint from the codimension 1
skeleton of $X$, where $X$ is regarded as a building modeled on
$h^*(A,W_{aff})$.
\end{defn}

The main result of this section is

\begin{thm}
\label{chaincondition} The folded path $p= g(\tilde{p})$ satisfies
the chain condition.
\end{thm}
\proof The proof of this theorem is mostly similar (except for the
projection $\P$ which causes extra complications) to the proof of
Theorem \ref{T2} in the previous section. We will prove Theorem
\ref{chaincondition} in two steps: We first establish it for the
paths $\tilde{p}$ which are {\em generic}. Then we use the {\em
compactness theorem} to prove it in general.

\begin{prop}
\label{chaincondition0} The conclusion of Theorem
\ref{chaincondition} holds for {\em generic} paths $\tilde{p}$.
\end{prop}
\proof If a point $\tilde{x}=\tilde{p}(t)$ is a regular point of
$X$, then
$$
dg_{\tilde{x}}: \Si_{\tilde{x}}(X)\to \Si_x(A), \quad x=g(\tilde{x})
$$
is an isometry. Thus the path $p$ trivially satisfies the chain
condition at the point $x$.

Therefore we assume that $\tilde{x}$ is a singular point of $X$. Since
$\tilde{p}$ is assumed to be generic, this point lies on exactly one
wall of $X$; moreover, $\tilde{p}$ is geodesic near $\tilde{x}$.

We first analyze what happens to the germ of $\tilde{p}$ at
$\tilde{x}$ under the retraction $f$. We suppose that the
restriction of $f$ to the germ $(\t{p}, \t{x})$ is not an isometry
(otherwise there is nothing to discuss).  Let $\tilde\zeta\in
\Si_{\tilde{x}}(X)$ denote the tangent to the geodesic
segment $\ol{\tilde{x}z}$. Let $\eta\in \Si_{\tilde{x}}(X\cap
\tilde{A})$ be the tangent vector $p'(t)$, $\tilde\xi:= -\tilde\eta$
(this vector is also tangent to the path $p$). Set $\be:=
d_{ref}(\tilde\zeta, \tilde\eta)$.

We obtain the  triangle $\tau=[\tilde\xi, \tilde\eta,
\tilde\zeta]$ in $\Si_{\tilde{x}}(X)$. The derivative of the
retraction $f$ at $\t{x}$ is a retraction of the spherical building
$\Si_{\t{x}}(X)$ into its apartment $S$, after identification of $S$
with the sphere $S_{x'}(A)$, $x':=f(\t{x})$ (see Lemma
\ref{consistent}). Define the following elements of $S$:
$$
\eta':=df_{\tilde{x}}(\tilde\eta), \quad \xi':=
df_{\tilde{x}}(\tilde\xi), \quad
\zeta':=df_{\tilde{x}}(\tilde\zeta).
$$

Therefore, according to Theorem \ref{T2}, the folded triangle
$\tau'=df_x(\tau)\subset S$ yields an $(S,W_{x'},
-\De_{x'})$-chain\footnote{Which is necessarily simple since
$\tilde{p}$ is assumed generic.}
$$
(\mu'=-\xi', \eta')
$$
such that $d_{ref}(\zeta, \eta')= d_{ref}(\t\zeta, \t\eta)=\be$.

Here $-\De_{x'}$ is a chamber in $(S,W_{x'})$ which contains
the unit tangent vector to the segment $\ol{x' z'}$ where $z'\in a$ is a regular point.

\begin{rem}
Note that our assumptions on $\tilde{p}$ imply that $(S,W_{x'})$
has a unique wall. If the corresponding wall in $A$ does not
pass through $z$, then the negative chamber $-\De_{x'}$  in $(S,W_{x'})$ is
uniquely determined by the condition that it contains the direction
tangent to
 $\ol{ x' z}$.
\end{rem}

Consider now the effect of the rest of the folding $g$ on the path
$\tilde{p}$ at $\tilde{x}$. Let $x:= g(\tilde{x})$. We identify $S$
with the unit tangent sphere at the point $x$.

The dilation $h$ clearly preserves the chain condition at $x'$
(since it acts trivially on the unit tangent sphere). The
restriction of the projection $\P=\P_{\Del}$ to the germ of
$hf(\tilde{p})$ at $hf(\tilde{x})$ is necessarily an isometry (since
$\tilde{p}$ is generic), hence it is given by an element $w\in
W=W_{sph}$. This element transforms the above chain to another
$(S,W_x, -\De_x)$--chain, where
$$
\De_x:= d(w\circ h) (\De_{x'}).
$$
What is left to verify is that the positive chamber $\De_x$ in this
complex contains a translate of the positive chamber $\De$. In case
when $x$ belongs to the interior of $\De$, the segment $\ol{ox}$ is
not contained in any wall and thus the negative chamber $-\De_x$ has
to contain the initial direction of the segment $\ol{xz}$ (see the
remark above). However this initial direction belongs to $-\Del$ and
thus $\De_x$ contains $\De$.

Consider the exceptional case when $x$ is on the boundary of $\De$.
It then belongs to a unique wall $H$ in the Coxeter complex
$(A,W_{aff})$ and this wall passes through the origin $o$. Rather
than trying to use Theorem \ref{T2} to verify the chain condition at
$x$, we give a direct argument. Let $\eta, \xi$ be the unit vectors
which are the images of $\eta', \xi'$ under
$$
d(w\circ h): \Si_{x'}(A)\to \Si_x(A).
$$
Since the path $p$ is entirely contained in $\De$, the vector
$p'_-(t)$ points {\em outside} of $\De$ and the vector $p'_+(t)$
points {\em inside}. The reflection $\si$ in the wall $H$ sends the
vector $-\xi=p'_-(t)$ the vector $\eta=p'_+(t)$. It is then clear
that the (simple) chain condition is satisfied at the point $x$.

\medskip
Lastly, we consider the points $\tilde{x}=\tilde{p}(t)$ for which
$f$ is an isometry on the germ of $\tilde{p}$ at $\tilde{x}$.
The point $x=g(\tilde{x})$ belongs to
a face of $\De$ contained in a wall $H$, and this is the only wall
of $(A,W_{aff})$ which passes through $x$. Then,
necessarily, the germ of the path $hf(\tilde{p})$ at
$hf(\tilde{x})$ is a geodesic.  We now simply repeat the
arguments of the exceptional case in the above proof (see also the
proof of Proposition \ref{chaincondition1}) to see that $p$
satisfies the chain condition at $x$. \qed

\medskip
We are now ready to prove Theorem \ref{chaincondition} for arbitrary
paths $\tilde{p}$. We will do so by approximating the path
$\tilde{p}$ via {\em generic paths}. Let $\la$ be an arbitrary
vector in $\t{A}$. We let $\tilde{q}_\la:=\tilde{p}+\la$ denote the
translation of the path $\tilde{p}$ by the vector $\la$. It is
clear, from the dimension count, that for an open and dense set of
vectors $\la$, the path $\tilde{q}_\la$ is {\em generic}.

Since the folding $g$ is continuous,
$$
p=g(\tilde{p})= \lim_{\la\to 0} g(\tilde{q}_\la).
$$
By the Proposition \ref{chaincondition0}, each $g(\tilde{q}_\la)$
satisfies the chain condition. Observe that the $\De$-lengths of the
paths $\tilde{p}+\la$ are independent of $\la$. Since $f$ and $\P$
preserve $\De$-lengths of piecewise-linear paths and the dilation $h$ changes them
by a fixed amount, we can apply the compactness theorem (Theorem
\ref{compactness}) to conclude that the limiting path $p$ satisfies
the chain condition as well. \qed

\medskip
We now verify that, at certain points, the folded path $p$ satisfies
the {\em maximal chain condition}.

\begin{prop}
\label{chaincondition1} Under the assumptions of Theorem
\ref{chaincondition} let $\tilde{x}=\tilde{p}(t)$ be such that the
folding $f$ restricts to an isometry on the germ $(\tilde{p},
\tilde{x})$. Then the path $p=g(\tilde{p})$ satisfies the {\em
maximal chain condition} at $x=g(\tilde{x})$.
\end{prop}
\proof Our proof follows Littelmann's arguments in his proof of the
PRV Conjecture  \cite{Littelmann1}. We fold the path
$q:=hf(\tilde{p})$ into $\De$ inductively.

We subdivide the interval $[0,c]$ as
$$
0=t_0 < t_1<...< t_k=c
$$
such that $[t_i, t_{i+1}]$ are maximal subintervals so that $q|[t_i,
t_{i+1}]$ is contained in a Weyl chamber of $W_{sph}$.

We first apply to $q$ an element $w_0\in W_{sph}$ which sends $q([0,
t_1])$ into $\De$, so we can assume that this subpath belongs to
$\De$. Assume that the restriction of $f$ to the germ $(\t{p}, \t{p}(t_1))$ is an isometry.
Let $\mu', \eta'$ be the vectors $q'_-(t_1), q'_+(t_1)$. Then
$$
\mu'\sim \eta',
$$
see Lemma \ref{positivity}. Set $x:= q(t_1)$. The image $\eta$ of
the vector $\eta'$ under $\P$ is obtained as
$$
dw_1(\eta'),
$$
where $w_1\in W_{sph}$ fixes the point $x$ and $\eta$ is the unique
vector in the $W_x$-orbit of $\eta'\in S_{x}$ which points inside
$\De$. Below we describe $w_1$ as a composition of reflections.

Let $R'$ denote the root subsystem in $R$ generated by the set of
simple roots $\Phi'$ which vanish at the point $x$. Let $\De'$
denote the positive chamber for $W_x$ defined via $\Phi'$. Then the
vector $\eta$ can be described as the unique vector in the
$W_x$-orbit of $\eta'$ (now, regarded as a vector in $V=T_o(A)$)
which belongs to the interior of $\De'$. According to Lemma
\ref{simplesubsystem},
$$
w_1= \tau_m\circ .. \circ \tau_1, \hbox{~where for each~} i, \quad
\tau_i=\tau_{\be_i}, \be_i\in \Phi',
$$
so that the sequence of vectors
$$
(\eta_0=\eta', \eta_1:=
\tau_1(\eta_{0}),...,\eta_m=\tau_m(\eta_{m-1})=\eta),$$ is a chain
in $(S, W_x,-\De)$ which is maximal as a chain in $(S,W_{sph},-\De)$.

We therefore apply the identity transformation to the path $q|[0,
t_1]$ and the element $w_1$ to the path $q|[t_1, c]$ to transform
the path $q$ to the new path
$$
q_1= q|[0, t_1] \cup w_1\circ q|[t_1, c].
$$
Clearly, $\P(q)= \P(q_1)$ and $\eta$ is the unit vector tangent to
$q_1|[t_1,c]$ at $x$. The above arguments therefore show that $q_1$
satisfies the {\em maximal chain condition} at the point $x$.

We then proceed to the next point $t_2$, $q_1(t_2)$ belongs to the
boundary of $\De$ and we transform $q_1$ to $q_2$ by
$$
q_2|[0,t_2]=q_1|[0,t_2], \quad q_2|[t_2, c]= w_2\circ q_1|[t_2,c],
$$
where $w_2$ is a certain element of $W_{sph}$ fixing
$q_1(t_2)=q_2(t_2)$. Therefore $\P(q_2)=\P(q_1)=\P(q)$ and we repeat the
above argument. \qed

\begin{defn}
Suppose that $P=\ol{zx}\cup p \cup \ol{yz}$ is a polygon in $A$,
where $p:[0,1]\to A$ is a piecewise-linear path such that $p(0)=x, p(1)=y$. We say
that $P$ {\em satisfies the chain condition} (resp. simple chain
condition, resp. maximal chain condition) if its subpath $p$
satisfies the chain condition (resp. simple chain condition, resp.
maximal chain condition).
\end{defn}

Therefore, as an application of Theorem \ref{chaincondition} we
obtain

\begin{cor}
\label{folded-chain}
 Suppose that $T=[\t{z}, \t{x}, \t{y}]\subset X$ is a
geodesic triangle, $\t{z}$ is a special vertex which belongs to an
alcove $a\subset A$. Let $\De\subset A$ be a Weyl chamber with the
tip $\t{z}=o$ and $P=Fold_{\De}(T)$ be the folding of $T$ into $\De$. Then
the folded triangle $P$ satisfies the chain condition.
\end{cor}

A converse to this corollary will be proven in Theorem
\ref{characterization}; the following is a {\em partial} converse to
Corollary \ref{folded-chain} (which is essentially contained in
\cite[Lemma 7.7]{KLM3}):

\begin{cor}
\label{chain-folded0} Let $\De\subset A$ be a Weyl chamber with  tip
$o$ in $X$. Suppose that a polygon $P= [o, x_1,...,x_{n}]\subset
\De$, satisfies the {\em simple chain} condition (at each vertex
$x_i, 0<i<n$) and
$$
p=\ol{x_1 x_2}\cup ... \cup \ol{x_{n-1} x_n}
$$
is a {\em billiard} path. Then $P$ unfolds to a geodesic triangle
$T\subset X$, i.e. $Fold_{\De}(T)=P$.
\end{cor}
\proof Let $f:=Fold_{\De}$.  We run the argument from the proof of
Theorem \ref{T2} in the reverse; the reader will observe that our
argument is essentially the same as in the proof of the Transfer
Theorem in \cite{KLM2}. Triangulating the polygon $P$ from the
vertex $o$ we obtain geodesic triangles $P_i=[o, x_i, x_{i+1}]$, $i=1, 2,...,n-1$.
Let $\xi_i, \zeta_i, \eta_i \in \Si_{x_i}(A)$ denote the unit tangent vectors to the
segments $\ol{x_i x_{i-1}}, \ol{x_i o}, \ol{x_i x_{i+1}}$
respectively.

We unfold $P$ inductively. Set $T_1:= P_1$; let $A_1:= A$,
this apartment contains the triangle $T_1$. Set $\t{x}_1:= x_1, \t{x}_2:= x_2$,

Suppose that we have constructed apartments $A_i\subset X$ and flat triangles
$T_i=[o, \t{x}_i, \t{x}_{i+1}]\subset A_i$, $i=1,...,m-1$,
so that $T_i$ is congruent to $P_i$ ($i=1,...,m-1$) and
$$
\angle( \t{\xi}_i, \t{\eta}_i)=\pi, i=1,...,m-1.
$$
Here $\t{\xi}_i, \t{\eta}_i, \t{\zeta}_i$
are directions in $\Si_{\t{x}_i}(X)$ which correspond to the directions $\xi_i, \eta_i, \zeta_i$ under the
congruences $T_i\to P_i$. Our goal is to produce a flat triangle $T_m\subset A_m\subset X$ so that the above
properties still hold.

\begin{figure}[tbh]
\centerline{\epsfxsize=5in \epsfbox{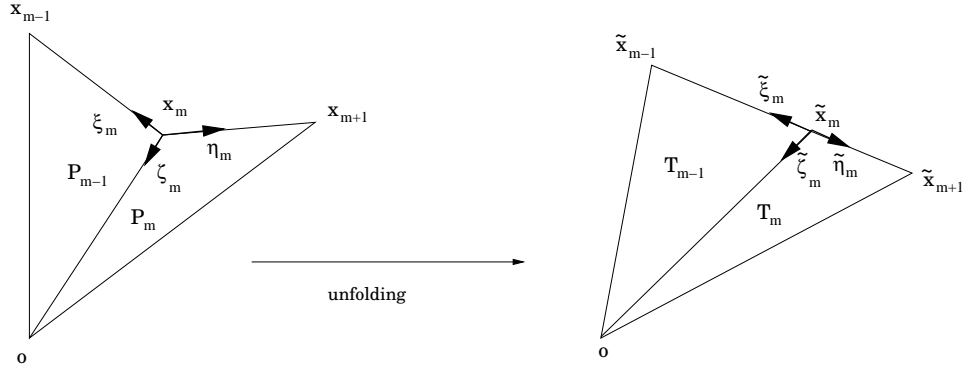}} \caption{\sl
Unfolding a broken triangle.} \label{unfold.fig}
\end{figure}

Since we have a simple chain $(-\xi_m, \eta_m)$ in
$(S_{x_m}, W_{x_m})$, it follows from Proposition \ref{simplechain}
that there exists a point $\tilde{\eta}_m\in \Si_{\t{x}_m}(X)$ so that
$$
d(\tilde{\eta}_m, \t{\xi}_m)=\pi, \quad d(\tilde{\eta}_m,
\zeta_m)=d_{ref}(\eta_m, \zeta_m).
$$
Let $A_m$ denote an apartment in $X$ which contains $\ol{o \t{x}_m}$ and
such that $\tilde{\eta}_m$ is tangent to $A_m$. Construct a geodesic segment
$\ol{\t{x}_m \t{x}_{m+1}}\subset A_m\subset X$ whose metric length equals the one of $\ol{x_m x_{m+1}}$
and whose initial direction is $\tilde{\eta}_m$. This defines a flat triangle
$$
T_m=[o, \t{x}_m, \t{x}_{m+1}]\subset A_m.
$$
It is clear from the construction that the triangle $T_m$ is congruent to $P_m$, in particular,
$d_{ref}(o, x_3)=d_{ref}(o, y_3)$. Observe also that
$$
\ol{\t{x}_{m-1} \t{x}_m}\cup \ol{\t{x}_m \t{x}_{m+1}}$$
is a geodesic segment (because $\angle(\tilde{\eta}_m, \t{\xi}_m)=\pi$). See Figure \ref{unfold.fig}.

Therefore, by induction we obtain a geodesic triangle $T=[o, x_1, \t{x}_n]\subset X$, which
is triangulated (from $o$) into flat geodesic triangles $T_i$ which are congruent  to $P_i$'s.
We claim that $f(T)=P$. For each $i$ the folding $f$ sends the
triangle $T_i$ to $P_i$, according to Proposition \ref{uniqueness}.
Therefore $f(T)=P$.

\medskip
As in the proof of Theorem \ref{T2}, the argument has to be modified
in case when $x_i=o$ for some $i$, $x_i=p(t_i)$. Then the vector
$\mu=p'_-(t_i)$ belongs to the negative chamber $-\De$ and the
vector $\la=p'_+(t_i)$ belongs to the positive chamber $\De$. Since
$p$ is a billiard path, there exists $w\in W_{sph}$ which sends
$\mu$ to $\la$. Now the chain and billiard conditions imply that
$x_i$ is the  only break-point in $p$. Thus we can take
$$
T:= \ol{ox_1} \cup p([0,t_i]) \cup w^{-1} p([t_i, 1]) \cup
w^{-1}(\ol{ox_n}).
$$
This degenerate geodesic triangle (it is contained in the geodesic
through the points $x_1, w(x_n)$) folds to $P$ under the projection
$\P_\De: A\to \De$. \qed

\medskip
The same argument as above proves the following generalization of
Corollary \ref{chain-folded0}

\begin{cor}
\label{chain-folded} Suppose that $P$ is a  polygon in $\De$ which
is the composition
$$
\ol{ox} \cup   p \cup  q \cup  \ol{yo}.
$$
Assume that paths $p, q$ satisfy the simple chain condition. Then
there exists a polygon $\tilde{P}\subset X$ of the form
$$
\ol{ox} \cup \tilde{p} \cup  \tilde{q} \cup \ol{\t{y}o}
$$
so that $f(\tilde{P})=P$, $f(\tilde{p})=p, f(\tilde{q})=q$,
$f(\t{y})=y$, and $\tilde{p}, \tilde{q}$ are geodesic paths.
\end{cor}

 We now use our analysis of the folded triangles
(polygons) to relate them to the Littelmann triangles (polygons).

\section{Littelmann polygons}
\label{LSpolygons}

\subsection{LS paths}
\label{LSpaths}

 Let $R$ be a root system on a Euclidean vector space
$V$, $W=W_{sph}$ be the finite Coxeter group associated with $R$,
let $W_{aff}$ denote the affine Coxeter group associated to $R$.
This root system $R$ is actually the {\em coroot system} for the one
considered by Littelmann in \cite{Littelmann2}. Accordingly, we will
switch weights to coweights, etc. We pick a Weyl chamber $\Del$ for $W$, this
determines the positive
roots and the simple roots in $R$. We get the Euclidean Coxeter
complex $(A,W_{aff})$, where $A$ is the affine space corresponding
to $V$. Given $x\in A$ let $W_x$ denote the stabilizer of $x$ in $W_{aff}$.

\medskip
Suppose we are given a  vector $\la\in  \De\subset V$, a
sequence of real numbers
$$
\ul{a}=( a_0=0< a_1 < ... < a_r=1),
$$
and a sequence of vectors in $W\la$
$$
\ul\nu= (\nu_1, ... , \nu_r), \hbox{~~so that~~} \nu_1 > ... > \nu_r
$$
with respect to the order in Definition \ref{order}.

\begin{defn}
The pair $(\ul\nu, \ul{a})$ is called a {\em real (billiard) path}
of the $\De$-length $\la$.
\end{defn}

\begin{defn}
(P. Littelmann \cite{Littelmann2}.) A real path of $\De$-length
$\la$ is called {\em rational} if  $\la$ is a coweight and all
numbers $a_i$ a rational.
\end{defn}

\begin{rem}
Littelmann uses the notion {\em path of type $\la$} rather than of
the $\De$-length $\la$.
\end{rem}

Set $a_i':= a_i -a_{i-1}$, $i=1, 2,...,r$. The data $(\ul\nu,
\ul{a})$ determines a piecewise-linear path $p\in {\mathcal P}$ whose restriction
to each interval $[a_{i-1}, a_i]$ is given by
\begin{equation}
\label{path}
 p(t)= \sum_{k=1}^{i-1} a_k' \nu_k + (t-a_{i-1})\nu_i,
t\in [a_{i-1}, a_i].
\end{equation}

Our interpretation of real and rational paths is the one of a broken
(oriented) geodesic $L$ in $V$. Each oriented geodesic subsegment of
$L$ is parallel to a positive multiple of an element of $W\la$, thus
$L$ is a billiard path. The break points of the above path are the
points
$$
x_1=a_1\nu_1, \ldots, x_i=x_{i-1}+ a_i' \nu_{i}, \ldots
$$
Since $\sum_i a_i'=1$, is clear that
$$
\Length (L)=\la,
$$
in the sense of the definition in section \ref{distances}. This
justifies our usage of the name $\De$-length $\la$ in the above
definitions, rather than Littelmann's notion of {\em type}.

Observe that given a piecewise-linear path $p(t)\in {\mathcal P}$ (parameterized
with the constant speed) one can recover the nonzero vectors
$\nu_i\in V$ and the numbers $a_i$ and $a_i'$.

\begin{defn}
\label{LSdef} (P. Littelmann \cite{Littelmann2}.)
A rational path
$p(t)$ is called an {\em LS path}\footnote{a Lakshmibai-Seshadri
path} if it satisfies further {\em integrality condition}:

For each $i=1,...,s-1$ there exists an $a_i$-chain for the pair
$(\nu_i, \nu_{i+1})$ (in the sense of Definition \ref{achain}).
\end{defn}

Observe that, since
$$
\eta \ge \tau \iff -\tau \ge -\eta,
$$
it follows that $p$ is an LS path if and only if $p^*$ is.

\begin{thm}
\label{locint}
 (P.\ Littelmann \cite[Lemma 4.5]{Littelmann2}) Each LS path belongs to
${\mathcal P}_{\Z,loc}$.
\end{thm}

Our next goal is to give a more geometric interpretation of LS
paths. Suppose that $p\in {\mathcal P}$ is a billiard path given by
the equation (\ref{path}), with the vertices
$$
0=x_0, x_1,...,x_r.
$$
At each vertex point $x_i, 0<i<r$,  we have unit tangent vectors
$\xi_i, \mu_i$ which are tangent to the segments $\ol{x_i x_{i-1}},
\ol{x_i  x_{i+1}}$. Note that at each vertex $x_i, 0<i<r$ we have
the restricted and unrestricted spherical Coxeter complexes; the
positive chamber $\De$ in $V$ determines positive chambers $\De_i$
in the restricted spherical complexes $(S_{x_i}, W_{x_i})$.

\begin{thm}
\label{T3} A billiard path $p(t)$ of $\De$--length $\la\in P(R^\vee)$
is an LS path if and only if it is a Hecke path which satisfies the
{\em maximal chain condition} (cf. Definition \ref{D3}): At each
vertex $x_i, 0<i<r$ there exists a $(S_{x_i}, W_{x_i},-\De_i)$-chain
between $-\xi_i$ and $\mu_i$ and this chain is maximal as a
$(S_{x_i}, W)$-chain.
\end{thm}
\proof Recall that given a nonzero vector $v\in V$, $\bar{v}$
denotes its normalization $v/|v|$.

\begin{figure}[tbh]
\centerline{\epsfxsize=4.5in \epsfbox{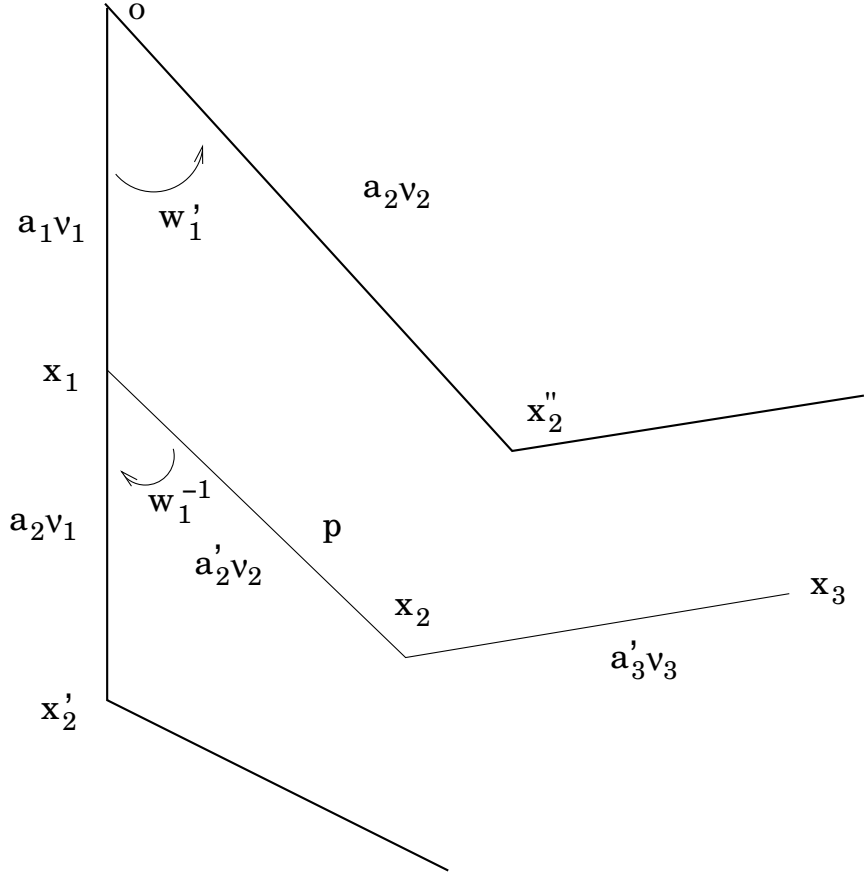}} \caption{\sl
Unbending a path.} \label{LS.fig}
\end{figure}

It is easy to see (and left to the reader) that if $p(t)$ is a
satisfies the above chain condition and $\la=\Length(p)$ is a
coweight, then all numbers $a_i$ are rational.

Consider the first break point $x_1=x_0+a_1\nu_1$ of the broken
geodesic path $p(t)$. Observe that
$$
\bar\nu_1=-\xi_1, \bar\nu_2=\mu_1\in S_{x_1}.
$$
According to Proposition \ref{eq}, existence of an $a_1$-chain for
the pair $(\nu_1, \nu_{2})$ is equivalent to existence of an
$(S_{x_1}, W_{x_1}, \De_1)$-chain, which is maximal in the
unrestricted Coxeter complex,
$$
(\bar\nu_1=\eta_{1,0}, \eta_{1,1},...,\eta_{1,s_1}=\bar\nu_2).
$$
Thus the path $p$ satisfies the maximal chain  condition at the
first break point $x_1$ if and only if it satisfies the {\em
integrality condition} as in Definition \ref{LSdef}, at the point
$x_1$.

We now proceed to the next break point $x_2=x_1+ a_2' \nu_2$. We
identify normalized vectors $\bar\nu_1, \bar\nu_2$ with unit vectors
in $S_{x_1}$. Note that if $p(t)$ is an LS path of the $\De$-length
$\la$, then there exists an element $w_1\in W_{x_1}$ which sends
$\bar\nu_1$ to $\bar\nu_2$. The same is true if $p$ is a Hecke path.

Set $x_2':= w_1^{-1}(x_2)$. Observe that, in both cases of an LS
path and a Hecke path,
$$
\ol{x_0 x_1}\cup \ol{x_1 x_2'}
$$
is a geodesic segment $\ol{x_0 x_2'}$; the corresponding directed
segment represents the vector $a_2 \nu_1$. Let $w_1'\in W_{sph}$
denote the linear part of $w_1$. Set $x_2'':=w_1'(x_2')$. We
translate the vectors $\bar\nu_2, \bar\nu_3$ to the unit tangent
sphere $S_{x_2''}$.  The directed segment $\ora{x_0 x_2''}$
represents the vector $a_2 \nu_2$. See Figure \ref{LS.fig}.

We are now again in position to apply Proposition \ref{eq} with
$a=a_2$: There exists a maximal chain
$$
(\bar\nu_2=\eta_{2,0}, \eta_{2,1},...,\eta_{2,s_2}=\bar\nu_3)
$$
if and only if there exists an $a_2$-chain for the pair
$(\nu_2,\nu_3)$. The product $w_1\circ (w_1')^{-1}$ is a translation
in $W_{aff}$ which carries $x_2''$ back to $x_2$. Therefore it
induces an isomorphism of the restricted Coxeter spherical complexes
$$
(S_{x_2''}, W_{x_2''})\to (S_{x_2}, W_{x_2})
$$
which carries positive chamber to positive chamber. Hence this
translation sends the chain $(\eta_{2,i})$ to a maximal chain in
$(S_{x_2}, W_{x_2})$.

We continue in this fashion: On the $i$-th step we ``unbend'' the
broken geodesic
$$
\ol{x_0 x_1}\cup ... \cup\ol{x_{i-1} x_i}
$$
to a directed geodesic segment $\ora{x_0 x_i'}$ representing the
vector $a_i\nu_1$, then apply an appropriate element $w_{i-1}'\in
W_{sph}$ to transform segment $\ol{x_0 x_i'}$  to $\ol{x_0 x_i''}$;
finally, appeal to Proposition \ref{eq} to establish equivalence
between the maximal chain condition and the LS path axioms. \qed

As a corollary of Theorem \ref{T3} we obtain

\begin{cor}
\label{C3} Let $T=[z,x,y]\subset X$ be a geodesic triangle and
$f=Fold_{z,h,\De}$ be a folding into the Weyl chamber. Set
$\be:=d_\De (x,y)$. Assume that $T'=f(T)$ is such that $f(x), f(y)$
and all break-points of the broken geodesic $f(\ol{xy})$ are special
vertices. Then $f(\ol{xy})$ is an LS path of the $\De$-length $k\be$.
Here $k$ is the conformal factor of the dilation $h$.
\end{cor}

\subsection{Root operators}
\label{operators}

With each simple root $\al\in \Phi$, Littelmann
\cite{Littelmann2} associates {\em raising and lowering} root
operators $e_\al$ and $f_\al$ acting ${\mathcal P}$ as follows.

Recall that given a path $p(t)$ and a root $\al$ we have the height function
$h_\al(t):=\al(p(t))$. The number $m_\al$ is the minimal value of
$h_\al$ on $[0,1]$.

If $m_\al>-1$ then $e_\al$ is not defined on $p$. Otherwise let
$t_1$ be the minimal $t$ for which $h_\al(t)=m_\al$ and let $t_0\in
[0, t_1]$ be maximal such that $h_\al(t)\ge m_\al+1$ for all $t\in
[0, t_0]$.

The operator $e_\al$ will not change the path $p$ for $t\in [0,
t_0]$ and, as far as $[t_1, 1]$ is concerned, the path $p|[t_1, 1]$
will change only by a translation in $W_{aff}$ along the line
$L_\al$ parallel to the vector $\al^\vee$. Thus it remains to
describe the path $q=e_\al(p)$ on $[t_0, t_1]$. If $h_\al$ were not
to have any local minima on $[t_0, t_1]$ then $q|_{[t_0, t_1]}$
would be obtained by the reflection
$$
q|_{[t_0, t_1]}:= \tau_\al\circ p|_{[t_0, t_1]}
$$
and we would set
$$
q:= p|_{[0,t_0]} *\tau_\al\circ p|_{[t_0, t_1]} * p|_{[t_1,1]}.
$$
(Here we treat the paths resulting from the restriction of $p$ to
subintervals of $[0,1]$ as elements of ${\mathcal P}$, according to
the convention in section \ref{prelim}.)

This is the definition of $e_\al$ of \cite{Littelmann1}, however the definition of
$e_\al$ which we will need in this paper is the more refined
one of \cite{Littelmann2}. Call a  subinterval $[s,u]\subset
[t_0, t_1]$ a {\em spike} if it is a maximal interval satisfying
$$
h_\al(s)=h_\al(u) =\min( h_\al \restr [s,u]).
$$
 Thus $h_\al\restr [t_0,t_1]$ is decreasing on the complement to the
union of {\em spikes}. The restriction of $q$ to each {\em spike} is
obtained from $p$ by a translation along $L_\al$. The restriction to
each subinterval disjoint from a {\em spike} is obtained by a
reflection. To be more precise, subdivide the interval $[t_0,t_1]$
into
$$
[t_0, s_1]\cup [s_1,s_2]\cup ... \cup [s_k, t_1],
$$
where the {\em spike} and {\em non-spike} intervals alternate. Observe
that $[t_0, s_1], [s_k, t_1]$ are not {\em spikes}. Then
$$
q:= p|_{[0,t_0]} *\tau_\al(p|_{[t_0, s_1]}) * p|_{[s_1,s_2]}*...
*\tau_\al( p|_{[s_k, t_1]}) * p|_{[t_1,1]}.
$$

Note that the operator $e_\al$ changes the geometry of the path $p$
by an isometry near every point $p(t)$ which is neither a point of
local minimum for $h_\al$ nor is a point where $h_\al(t)=m_\al-1$.
Otherwise the local change is done by a ``bending''  with respect to
a hyperplane parallel to $H_\al$. These hyperplanes are not
necessarily walls of $W_{aff}$. However, if all local minimal values
of $h_\al$ belong to $\Z$, these hyperplanes are indeed walls and we
obtain:

For each path $p\in {\mathcal P}_{\Z,loc}$, for each simple root
$\al$, the path $q=e_\al(p)$ satisfies the following: The interval
$[0,1]$ can be subdivided into subintervals $[s_i, s_{i+1}]$ such
that the restriction $q|[s_i, s_{i+1}]$ is obtained from the
restriction of $p$ by post-composition with an element of $W_{aff}$.

\medskip
The lowering operators $f_\al$ are defined analogously to the
raising operators, we refer the reader to \cite{Littelmann2} for the
precise definition. (See however Property 1 below.) At this stage we
note only that $f_\al$ is undefined on $p$ iff $m_\al> h_\al(1)-1$.
Let ${\mathcal E}$ be the semigroup generated by $e_\al$'s,
${\mathcal F}$ be the semigroup generated by $f_\al$'s and
${\mathcal A}$ be the semigroup generated by all root operators. The
semigroups contain the identity operator by default. For each
$\phi\in {\mathcal A}$ let $Dom(\phi)$ denote the domain of $\phi$.

\begin{rem}
In fact, Littelmann extends the operators $f_\al, e_\al$ to the
entire ${\mathcal P}$ by declaring $f_\al(p)={\mathbf 0}$ for all
$p$ for which $f_\al$ is undefined. However we will not need this
extension in the present paper.
\end{rem}

Below we list certain properties of the root operators. Most of them
are either clear from the definition or are proven in
\cite{Littelmann2}. Most proofs that we present are slight
modifications of the arguments in \cite{Littelmann2}.

\medskip
{\bf Property 1.} (P. Littelmann, \cite[Lemma 2.1
 (b, e)]{Littelmann2}.)
$$
e_\al\circ f_\al(p)= p, \hbox{~~if~~} p\in Dom(f_\al),
$$
$$
f_\al\circ e_\al(p)= p, \hbox{~~if~~} p\in Dom(e_\al),
$$
$$
e_\al(p^*)= (f_\al(p))^*, \quad (e_\al(p))^*= f_\al(p^*),
$$
the latter could be taken as the definition of $f_\al$.

\medskip
{\bf Property 2.} \cite[Lemma 2.1]{Littelmann2}. For each
$p\in Dom(e_\al)\cap {\mathcal P}$,
$$
m_\al(e_\al(p))= m_\al(p)+1,
$$
$$
p\in Dom(e^N_\al) \iff N < |m_\al|.
$$

\medskip
{\bf Property 3.} Suppose that $p$ is a path in ${\mathcal P}_\Z$
which does not belong to the domain of any $e_\al, \al\in \Phi$.
Then $p$ is contained in $\Del$. Indeed, for each simple root $\al$
we have to have $m_\al(p)>-1$. Since $p\in {\mathcal P}_\Z$,
$m_\al(p)=0$. Thus $p\in {\mathcal P}^+$.

\medskip
{\bf Property 4.}  \cite[Proposition 3.1 (a, b)]{Littelmann2}.
For each $\al \in \Phi$, $Dom(f_\al)\cap {\mathcal P}_\Z$ is open
and $f_\al|{\mathcal P}_\Z$ is continuous.

\medskip
{\bf Property 5.} \cite[\S 7, Corollary 1 (a)]{Littelmann2}.
Let $p\in {\mathcal P}^+$ and $\phi$ be a composition of lowering
operators defined on $p$. Then $\phi(p)\in {\mathcal P}_\Z$ .

\medskip
{\bf Property 6.} Combining Properties 4 and 5 we conclude that for
each ${\mathbf f}\in {\mathcal F}$, $Dom({\mathbf f})\cap {\mathcal
P}^+$ is open and ${\mathbf f}|{\mathcal P}^+$ is continuous.

\medskip
{\bf Property 7.} \cite[Corollary 3, Page 512]{Littelmann2}.
$p$ is an LS path of the $\De$-length $\la$ if and only if there
exists ${\mathbf f}\in {\mathcal F}$ such that
$$
p={\mathbf f}(\pi_\la)
$$

\medskip
{\bf Property 8.} \cite[Corollary 2(a), page
512]{Littelmann2}. The set of LS paths of the given $\De$-length is
stable under ${\mathcal A}$.

\medskip
{\bf Property 9.} Suppose that $p\in {\mathcal P}$, $t\in [0,1]$,
$\al\in \Phi$ and $x:=p(t)$ satisfy
$$
\al(x)\in \Z,  \quad p'_-(t)\gtrsim_{W_x} p'_+(t).
$$
Then the path $q=e_\al(p)$ also satisfies
$$
q'_-(t)\gtrsim_{W_y} q'_+(t),
$$
for $y=q(t)$. \proof If $h_\al(t)\ne m_\al, m_\al-1$, the germs of
the paths $p$ and $q$ at $t$ differ by a translation. Thus the
conclusion trivially holds in this case. The same argument applies
if $h_\al(t)=m_\al-1$ and
$$
\al(p_-'(t)) \ge 0, \al(p_+'(t))\le 0.
$$
The nontrivial cases are:

1. $h_\al(t)=m_\al-1$, $\al(p'_-(t))\le 0$, $\al(p'_+(t))\le 0$. In
this case the assertion follows from Lemma \ref{easy} with
$\nu=p'_-(t), \mu=p'_+(t)$.

2. $h_\al(t)=m_\al$, $\al(p'_-(t))\le 0$, $\al(p'_+(t))\ge 0$. In
this case the assertion follows from Lemma \ref{hard} with
$\nu=p'_-(t), \mu=p'_+(t)$. \qed

\medskip
{\bf Property 10.}  Suppose that $p=p_1*p_2$ where $p_1\in {\mathcal
P}^+$. Then for each $e\in {\mathcal E}$ defined on $p$ we have
$$
e(p)= p_1* e(p_2).
$$
\proof It suffices to prove this for $e=e_\al$, $\al\in \Phi$. In
the latter case it follows directly from the definition of the
operator $e_\al$. \qed

The next property is again clear from the definition:

\medskip
{\bf Property 11.} Suppose that
$$
{\mathbf e}= e_{\be_m}\circ ... \circ e_{\be_1}
$$
where $\be_i\in \Phi$, $p\in Dom({\mathbf e})$. Set
$$
p_i:= e_{\be_i}\circ ...\circ e_{\be_1}(p), i=1,...,m.
$$
Then for each $T\in [0, 1]$, the sequence of vectors
$$
(p'_+(T), (p_1)_+'(T),..., (p_m)_+'(T)),
$$
after deleting equal members, forms a chain.

\begin{lem}
\label{finiteness} Given a path $p$ there are only finitely many
operators ${\mathbf e}\in {\mathcal E}$ which are defined on $p$.
\end{lem}
\proof Break the path $p$ as the concatenation
$$
p_1*...*p_s.
$$
of geodesic paths each of which is contained in a single alcove and
let $T_i\in [0, 1]$ be such that $p(T_i)=p_i(1/2)$; set $T_0:=0$.
Then for each $e_\al\in {\mathcal E}$ defined on $p$ there exists
$i$ such that the derivatives of ${e}_\al(p)$ and $p$ at $T_i$ are
not the same. Moreover,
$$
q={e}_\al(p)= q_1*...*q_s,$$ where each $q_i$ is a geodesic path
contained in an alcove. Consider the vector
$$
L(p):= (\ell( p'(0)), \ell(p'(T_1)),..., \ell(p'(T_s))) \in (\N\cup
\{0\}) ^{s+1},
$$
where $\N^{s+1}$ is given the lexicographic order and $\ell$ is the length
function on $W_{sph}$-orbits induced from the word metric on
$W_{sph}$ as in Proposition \ref{wordlength}. Then, by combining Proposition
\ref{wordlength} and Property 11 above, for each $\al\in \Phi$,
$$
L(e_\al(p))< L(p).
$$
Lemma follows. \qed

\subsection{Generalized LS paths}
\label{GLS}

 In this paper we will need two generalizations of the
concept of an LS path; the first one will be needed for the proof of
the saturation theorem (section \ref{saturationsection}), the second
will be used in section \ref{unfoldingLS} for the proof of the
unfolding theorem. Although we will use the name {\em generalized LS
path} for both generalizations, it will be clear from the context
which generalization is being referred to.

\medskip
{\bf The first generalization, ${\mathcal L}{\mathcal S}_1$}.
\footnote{This generalization of LS paths will be used in the proof of the saturation theorem.}

\medskip
Suppose we are given a collection of LS paths $p_i$ of the
$\De$-length $\la_i\in \De\cap P(R^\vee)$, $i=0,...,m$. We will use
the notation
$$
\ul{\la}=(\la_0,...,\la_m)
$$
and
$$
\la:= \sum_{i=0}^m \la_i.
$$

\begin{rem}
Actually, for our main application it will suffice to consider
$\la_i$'s which are multiples of the fundamental coweights
$\varpi_i$. Therefore such paths are automatically {\em
generalized Hecke paths} as defined in Definition \ref{gH}.
\end{rem}

\begin{defn}
\label{LS1} The concatenation
$$
p= p_0*p_1*...*p_m
$$
will be called a {\em generalized LS path} with $\length
(p)=\ul{\la}$, if for each $i=0,...,m-1$
$$
p'_i(1) \gtrsim  p'_{i+1}(0).
$$
The set of such generalized LS paths will be denoted ${\mathcal
L}{\mathcal S}_1$.
\end{defn}

This definition is a very special case of the one used by Littelmann
in \cite{Littelmann3} under the name of a {\em locally integral
concatenation}.

Recall that according to the definition of $\De$-length,
$$
\la= \Length(p).
$$
Observe that each LS path $p$ satisfies the above definition, since
for each $i$
$$
p'_i(1) \ge p'_{i+1}(0).
$$

\begin{ex}
Suppose that $u, v$ are dominant coweights. Then $p=\pi_u*\pi_v$ is
a generalized LS path.
\end{ex}

\medskip
{\bf The second generalization ${\mathcal L}{\mathcal S}_2$}.
\footnote{This notion of generalized LS path will be used only to unfold Hecke paths.}

\medskip
Suppose that $p_1, p_2\in {\mathcal P}$ appear as
$$
p_1= \tilde{p}_1|_{[0,a]}, \quad 0<a<1,
$$
$$
p_2= \tilde{p}_2|_{[b,1]}, \quad 0<b<1,
$$
where $\tilde{p}_1, \tilde{p}_2$ are LS paths, $a, b\in \Q$. (See
section \ref{paths} for the definition of $\tilde{p}_1|_{[0,a]}$
and $\tilde{p}_2|_{[b,1]}$.) Define the path $p:= p_1*p_2$. Assume
that
$$
p_1'(1)\rhd p_2'(0);
$$
in other words, if $t$ is such that $p(t)= p_1(1)$ then
$$
p'_-(t) \rhd p_+(t).
$$

\begin{defn}
The concatenation $p$ will be called a {\em generalized LS path}
if the {\em concatenation point} $p_1(1)$ is a regular
point\footnote{I.e. it does not belong to any wall.} of
$(A,W_{aff})$ and $p(1)\in P(R^\vee)$.

The set of such generalized LS paths will be denoted ${\mathcal
L}{\mathcal S}_2$.
\end{defn}

\begin{ex}
\label{e2} Suppose that $u\in\De, v\in V$ are such that $u, v \in
P(R^\vee)\otimes \Q$, $u+v\in P(R^\vee)\cap \De$ and the head of the
vector $u$ is a regular point in $(A,W_{aff})$. Then $p=\pi_u*
\pi_v\in {\mathcal L}{\mathcal S}_2$.
\end{ex}

This definition is again a very special case of the one given by
Littelmann in \cite[5.3]{Littelmann2}. Littelmann does not assume
that $p_1(1)$ is regular, but instead imposes certain chain
conditions at this point.

\medskip
{\bf Properties of generalized LS paths:}

\medskip
{\bf Property 0.} If $p\in {\mathcal L}{\mathcal S}_1$ then $p^*$ is
also in ${\mathcal L}{\mathcal S}_1$.

\proof Represent $p$ as a concatenation $p_1*...*p_m$ of LS paths as
in the Definition \ref{LS1}. Then
$$
p^*=(p_m^*)*...*(p_1^*),
$$
where each path $p_i^*$ is again an LS path. Now the assertion
follows from Lemma \ref{iff}. \qed

\medskip
{\bf Property 1.} ${\mathcal L}{\mathcal S}_2$ is stable under the
root operators  \cite[Lemma 5.6, 2-nd assertion]{Littelmann2}.
In particular, suppose that $p$ is as in Example \ref{e2}. Then for
each ${\mathbf f}\in {\mathcal F}$  (defined on $p$), ${\mathbf
f}(p)$ is a generalized LS path of the $\De$-length $\la=u+v$.

{\bf Property 2.}
 ${\mathcal L}{\mathcal S}_1$ is stable
 under the root operators.

\proof Suppose that $p=p_1*...*p_m$ is a concatenation of LS paths
as above and $e_\al$ is a raising operator. In particular, for each
$i$ we have a vector $u_i$ so that
$$
p'_i(1) \ge u_i \sim p'_{i+1}(0).
$$
For each $i$, $e_\al(p_i)$ is again a Littelmann path. Therefore
$e_\al(p)$ is a concatenation of LS paths $q_1*...*q_m$. We have to
verify that for each $i$ there is a vector $v_i\in V$ so that
$$
q'_i(1) \ge v_i \sim q'_{i+1}(0).
$$
This however follows from the Property 9 in the previous section. To
check that ${\mathcal L}{\mathcal S}_1$ is preserved by $f_\al$ we
use that $q\in {\mathcal L}{\mathcal S}_1\iff q^* \in {\mathcal
L}{\mathcal S}_1$ and
$$
f_\al(p)= (e_\al(p^*))^*.  \qed
$$

\medskip
{\bf Property 3.} ${\mathcal L}{\mathcal S}_1$ and ${\mathcal
L}{\mathcal S}_2$ are contained in ${\mathcal P}_{\Z,loc}$. For
${\mathcal L}{\mathcal S}_1$ it is immediate since, by Theorem
\ref{locint}, the set of LS
paths is contained in ${\mathcal P}_{\Z,loc}$.
For ${\mathcal L}{\mathcal S}_2$ it is a special case
of \cite[Lemma 5.5]{Littelmann2}.

\medskip
{\bf Property 4.} Suppose that $p\in {\mathcal L}{\mathcal S}_1$.
Then there exists an element ${\mathbf e}\in {\mathcal E}$
defined on $p$ such that $q={\mathbf e}(p)\in {\mathcal P}^+$.

\proof If $m_{\al_1}(p)\le -1$ then we apply a power
$e_{\al_1}^{k_1}$ to $p$ so that $q_1:=e_{\al_1}^{k_1}(p)$ satisfies
$m_{\al_1}(q_1)>-1$. However, since $e_{\al_1}^{k_1}(p)\in {\mathcal
L}{\mathcal S}_1\subset {\mathcal P}_{\Z,loc}$, it follows that
$q_1\in {\mathcal P}_{\Z,loc}$ and so $m_{\al_1}(q_1)=0$. We then
apply a power of $e_{\al_2}$ to $q_1$, etc. According to Lemma
\ref{finiteness}, this process must terminate. Therefore, in the end
we obtain a path
$$
{\mathbf e}(p)=q
$$
which does not belong to the domain of any raising operator. Since
$q\in {\mathcal P}_{\Z,loc}$ it follows that $q$ is entirely
contained in $\De$. \qed

Recall  \cite[Chapter VI, section 10]{Bourbaki} that if a root system
$R$ spans $V$ then each dominant coweight $\la\in \De$ is a
positive integral combination
$$
\la= \sum_{i=1}^l n_i \varpi_i,
$$
where $\varpi_i$ are fundamental coweights. This assertion (as it stands) is false without the above assumption
on $V$. In the general case we  have
$$
\la= \la'+ \sum_{i=1}^l n_i \varpi_i,
$$
where $\la'\in V'$, $n_i\in \N\cup \{0\}$.
As an alternative the reader can restrict the discussion to semisimple groups only,
when $V'=0$.

\begin{conv}
\label{con} From now on we will be assuming that in Definition
\ref{LS1}
$$
\la_j=k_j \varpi_j, \quad k_j\in \N,
$$
for each $j=1,...,m$, where $\varpi_j$ is the $j$-th fundamental coweight, and
$$
\la_0\in V'.
$$
Then the subpath $p_0$ is necessarily geodesic.
\end{conv}

\begin{lem}
\label{positiv} Suppose that $p\in {\mathcal L}{\mathcal S}_1 \cap
{\mathcal P}^+$ is a generalized LS path with $\length(p)=\ul{\la}$.
Then
$$
p=\pi_{\la_0}*...*\pi_{\la_m}.
$$
\end{lem}
\proof Represent $p$ as the concatenation of maximal LS subpaths,
$p=p_0*p_1...*p_m$. The geodesic subpath $p_0$ clearly equals $\pi_{\la_0}$.
Since $p_1$ is an LS path and $p_1'(0)\in \De$, we
see that $p_1$ is a geodesic path (see Corollary \ref{end}) which
therefore equals $\pi_{\la_1}$. Moreover, because
$$
p_1'(1)\ge u_1 \sim p'_2(0),
$$
it follows that $u_1=p_1'(1)$ and thus $p_1'(1)=\la_1\sim p_2'(0)$.
Let $x_1:= p_1(1)$; this point lies on the boundary face of $\De$
which does not contain $\la_2$. Note that the vector $p_2'(0)$,
regarded as an element of $T_{x_1}(A)$, points inside the Weyl
chamber $\De$ (for otherwise $p$ is not contained in $\De$). On the
other hand, since $\la_1\sim p_2'(0)$, the vector $p_2'(0)$ belongs
either to $\De$ or to the Weyl chamber
$$
\tau_{\be_2}(\De)
$$
adjacent to $\De$. Since $p\in {\mathcal P}^+$, it is clear that
$p_2'(0)\in \De$. Thus $p_2$ is the geodesic path $\pi_{\la_2}$.
Continuing in this fashion we conclude that
$$
p=\pi_{\la_0}*\pi_{\la_1}*...*\pi_{\la_m}. \qed
$$

\begin{thm}
\label{pathlow}
 Suppose that $p$ is a generalized LS path in the sense of ${\mathcal L}{\mathcal S}_1$
 with $\length(p)=\ul{\la}$
(satisfying convention \ref{con}). Then there exists ${\mathbf f}\in
{\mathcal F}$ such that
$$
p={\mathbf f}(\pi_{\la_0}*...*\pi_{\la_m})
$$
\end{thm}
\proof If $p\in {\mathcal P}^+$ then we are done. Otherwise, by
combining Lemma \ref{positiv} with the Property 4, we find an
${\mathbf e}\in {\mathcal E}$,
$$
{\mathbf e}= e_{\al_1}^{k_1}\circ ... e_{\al_n}^{k_n}
$$
 such that $p\in Dom({\mathbf e})$ and ${\mathbf e}(p)=q\in
{\mathcal P}^+$. Therefore
$$
q=\pi_{\la_0}*\pi_{\la_1}*...*\pi_{\la_m}
$$
and thus the composition
$$
{\mathbf f}= f_{\al_n}^{k_n}\circ ... f_{\al_1}^{k_1}
$$
satisfies $p={\mathbf f} {\mathbf e}(p)={\mathbf
f}(\pi_{\la_0}*\pi_{\la_1}*...\pi_{\la_m})$. \qed

\subsection{Path model for the representation theory of Lie groups}
\label{LSrep}

Suppose that
$$p(t)\in {\mathcal L}{\mathcal S}_1 \cup {\mathcal L}{\mathcal S}_2$$
is a generalized LS path with $\Length(p)=\be$ and
$$\length(p)=\ul{\la}.$$
Suppose that $\al\in P(R^\vee)$ is such that $\al+p(t)$ is contained
in $\De$. Then $\al$ and $p$ define a polygon
$$
P:=\ol{o y_0} \cup (p+\al)\cup \ol{y_n o} \subset \De$$
 where $\al=\ov{o y_0}$, $y_n=\al+p(1)$. Let $\ga$ denote the vector
$\ov{o y_{n}}$; then $\ga$ is also a dominant coweight. Recall that
the contragredient dominant coweight $\ga^*\in \De$ is obtained by
projecting the vector $-\ga$ to the Weyl chamber $\De$ by the
projection $\P: V\to \De$.

\begin{defn}
1. A polygon $P$ above will be called a {\em (broken) Littelmann
polygon} with the $\De$-side lengths $\al, \be, \ga^*$.

2. If $p(t)$ is an LS path then $P$ will be called a {\em (broken)
Littelmann triangle} with the $\De$-side lengths $\al, \be, \ga^*$.
\end{defn}

Pick a lattice $L$ such that
$$
Q(R^\vee)\subset L\subset P(R^\vee).
$$
Then there exists a unique connected semisimple complex Lie group
$G^\vee$  with the root system $R^\vee$ and the character lattice
$L$ of the maximal torus $T^\vee\subset G^\vee$.
 Recall that irreducible representations $V$ of $G^\vee$ are
parameterized by their dominant weights, $V=V_\la$, $\la\in \De\cap
L$.

Pick a path $q\in {\mathcal P}^+$ such that $q(1)=\be$. Then,
according to \cite[Decomposition formula, Page 500]{Littelmann2} we
have

\begin{thm}
\label{Lit2} The tensor product $V_\al\otimes
V_\be$  contains $V_\ga$ as a subrepresentation if and only if there
exists a path $p\in {\mathcal F}(q)$ such that $\pi_{\al}*p\in
{\mathcal P}^+$ and $\pi_{\al}*p(1)=\ga$.
\end{thm}

\begin{rem}
Littelmann works with simply-connected group $G^\vee$ and weights
$\al, \be,\ga$ in $P(R^\vee)$. The statement for
non-simply-connected groups trivially follows from the
simply-connected case.
\end{rem}

In particular, since $p$ is an LS paths of the $\De$-length $\be$ if
and only if $p\in {\mathcal F}(\pi_\be)$, it follows that

\begin{thm}
\cite{Littelmann1, Littelmann2}. \label{Lit} The tensor
product $V_\al\otimes V_\be$  contains $V_\ga$ as a
subrepresentation if and only if there exists a (broken) Littelmann
triangle in $\Del\subset V$, with the $\De$-side-lengths $\al, \be,
\ga^*$.

In other words, $V_\ga\subset V_\al\otimes V_\be$ if and only if
there exists an LS path $p$ of $\De$-length $\be$ such that

1. $\pi_\al *p\in {\mathcal P}^+$.

2. $p(1)+\al=\ga$.
\end{thm}

We will apply Theorem \ref{Lit2} as follows.
Represent the vector $\be$ as the integer linear combination of
fundamental coweights
$$
\be =\sum_{i=1}^n k_i \varpi_i.
$$
We reorder the fundamental coweights so that $k_i>0$ for all
$i=1,...,m$ and $k_i=0, i\ge m+1$. Set $\la_i:= k_i \varpi_i, 1\le
i\le m$ and let $\ul{\la}= (\la_1,...,\la_m)$. Therefore the path
$$
\pi_{\ul\la}:= \pi_{\la_1}*... \pi_{\la_m}
$$
belongs to ${\mathcal L}{\mathcal S}_1$ and $\pi_{\ul\la}(1)=\be$.
Then

\begin{cor}
The tensor product $V_\al\otimes V_\be$ contains $V_\ga$ as a
subrepresentation
 if and only if there exists a generalized LS path $p$ so that

1. $\length(p)=\ul{\la}$.

2. $\pi_{\al}*p\in {\mathcal P}^+$.

3. $\pi_{\al}*p(1)=\ga$.
\end{cor}
\proof Set
$$
q=\pi_{\ul{\la}}.$$
 According to Theorems \ref{pathlow} and Property 2 of
generalized LS paths (section \ref{GLS}), $p\in {\mathcal P}$ is a
generalized LS path with $\length(p)=\ul\la$ if and only if $p\in
{\mathcal F}(q)$. Now the assertion follows from Theorem \ref{Lit2}.
\qed

Combining Corollary \ref{folded-chain}, Theorem \ref{T3} and Theorem
\ref{Lit} we obtain

\begin{cor}
Suppose that $X$ is a thick Euclidean building modeled on the
Coxeter complex $(A,W_{aff})$. Let $\al, \be, \ga^*\in L$ be
dominant coweights. Suppose that $a\subset A$ is an alcove
containing a special vertex $o$, $T=[o, x, y]\subset X$ is a
geodesic triangle with the special vertices  and the $\De$-side
lengths $\al, \be, \ga^*$. Assume also that the broken side
$Fold_{a,A}(\ol{x y})$ of the folded triangle
$$
Fold_{a,A}(T)
$$
has breaks only at the special vertices of $A$. Then

1. The folded triangle $T'=Fold_{o,id,\De}(T)\subset \De$ is a
Littelmann triangle.

2. $V_\ga\subset V_\al\otimes V_\be$.
\end{cor}
\proof Indeed, according to Corollary \ref{folded-chain}, the folded
triangle $T'$ satisfies the chain condition. Each break point $x_i$
on the broken side of $T'$ is

1. Either a special vertex, in which case it satisfies maximal chain
condition by Remark \ref{R1}, or

2. $Fold_{a,A}(\ol{xy})$ is geodesic at the point corresponding to
$x_i$, so the chain at $x_i$ can be chosen to be maximal by
Proposition \ref{chaincondition1}.

Hence Theorem \ref{T3} implies that $T'$ is a Littelmann triangle.
The second assertion now follows from Theorem \ref{Lit}. \qed

\medskip
Of course, the assumption that the break points occur only at the
special vertices is very restrictive. In section
\ref{saturationsection} we will get rid of this assumption at the
expense of dilation of the side-lengths.

\section{Unfolding}
\label{last}

The goal of this section is to establish an intrinsic
characterization of folded triangles as the broken billiard
triangles satisfying the chain condition. We first prove this
characterization for Littelmann triangles and then, using this, give
a general proof.

Throughout this section we assume that $X$ is a thick locally
compact Euclidean building modeled on the Coxeter complex $(A,W)$,
$\Del\subset A$ is a Weyl chamber with tip $o$. Let $g: X\to
\De$ denote the folding $Fold_{\De}$.

Let $T\subset \Del$ be a  billiard triangle which is the
composition
$$
T=\ov{ox}\cup r \cup \ol{y o},
$$
where $r(t)=p(t)+\al$, $\al=\ov{ox}$ and $p\in {\mathcal P}$ is a Hecke path. Thus $T$
has the geodesic sides $\ol{ox}, \ol{oy}$ and the broken side
$r$. We set $\ga:= \ov{o y}$ and let $\be\in \De$ denote
$\De$-length of the path $p$.

\subsection{Unfolding Littelmann triangles}
\label{unfoldingLS}

\begin{thm}
\label{LSfolded}
 Suppose that, in addition, $T$ is a {\em Littelmann triangle},
i.e. $\al, \ga\in L\subset P(R^\vee)$ and $p$ is an LS path. Then
$T$ can be {\em unfolded} in $X$, i.e. there exists a geodesic
triangle $\t{T}\subset X$ such that $g(\t{T})=T$.
\end{thm}
\proof Here is the idea of the proof: We know that billiard
triangles in $\De$ satisfying the {\em simple chain condition} can be
unfolded to geodesic triangles in $X$, see Corollary
\ref{chain-folded}. Littelmann triangle $T$ is billiard, satisfies
the chain condition, but not necessarily the {\em simple chain condition}.
Our goal is to approximate $T$ by Littelmann polygons $P_\eps$,
$\lim_{\eps\to 0} P_\eps=T$, which satisfy the simple chain condition.
We then unfold each $P_\eps$ to a geodesic quadrilateral
$\t{T}_\eps\subset X$. Since $X$ is locally compact, there is a
convergent sequence $\t{T}_{\eps_j}$ whose limit is a geodesic triangle
$\t{T}$ which folds to $T$. Below is the detailed argument.

\begin{figure}[tbh]
\centerline{\epsfxsize=3.5in \epsfbox{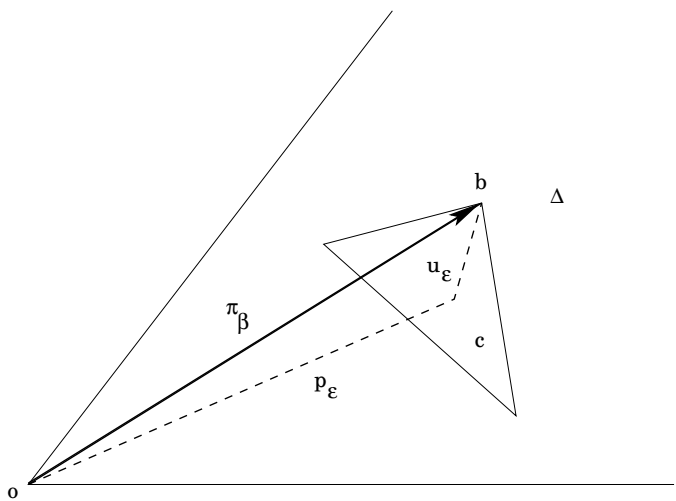}} \caption{\sl
Approximation.} \label{unfolding.fig}
\end{figure}

\medskip
Consider the geodesic path $\pi_\be=\ol{ob}\in {\mathcal P}^+$.
Since $p$ is an LS path with the $\De$-length $\be$, according to
Property 7 in section \ref{operators}, there exists a composition
$\phi\in {\mathcal F}$ of lowering operators so that
$$
\phi(\pi_\be)= p.
$$

Let $c\subset \De$ denote an alcove which contains the germ of the
segment $\ol{ob}$ at $b$. Pick a point $u$ in the interior of
$\ol{ob}\cap c$. Then for each $\eps>0$ there exists a
point $u_\eps\in int(c)\cap P(R^\vee)\otimes \Q$ such that

1. $|u-u_\eps|<\eps$.

2. The segments $\ol{ou_\eps}, \ol{u_\eps b}$ do not pass
   through any point of intersection of two or more walls
   (except for the end-points of these segments).

Observe that $u_\eps$ is a regular point in $(A,W_{aff})$, i.e. it
does not belong to any wall.

In other words, the path
$$
\hat{p}_\eps:= \ol{ou_\eps}\cup \ol{u_\eps b} \in {\mathcal P}
$$
is {\em generic}. Parameterize $\hat{p}_\eps$ with the constant speed
so that $\hat{p}_\eps(t_\eps)=u_\eps$. See Figure \ref{unfolding.fig}.

Then the path $\hat{p}_{\eps}$ belongs to ${\mathcal P}^+$; clearly it
is also a generalized LS path: $\hat{p}_\eps\in {\mathcal L}{\mathcal
S}_2$. Moreover,
$$
\lim_{\eps\to 0} \hat{p}_{\eps}= \pi_\be.
$$
Therefore, according to Property 6 of the root operators (see
section \ref{operators}), the operator $\phi$ is defined on all
$\hat{p}_{\eps}$ for $\eps$ sufficiently small and
$$
\lim_{\eps\to 0} \phi(\hat{p}_{\eps})= \phi(\pi_\be)=p.
$$
Set $p_\eps:= \phi(\hat{p}_\eps)$. Since $\hat{p}_\eps$ was generic, the
path $p_\eps$ is generic as well. By construction, for each
sufficiently small $\eps$,
$$
p_\eps(1)=p(1).
$$
Observe also that the germ of the path $p_{\eps}$ at the point
$p_\eps(t_\eps)$ is isomorphic (via an element of
$W_{aff})$ to the germ of $\hat{p}_{\eps}$ at $u_\eps$ (since $u_\eps$
is regular). Similarly, $p_\eps$ is the composition of the path
$p_\eps|[0, t_\eps]$ with the path that belongs to the
$W_{aff}$-orbit of $\ol{u_\eps b}$.

For each $\eps$ we form a new polygon $P_\eps$ by replacing
the broken side $r(t)=\al+p(t)$ (in $T'$) with the path $\al+p_\eps(t)$.
Clearly,
$$
\lim_{\eps\to 0} P_\eps=T.
$$

To simplify the notation we now fix $\eps>0$ and let $q:= p_\eps$.

\begin{lem}
For all sufficiently small $\eps$, the polygon $P_\eps$ is
contained in $\De$.
\end{lem}
\proof Suppose that $\la$ is a simple root which is negative at some
point of the path $\al+q(t)$.

Since $\la$ is nonnegative on the limiting path $\al+p$, the minimum
of the function $J_\la(t):= \la(q(t)), t\in [0,1]$, converges to zero as $\eps\to 0$.
However, as a generalized LS path, $q$ belongs to ${\mathcal P}_\Z$
(see Property 5 in section \ref{operators}). Since $\al\in
P(R^\vee)$, it follows that the minimum of  $J_\la(t)$ is an
integer. Hence it has to be equal to zero for all sufficiently small
values of $\eps$. Contradiction. \qed

Since each $q$ is a generalized LS path and $P_\eps\subset \De$,
the polygon $P_\eps$ is a Littelmann polygon. Moreover, since $q$
is {\em generic}, the polygon $P_\eps$ satisfies the {\em simple
chain condition}. Thus

1. For each $t\in [0,t_{\eps})$ either $q$ is smooth at $t$ or
$$
(q'_-(t), q'_+(t))
$$
is a chain of length 1: At $m=q(t)$ the above tangent vectors are
related by a single reflection in $W_{m}$. The fixed-point set of
this reflection is the unique wall passing through $m$.

2. The subpath $q([t_\eps,1])$ in $q$ is a geodesic segment and
$$
\del(\eps):= \angle(\ov{u_\eps o}, \ov{u_\eps b})= \pi- \angle(
q'_-(t(\eps)), q'_+(t(\eps))).
$$

Now we are in position to apply Corollary \ref{chain-folded} and
unfold $P_\eps$ in $X$: For each $\eps$ there exists a geodesic
quadrilateral $\t{T}_\eps$ (with one vertex at $o$) in $X$ such that
$$
g(\t{T}_\eps)= P_\eps.
$$
Let $\t{z}=\t{z}_\eps$ denote the point of $\t{T}_\eps$ which maps to $z=q(t_\eps)$
under the folding map $f$. Since $z$ is a regular point, the point
$\t{z}$ is regular as well and the angle between the sides of $\t{T}_\eps$
at $\t{z}$ is the same as the angle between the sides of $P_\eps$ at
$z$, i.e. equals $\del(\eps)$.

Since the building $X$ is locally compact, the sequence of
quadrilaterals $\t{T}_\eps$ subconverges to a geodesic quadrilateral
$T\subset X$ which is a geodesic triangle since
$$
\lim_{\eps\to 0}\del(\eps)=\pi.
$$
By continuity of the folding $g: X\to\De$,
$$
g(T)= \lim_{\eps\to 0} P_{\eps} =T. \qed
$$

In the above proof we assumed that the polygon $T$ is entirely
contained in $\De$. This assumption can be weakened. Let $f: X\to A$
denote the folding $Fold_{a,A}$ into the apartment $A$, where $a$ is
an alcove containing $o$. Suppose that $T\subset A$ is as above, so
that $\al, \ga\in P(R^\vee)$, $p$ is a billiard path, $r=p+\al$. Define two subsets $J, J'\subset I=[0,1]$:
$$
J:= cl( r^{-1}(int( \De))), \quad J':= cl( r^{-1}(int(
V\setminus \De))).
$$
Clearly, $I=J\cup J'$ and the set $J\cap J'$ is finite.

We assume that for each $t\in J$ the germ of $p$ at $t$ satisfies
the maximal chain condition, and for each $t\in J'$ the germ of $p$
at $t$ is geodesic.

\begin{thm}
\label{addendum} Under the above assumptions the polygon $T$ can be
unfolded to a geodesic triangle in $X$ via the retraction $f$.
\end{thm}
\proof Recall that unfolding of $T$ is a local problem of behavior of
the path $r$ at the break-points (Lemma \ref{locality}), which
in our case all occur inside $\De$.

We first replace $T$ with the polygon $P=\P(T)$, where
$\P=\P_{\De}$ is the projection of $A$ to the Weyl chamber $\De$.
Then, analogously to the proof of Proposition \ref{chaincondition1},
the new polygon $P$ still satisfies the maximal chain condition.
Therefore, according to the previous theorem, the polygon $P$
unfolds in $X$ via the folding map $g=Fold_{\De}: X\to \De$. However
this means that the unfolding condition (stated in Lemma
\ref{locality}) is satisfied at each break-point of the polygon $T$
(since the germs of $r$ and of $\P(r)$ are the same). Hence
the original polygon $T$ unfolds to a geodesic triangle in $X$ via $f: X\to A$. \qed

\medskip
Let $\ul{G}$ be a connected split semisimple algebraic group with
the root system $R$ and the cocharacter
lattice $L$ of the maximal torus $\ul{T}\subset \ul{G}$. We let
$\ul{G}^\vee$ denote its Langlands' dual and set
$$
G^\vee:=\ul{G}^\vee(\C).
$$
We assume that $\al,\be,\ga^*\in L$ are dominant weights of $G^\vee$
such that
$$
(V_\al\otimes V_\be \otimes V_{\ga^*})^{G^\vee}\ne 0,$$
equivalently,
$$
V_{\ga}\subset V_\al\otimes V_\be.
$$

As a corollary of Theorem \ref{LSfolded} we get a new proof of

\begin{thm}
[Theorem 9.17 in \cite{KLM3}, also proven in \cite{Haines}] \label{klm3}
  Under the above assumptions,
in the thick Euclidean building $X$ there exists a geodesic triangle
with special vertices and the $\De$-side lengths $\al, \be, \ga^*$. In other words,
$$
n_{\al,\be}(\ga)\ne 0 \Rightarrow m_{\al,\be}(\ga)\ne 0.
$$
\end{thm}
\proof Since
$$
V_{\ga}\subset V_\al\otimes V_\be,
$$
according to Littelmann's Theorem \ref{Lit}, there exists a
Littelmann triangle $T'\subset \De$, as in Theorem \ref{LSfolded}.
Let $p\in {\mathcal P}$ denote the LS path (of the $\De$-length
$\be$) representing the broken side of $T'$; $p=\phi(\pi_\be)$,
where $\phi\in {\mathcal F}$ is a composition of lowering operators.
Thus, by Theorem \ref{LSfolded}, there exists a triangle $T=[o, x,
y]\subset X$ such that $Fold_\De(T)=T'$. Therefore, by the
definition of folding,
$$
d_{ref}(o, x)= d_{ref}(o,x')=\al, \quad d_{ref}(o, y)=
d_{ref}(o,y')=\ga.
$$
Since we assumed that $\al, \be\in L\subset P(R^\vee)$ then $x, y$
are special vertices of $x$. Since folding preserves the
$\De$-length,
$$
d_\De(x,y)= \Length(p)=\be. \qed
$$

\subsection{Characterization of folded triangles}
\label{characterizationsection}

The goal of this section is to extend the results of the previous
one from the case of Littelmann triangles to general broken
triangles satisfying the chain condition.

\begin{thm}
\label{chain-unfolding} Suppose that $p\in {\mathcal P}$ is a Hecke
path, $\al=\ov{o u}\in \De$ is such that the path $q:= \al+p$ is
contained in $\De$. Define the billiard triangle $T':=\ol{o u}\cup q
\cup \ol{ q(1) o}$.
 Then $T'$ can be unfolded in $X$.
\end{thm}
\proof The idea of the proof is that the set of unfoldable billiard
paths is closed, thus it suffices to approximate $p$ by unfoldable
paths. We first prove the theorem in the case when $o$ does not
belong to the image of the path $q$.

According to Lemma \ref{locality}, unfolding of a path is a purely
local matter. Therefore the problem reduces to the case when $q$ has
only one break-point, $x=q(t_1)$. If $p$ were an LS path, we would
be done. In general it is not, for instance, because it might fail
the {\em maximal chain condition}. We resolve this difficulty by
{\em passing to a smaller Coxeter complex} and a smaller building.

Let $R_x$ denote the root subsystem in $R$ which is generated by the
roots corresponding to the walls passing through $x$. This root
system determines a Euclidean Coxeter complex where the stabilizer
of the origin is a finite Coxeter group $W'_{sph}$ which is
conjugate to the group $W_x$ via the translation by the vector
$\ov{ox}$. Let $\De_x$ denote the positive Weyl chamber of $(V,
W'_{sph})$ (the unique chamber which contains $\De$). Let $\xi$,
$\eta$ and $\zeta$ denote the normalizations of the vectors
$-p'_-(t_1), p'_+(t_1), \ov{xo}$.

Then, since $p$ satisfies the chain condition, there exists an
$(S,W_x,\De_x)$-chain
$$
(\nu_0,...,\nu_m), \quad \nu_0=-\xi, \quad \nu_m=\eta, \quad
\nu_i=\tau_i(\nu_{i-1}), \quad 1\le i\le m.
$$
Our first observation is that although this chain may fail to be a
maximal chain with respect to the unrestricted Coxeter complex
$(S,W_{sph})$, we can assume that it is maximal with respect to the
restricted Coxeter complex $(S,W_x)$.

Next, the initial and final points of $q$ may not belong to
$P(R^\vee_x)$. Recall however that {\em rational points} are dense
in $S$, see Lemma \ref{density}; therefore, there exists a sequence
of rational points (with respect to $R_x$) $\xi_j\in S$ which
converges to $\xi$. Thus, using the same reflections $\tau_i$ as
before, we obtain a sequence of rational chains $(\nu_i^j),
i=0,...,m$, where $\nu_0^j=-\xi_j, \nu_m^j=\eta_j$. We set
$\zeta_j:= \zeta$.

Hence for each $j$ there exists a number $c=c_j\in \R_+$ so that the
points
$$
x_j=x+ c\xi_j, \quad y_j:= x+ c\eta_j
$$
belong to $P(R^\vee)$.  We define a sequence of paths
$$
q_j:=\ol{x_j x}\cup \ol{x y_j} \in \t{\mathcal P}.
$$
Our next goal is to choose the sequence $\xi_j$ so that the germ of
each $q_j$ at $x$ is contained in $\De$. If $x$ belongs to the
interior of $\De$ then we do not need any restrictions on the
sequence $\xi_j$. Assume therefore that $x$ belongs to the boundary
of $\De$. Let $F$ denote the smallest face of the Coxeter complex
$(V, W_{sph})$ which contains the point $x$ and let $H$ denote the
intersection of all walls through the origin which  contain $x$. It
is clear that $F$ is a convex homogeneous polyhedral cone contained
in $H$ and $x$ belongs to the interior of $F$ in $H$. If $w\in W_x$
is such that $w(-\xi)= \eta$ then $w$ fixes $H$ (and $F$) pointwise.

By Lemma \ref{density}, applied to the root system $R_x$, there
exists a sequence of unit rational vectors $\xi_j$ and positive
numbers $\eps_j$ converging to zero so that points $x\pm \eps_j
\xi_j$ belong to $int_H(F)$; therefore the sequence $w(x+ \eps_j
\xi_j)$ is also contained in $int_H(F)$.

Using this sequence $\xi_j$ we define the paths $q_j$; clearly the
germ of $q_j$ at $x$ is contained in $int_H(F)\subset \De\subset
\De_x$.

\begin{rem}
Note that, typically, the sequence $(c_j)$ is unbounded and the
paths $q_j$ are not contained in $\De_x$.
\end{rem}

We let $p_j\in {\mathcal P}$ denote the path $q_j- q_j(0)$. Then
each $p_j$ is an LS path with respect to the root system $R_x$:
Integrality and the maximal chain condition now hold. Set
$$
\la_j:= \hbox{length}_{\De_x}(p_j).
$$

\begin{rem}
Observe that,
$$
\lim_j \bar\la_j = \bar\la\in \De,
$$
where $\la$ is the $\De_x$-length of $p$.
\end{rem}

Therefore, according to Theorem \ref{addendum} for each $j$ the path
$q_j$ is unfoldable in a thick Euclidean building $X_x$ modeled on
the Coxeter complex
$$
(A, W'_{aff}), \hbox{~~where~~} W'_{aff}= V\ltimes W_x.
$$
This means that there exists a geodesic path $\tilde{q}_j$ in $X_x$,
whose $\De_x$-length is $\la_j$, and which projects to $q_j$ under
the folding $X_x\to A$.

Let $z_j\in \t{q}_j$ be the points which correspond to the point $x$
under the folding map
$$
\t{q}_j \to q_j.
$$
Thus the ``broken triangle'' $[\xi_j, \zeta, \eta_j]$ in $S_x$
unfolds in $\Si_{z_j}(X_x)$ into a triangle $[\t\xi_j,
\t\zeta_j,\t\eta_j]$ such that
$$
d_{ref}(\t\xi_j, \t\zeta_j) = d_{ref}(\xi_j, \zeta_j),
$$
$$
d_{ref}(\t\eta_j, \t\zeta_j) = d_{ref}(\eta_j, \zeta_j),
$$
$$
d(\t\eta_j,\t\zeta_j)=\pi.
$$

Observe that the metric distance from $o$ to $z_j$ is uniformly
bounded. Since $X_x$ is locally compact, the sequence of buildings
$\Si_{z_j}(X_x)$ subconverges to the link of a vertex $u\in
X_x\subset X$.

\begin{rem}
The spherical buildings $\Si_{z_j}(X_x)$, $\Si_{u}(X_x)$ have to be
modeled on the same spherical Coxeter complex $(S,W_x)$, since the
structure group can only increase in the limit and the structure
group at $z_j$ was already maximal possible, i.e. $W_x$.
\end{rem}

Accordingly, the triangles $[\t\xi_j, \t\zeta_j,\t\eta_j]$
subconverge to a triangle $[\t\xi, \t\zeta,\t\eta]$ whose refined
side-lengths are
$$
d_{ref}(\xi, \zeta),   d_{ref}(\zeta, \eta),\pi.
$$
This shows that the triangle $[\xi,\zeta,\eta]$ can be unfolded in a
building which is modeled on $(S,W_x)$. We now apply the {\em
locality} lemma \ref{locality} to conclude that the path $q$ can be
unfolded in $X$ to a geodesic path. Thus the broken triangle $T'$
unfolds to a geodesic triangle as well.

\medskip
If $o$ belongs to the image of $q$ we argue as follows. The path
$q$, as before, has only one break point, which in this case occurs
at the origin:
$$
q= \ol{zo}\cup \ol{oy}.
$$
There exists an element $w\in W_{sph}$ which sends the vector
$\eta=-\xi$ to $\xi$, where $\eta$ is the normalization of the
vector $\ov{oy}$. Then consider the geodesic path
$$
\tilde{q}:= \ol{zo}\cup w(\ol{oy}).
$$
It is clear that $g(\tilde{q}) = \P(\tilde{q}) =q$. \qed

By combining Theorem \ref{chain-unfolding} and Corollary
\ref{folded-chain} we obtain the following

\begin{thm}
[Characterization of folded triangles]\label{characterization} A
polygon $P\subset \De$ of the form
$$
\ol{ox} \cup (p+\al) \cup \ol{yo}, \hbox{~where~}\al=\ov{ox},
\ga=\ov{oy}\in \De,
$$
can be unfolded to a geodesic triangle in $X$ if and only if $p$ is
a Hecke path, i.e. a billiard path which satisfies the chain
condition.
\end{thm}

\section{Proof of the saturation theorem}
\label{saturationsection}

We first prove Theorem \ref{3->4} formulated in the Introduction. Part 1 of Theorem was proven
in \cite{KLM3}, so we prove Part 2.
Let $X$ be a (thick) Euclidean building of rank $r$ modeled on a
discrete Coxeter complex $(A,W)$; the building $X$ is the Bruhat-Tits building associated with
the group $G=\ul{G}(\K)$. Then the assumption  that $m_{\al,\be}(\ga)\ne 0$ is equivalent to
the assumption that there exists a geodesic triangle $T=[\t{x}, \t{y}, \t{z}]\subset X$,
where $\t{x}, \t{y}, \t{z}$ are special vertices of
$X$ and whose $\De$-side-lengths are $\al, \be, \ga^* \in \De\cap P(R^\vee)$.

Recall that there exists an apartment $\tilde{A}\subset X$ which
contains the segment $\ol{\t{x} \t{y}}$. We let $\t{W}$ denote the
affine Weyl group operating on $\t{A}$. Our first step is to replace
the geodesic triangle $T$ with a geodesic polygon
$$
\tilde{P} := [\tilde z, \t{x}=\tilde{x}_1,..., \tilde{x}_n,
\tilde{x}_{n+1}=\tilde{y}]
$$
as follows. We now treat the point $\t{x}$ as the origin in the
affine space $\t{A}$. Let $\tilde\De\subset \tilde{A}$ be a Weyl chamber in $(\t{A},\t{W})$,
so that $\t\De$ has its tip at $\tilde{x}$ and  $\ol{\t{x}\t{y}}\subset \t\De$.

Consider the vectors $\varpi_1,...,\varpi_r\in \t\De$ which are the
fundamental coweights of our root system. Then the vector
$\ov{\t{x}\t{y}}$ is the integer linear combination
$$
\ov{\t{x}\t{y}}= \sum_{i=1}^r n_i \varpi_i, n_i\in \N\cup \{0\}.
$$
Accordingly, we define a path $\t{p}$ in $\t\De$ with the initial
vertex $\t{x}$ and the final vertex $\t{y}$ as the concatenation
$$
\t{p}=\pi_{\ul{\la}}= \pi_{\la_1}...*\pi_{\la_r}= \t{p}_1\cup..\cup
\t{p}_r,
$$
where $\la_i:=n_i \varpi_i$, $\ul{\la}=(\la_1,...,\la_r)$. Observe
that the path  $\t{p}$ satisfies the assumptions of Theorem
\ref{chaincondition}. Moreover, $\t{p}$ is a generalized Hecke
path.

\medskip
Next, let $A\subset X$ be an apartment containing $\ol{\t{z}\t{x}}$,
$a\subset A$ be an alcove containing $\t{z}$; consider the
retraction $f=Fold_{a,A}: X\to A$. This retraction transforms
$\t{P}$ to a polygon $\hat{P}:=f(\t{P})\subset A$ which has geodesic
sides $f(\ol{\t{z}\t{x}})$, $f(\ol{\t{y}\t{z}})$. Note that the
break-points in
$$
\hat{p}:=f(\t{p})= \hat{p}_1\cup...\cup \hat{p}_r,
$$
are the images of the vertices $\t{x}_i$ of $\t{P}$ which are
break-points $\t{p}$, but in addition we possibly have break-points
within the segments $f(\t{p}_i)$. The latter can occur only at the
values of $t$ for which the geodesic segments of $\t{p}_i$ intersect
transversally the walls of $(\t{A},\t{W})$. Since $\t{x}$ is a special vertex and the edges of
each $\t{p}_i$ are parallel to multiples of $\varpi_i$, it follows that each segment $\t{p}_i$ is
contained in the 1-skeleton of $X$. Thus the break-points of $f(\t{p}_i)$ are
automatically vertices of $\t{A}$. We subdivide the path $\t{p}$ so
that all break-points of $\hat{p}$ are the images of the vertices of
$\t{P}$.

Let $k=k_R$ be the saturation constant of the root system $R$. Then,
according to Lemma \ref{kexists}, for each vertex $v\in A$, the
point $kv\in A$ is a special vertex of $A$. Therefore, applying a
dilation $h\in Dil(A,W)$ with the conformal factor $k$ to the polygon
$\hat{P}$, we obtain a new polygon $k\cdot \hat{P}=h(\hat{P})$,
whose vertices are all special vertices of $A$. Thus we can identify
the Weyl chamber $\De$ with a chamber in $A$ whose tip $o$ is at the
vertex $h(\tilde{z})$ and which contains the geodesic segment
$h(\ol{\tilde{z}\tilde{x}})$.

Lastly, let $P=\P(k \hat{P})=g(\t{P})$ denote the projection of the
polygon $k\hat{P}$ to the Weyl chamber $\De$, where
$g=Fold_{z,h,\De}$. We set $x:=g(\t{x}), y:= g(\t{y}), x_i:=
g(\t{x}_i)$, $p:= g(\t{p})$, etc.

\begin{prop}
$P=\ol{ox} \cup p \cup \ol{yo}$ is a Littelmann polygon such that
$$
\length(p)= k\cdot \length(\t{p}).
$$
\end{prop}
\proof We have to show that the path $p$ satisfies the chain
condition with maximal chains at each vertex.

The chain condition at each vertex of the path $p$ follows
immediately from Theorem \ref{chaincondition}.

The maximality condition is immediate for the break-points which
occur at the special vertices of $A$, in particular, for all
break-points which are images of the break-points of $\hat{p}$. The
remaining break-points are the ones which occur at the points $\P
(\hat{x}_i)$, where $\hat{x}_i$ are smooth points of $k\hat{p}$ at
which this path transversally intersects the  walls of $A$ passing
through $o$. However at these points the maximality condition
follows from Proposition \ref{chaincondition1}.

The second assertion of the proposition was proven in Lemma
\ref{invarianceoflength}.  \qed

\begin{rem}
\label{miniscule}
 Observe that in the case when $\be$ is the sum of
{\em minuscule fundamental coweights}, the multiplication by $k$ in the above
proof is unnecessary since all the vertices of the polygon $\hat{P}$
(and hence $P=\P(\hat{P})$) are already special. Thus, in this case,
the polygon
$$P=Fold_{\De}(\t{P})$$
is a Littelmann polygon.
\end{rem}

\medskip
The above proposition shows that the polygon $P$ is a Littelmann polygon in $\De$,
which has two geodesic sides having the $\De$-lengths $k\al, k\ga^*$ and
the concatenation of the remaining sides equal to a generalized LS path of the $\De$-length $k\be$.
Therefore, according to Littelmann's theorem (see Theorem
\ref{Lit}),
$$
 V_{k\ga}\subset V_{k\al}\otimes V_{k\be}.
\qed
$$
This concludes the proof of Theorem \ref{3->4}. \qed

\begin{cor}
\label{corQ3->Q4} Suppose that $\al,\be,\ga\in L\subset P(R^\vee)$
are dominant weights for the complex semisimple Lie group
$\ul{G}^\vee(\C)$, such that $\al+\be+\ga\in Q(R^\vee)$ and that
there exists $N\in \N$ so that
$$
 V_{N\ga}\subset V_{N\al}\otimes V_{N\be}.
$$
Then for the saturation constant $k=k^2_R$ we have
$$
V_{k\ga}\subset V_{k\al}\otimes V_{k\be} .
$$
\end{cor}
\proof Let $X$ be the Euclidean (Bruhat-Tits) building associated to
the group $\ul{G}(\K)$. Then the assumption that
$$
 V_{N\ga}\subset V_{N\al}\otimes V_{N\be}
$$
implies that $(N\al,N\be,N\ga^*)$ belongs to $D_3(X)$ (see
\cite[Theorem 9.17]{KLM3}, or \cite[Theorem 10.3]{KLM3},  or Theorem
\ref{klm3}). Since $D_3(X)$ is a homogeneous cone and $N> 0$, $(\al,\be,\ga^*)\in
D_3(X)$ as well. Moreover, according to Theorem \ref{klm2}, since
$\al, \be, \ga\in P(R^\vee)$ and $\al+\be+\ga\in Q(R^\vee)$, there
exists a triangle $T\subset X$ with the $\De$-side-lengths
$\al,\be,\ga^*$, whose vertices are also vertices of $X$. Now the
assertion follows from Part 3 of theorem \ref{3->4}. \qed

Using Remark \ref{miniscule} we also obtain:

\begin{thm}
\label{mini3->4}
Let $X$ be a building as above. Suppose that $T=[x,y,z]$ is a
geodesic triangle in $X$ with the $\De$-side-lengths
$(\al,\be,\ga^*)$, which are dominant weights of $G^\vee$ and so that
one (equivalently, all) vertices of $T$ are special and at least one
of the weights $\al,\be,\ga$ is the sum of minuscule weights. Then
$$
V_{\ga} \subset V_{\al}\otimes V_{\be}.
$$
\end{thm}

This theorem was originally proven by Tom Haines in the case when
all the weights $\al, \be, \ga$ are sums of minuscules.

\begin{conjecture}
(T. Haines) Suppose that $\al, \be, \ga$ are sums of minuscule
weights. Then, in the above theorem, the assumption that one vertex
of $T$ is special can be replaced by $\al+\be+\ga\in Q(R^\vee)$.
\end{conjecture}

Note that (among irreducible root systems) the root systems $G_2, F_4, E_8$ have no minuscule weights,
$B_n, C_n, E_7$ have exactly one minuscule weight and the root systems $A_n, D_n, E_6$ have more than 1 minuscule weights.
For the root system $A_n$ Haines conjecture follows from the saturation theorem. For $D_n$ and $E_6$ it would
follow from the affirmative answer to Question \ref{simplyconjecture}.

\begin{figure}[tbh]
\centerline{\epsfxsize=5in \epsfbox{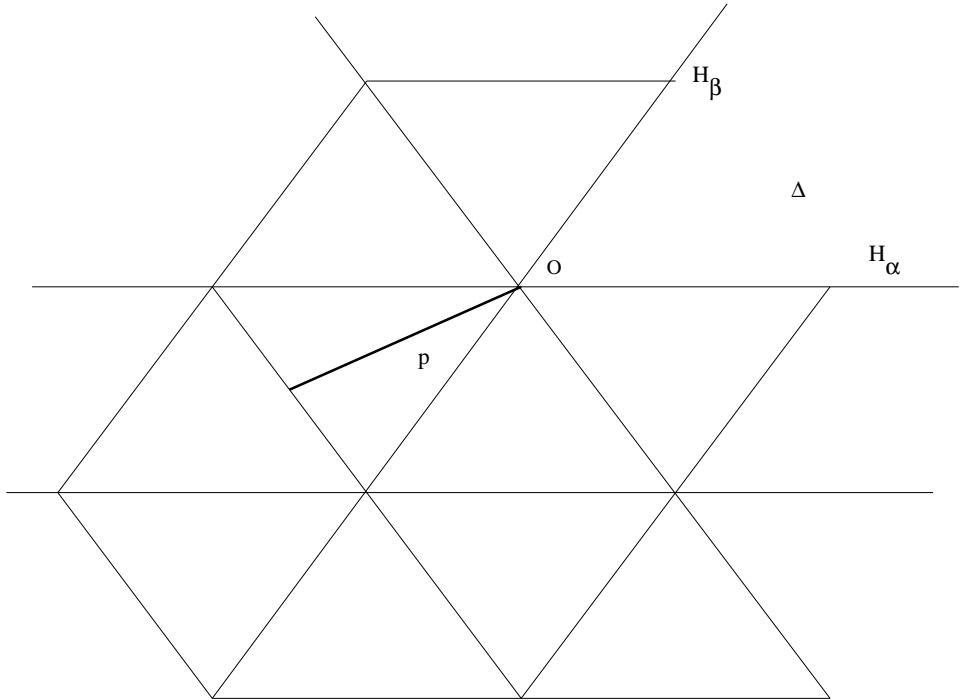}} \caption{\sl A Hecke
path which does not satisfy the integrality condition.}
\label{hecke}
\end{figure}

\begin{prop}
Suppose that the root system $R$ has exactly one minuscule coweight
$\la$. Then the above conjecture holds for $R$.
\end{prop}
\proof Let $(A,W)$ denote the Euclidean Coxeter complex
corresponding to the root system $R$ and let $X$ be a thick Euclidean building modeled on $(A,W)$.
Given $(\al,\be,\ga^*)\in D_3(X)$ such that $\al, \be, \ga\in P(R^\vee), \al+\be+\ga\in Q(R^\vee)$
we have to construct a geodesic triangle $T=[o,x,y]\subset X$ with special vertices and
the $\De$-side lengths $(\al,\be,\ga^*)$. Clearly, it suffices to treat the case when the root system $R$ is
irreducible and spans $V$. Therefore $R$ has type $B_n, C_n$, or $E_7$. In particular, the index of connection $i$ of $R$ equals $2$
and $-1\in W_{sph}$. In particular, $\ga=\ga^*$.
Let $\la$ denote the unique minuscule coweight of $R$ and let $\La$ denote the span in $\la$ in $V$.

Observe that $\la$ does not belong to the coroot lattice $Q(R^\vee)$ and thus, since $i=2$,
$$
\N\cdot \la\cap Q(R^\vee)= 2\N \cdot \la.
$$
Suppose now that $\al=a\la, \be=b\la, \ga=b\la$, where $a, b, c\in \N\cup \{0\}$ and
$$
(\al,\be,\ga)\in D_3(X), \quad \al+\be+\ga\in Q(R^\vee).
$$
Thus $a+b+c$ is an even number and the triple $(a,b,c)$ satisfies
the ordinary metric triangle inequalities.

Let $(A',W')=(\R, 2\Z\ltimes \Z/2)$ denote the rank 1 Coxeter complex;  its
vertex set equals $\Z$. The positive Weyl chamber in $(A',W')$ is $\R_+$ and
we can identify $\De'$-distances with the usual metric distances. Let $X'$ denote a thick
building  which is modeled on $(A',W')$ (i.e. a simplicial tree with edges of unit length and thickness $\ge 3$).
Then the above properties of $a, b, c$ imply that $X'$ contains a triangle $T'=[o', x', y']$ with the metric
side-lengths $a, b, c$. If this triangle is contained in a single apartment $A'\subset X'$, we send $T'$ to
a geodesic triangle $T\subset X$ via the isometry $A'\to \La\to A\to X$. If not, we
obtain a folded (Hecke) triangle $P=f'(T')=[o', x', u', f'(y')]\subset \De'$. Note that the unit tangent directions
$\xi', \eta'\in S_{u'}(A')$ to the segments $\ol{u'x'}, \ol{u' f(y')}$ are antipodal. Now embed
the apartment $A'$ into $A\subset X$ via the isometry $\iota$ that sends $A'$ to $\La$, $o'$ to $o$, $1$ to $\la$ (the latter
is a special vertex). Then the point $u:=\iota(u')$ is also a special vertex in $A$.
We claim that the resulting  broken triangle $[o, x, u, y]\subset A$ is a Hecke triangle in $A$.
Indeed, the directions $\xi, \eta$ at $\Si_u$ which are images of $\xi', \eta'$ under $\iota$ are
antipodal and $\eta\in \De$. Therefore, since $u$ is a special vertex and  $-1\in W_{sph}$, according to Lemma
\ref{constructingachain} $\xi\ge \eta$. Thus $\ol{xu}\cup \ol{uy}$ is a Hecke path. It follows that
the broken triangle $[o, x, u, y]$ is a Hecke triangle and hence it unfolds to a geodesic triangle $T$ in $X$.
The triangle $T$  has special vertices and $\De$-side lengths $(\al,\be,\ga)$.

Below is an alternative to the above argument. Let $\Si\subset \De^3$ denote the collection of
triples of dominant weights $\tau, \eta, \mu$ such that
$$
(V_\tau \otimes V_\eta \otimes V_\mu)^{G^\vee}\ne 0,
$$
where $G^\vee$ is assumed to be simply-connected.
This set is an additive semigroup, see for instance  \cite[Appendix]{KLM3}.
It suffices to prove that $(\al,\be,\ga)\in \Si$.

Set $t:= \frac{1}{2}(a+b-c)$. (The number $t$ is the metric length of the ``leg'' of the geodesic triangle
$T'\subset  X'$ in the above argument, the leg which contains the vertex $x'$.) Since $a+b+c$ is even, the number $t$
is an integer. Set
$$
a_1:= a-t, b_1:= b-t, c_1:= c, \quad \al_1:= a_1\la, \be_1:= b_1\be, \ga_1:= c_1\ga.
$$
Then $c_1=a_1+b_1$ and the metric triangle inequalities for $a, b, c$ imply that $t\ge 0, a_1\ge 0, b_1\ge 0$.
Thus $\al_1, \be_1, \ga_1$ are still dominant weights of $G^\vee$ and they satisfy
$$
\ga_1= \al_1+\be_1.
$$
Then, since $-1\in W_{sph}$ and $\ga_1=\ga_1^*$, we have: $(\al_1,\be_1,\ga_1)\in \Si$. Moreover,
$(t\la, t\la, 0)$ also clearly belongs to $\Si$ and we have
$$
(\al,\be,\ga)= (\al_1,\be_1,\ga_1) + (t\la, t\la, 0)= (\al,\be,\ga). \qed
$$

\begin{ex}
There exists a Hecke path $p\in {\mathcal P}$ such that $p(1)\in
P(R^\vee)$, however for the saturation constant $k=k_R$, the path
$k\cdot p$ is not an LS path.
\end{ex}
\proof Our example is for the root system $A_2$, in which case
$k=1$. We will give an example of a Hecke path $p\in {\mathcal P}$
such that $p(1)\in P(R^\vee)$ but $p$ does not belong to ${\mathcal
P}_\Z$. Since, according to Theorem \ref{locint}, each LS path
belong to ${\mathcal P}_\Z$, it proves that $p$ is not an LS path.

The Hecke path $p$ in question has $\De$-length $\varpi_1+\varpi_2$,
where $\varpi_1, \varpi_2$ are the fundamental coweights; the
break-point of $p$ occurs at the point $-(\varpi_1+\varpi_2)/2$
where the path $p$ backtracks back to the origin.
Thus, for the simple roots $\al$ and $\be$, the minimum of the
functions $\al(p(t)), \be(p(t))$ equals $-1/2$. See Figure
\ref{hecke}.  \qed

\bigskip
\noindent Michael Kapovich: \newline Department of Mathematics,
\newline University of California, \newline Davis, CA 95616, USA
\newline
 kapovich$@$math.ucdavis.edu

\smallskip
\noindent John J. Millson: \newline Department of Mathematics,
\newline University of Maryland, \newline College Park, MD 20742,
USA
\newline
 jjm$@$math.umd.edu


\newpage
\begin{thebibliography}{BaBE}
\addcontentsline{toc}{section}{Bibliography}




\bibitem[Ba]{Ballmann}
W.\ Ballmann, ``Lectures on spaces of nonpositive curvature. With an
appendix by Misha Brin.'' DMV Seminar, vol. 25. Birkhauser Verlag,
Basel, 1995.

\bibitem[BK]{BK}
P. Belkale, S. Kumar, {\em Eigencone, saturation and Horn problems for symplectic and odd orthogonal groups},
 Preprint arXiv:0708.0398, 2007.

\bibitem[BS]{BS}
A.\ Berenstein and R.\ Sjamaar, {\em Coadjoint orbits, moment
polytopes, and the Hilbert-Mumford criterion},  Journ. Amer. Math.
Soc., vol. {\bf 13} (2000), no. 2, p. 433--466.

\bibitem[B]{Borel}
A. Borel, ``Linear Algebraic Groups'', Springer Verlag,
Graduate Texts in Mathematics , Vol. 126.
2nd printing, 1991.


\bibitem[Bo]{Bourbaki}
N. Bourbaki, ``Lie groups and Lie algebras", Chap. 4, 5, 6. Springer Verlag, 2002.


\bibitem[Br]{Brown}
K. Brown,``Buildings'', Springer-Verlag, New York, 1989.


\bibitem[BT]{BT}
F.\ Bruhat\ and\ J. Tits, {\em Groupes reductifs sur un corps
local,} Publ. Math. IHES, No. 41, (1972), p.
5--251.

\bibitem[D]{Demazure}
M. Demazure, {\em Sch\'emas en groupes r\'eductifs},
Bull. Soc. Math. France,  {\bf 93}  (1965) p. 369--413.

\bibitem[DW]{DW}
H.\ Derksen and J.\ Weyman, {\em Semi-invariants of quivers and saturation for Littlewood-Richardson coefficients}
J. Amer. Math. Soc.  {\bf 13}  (2000),  no. 3, p. 467--479.


\bibitem[G]{Garrett}
P. Garrett, ``Buildings and classical groups,'' Chapman \& Hall, London, 1997.


\bibitem[GL]{GL}
S. Gaussent and P. Littelmann, {\em LS-Galleries, the path model and
MV-cycles},  Duke Math. J.  127  (2005),  no. 1, p. 35--88.


\bibitem[Gro]{Gross}
B.\ Gross, {\em On the Satake isomorphism}, In: ``Galois
representations in arithmetic algebraic geometry
              (Durham, 1996)'',
London Math. Soc. Lecture Notes, vol. {\bf 254}, (1998) p. 223--237.


\bibitem[Ha]{Haines}
T.\ J.\ Haines, {\em Structure constants for Hecke and
representations rings}, IMRN vol. {\bf 39} (2003), p. 2103--2119.


\bibitem[KKM]{KKM}
M. Kapovich, S. Kumar and J. J. Millson, {\em Saturation and irredundancy for $Spin(8)$},
To appear in  Pure and Applied Mathematics Quarterly.

\bibitem[KLM1]{KLM1}
M. Kapovich, B.\ Leeb and J.\ J.\ Millson, {\em Convex functions on
symmetric spaces, side lengths of polygons and the stability
inequalities for weighted configurations at infinity}, Preprint,
June 2004.


\bibitem[KLM2]{KLM2}
M. Kapovich, B.\ Leeb and J.\ J.\ Millson,  {\em Polygons in buildings and
their side-lengths}, Preprint, 2004.


\bibitem[KLM3]{KLM3}
M. Kapovich, B.\ Leeb and J.\ J.\ Millson,  {\em Polygons in symmetric
spaces and buildings with applications to
  algebra},  Memoirs of AMS, to appear.


\bibitem[KL]{KleinerLeeb}
B.\ Kleiner and B.\ Leeb, {\em Rigidity of quasi-isometries for
symmetric spaces and Euclidean buildings}, Publ. Math. IHES, vol.
{\bf 86} (1997), p. 115--197.




\bibitem[KT]{KnutsonTao}
A. Knutson\ and\ T. Tao, {\em The honeycomb model of ${\rm
GL}_n(\C)$ tensor products. I. Proof of the saturation conjecture,}
J. Amer. Math. Soc., vol. {\bf 12} (1999), no.~4, p. 1055--1090.




\bibitem[L1]{Littelmann1}
P.\ Littelmann, {\em  A Littlewood-Richardson rule for symmetrizable
Kac-Moody  algebras},
Invent. Math. 116 (1994), no. 1-3, p. 329--346.


\bibitem[L2]{Littelmann2}
P.\ Littelmann, {\em Paths and root operators in representation
theory}, Annals of Math. (2)  142 (1995) no. 3, p. 499--525.


\bibitem[L3]{Littelmann3}
P.\ Littelmann, {\em Characters of representations and paths in
${\mathfrak H}^*_{\R}$}, In: ``Representation theory and automorphic
forms'' (Edinburgh, 1996), p. 29--49, Proc. Sympos. Pure Math., 61,
Amer. Math. Soc., Providence, RI, 1997.


\bibitem[Ron]{Ronan}
M. Ronan,``Lectures on buildings'', Perspectives in Mathematics, 7.
Academic Press, Inc., 1989.


\bibitem[Rou]{Rousseau}
G. Rousseau, {\em Euclidean buildings}, Lectures at \'Ecole d'\'et\'e de
math\'ematiques ``Nonpositively curved geometries, discrete groups and
rigidities'', Institut Fourier, Grenoble, 2004.

\bibitem[Sat]{Satake}
I.\ Satake, {\em Theory of spherical functions on reductive algebraic groups
over $p$-adic fields}, Publ. Math. IHES, vol. {\bf 18} (1963), p. 1--69.

\bibitem[Sc]{Schwer}
C. Schwer, {\em Galleries, Littlewood polynomials and
structure constants of the spherical Hecke algebra},
Int. Math. Res. Not.  2006, Art. ID 75395, 31 pp.


\end{thebibliography}
\end{document}